\documentclass[journal]{new-aiaa}

\usepackage[utf8]{inputenc}
\usepackage{newtxtext,newtxmath}

\usepackage{multirow}
\usepackage{color}
\usepackage{indentfirst}
\usepackage{graphicx}

\usepackage{amsmath,amsfonts,amssymb,amsthm,amsbsy}
\usepackage{enumitem}
\usepackage[ruled,linesnumbered]{algorithm2e}
\usepackage{subcaption}

\theoremstyle{plain}

\theoremstyle{definition}

\theoremstyle{remark}
\newcounter{REM}
\newtheorem{rem}[REM]{Remark}

\newcommand{\btheta}{\boldsymbol{\theta}}
\newcommand{\bpi}{\boldsymbol{\pi}}
\newcommand{\bmu}{\boldsymbol{\mu}}
\newcommand{\bPhi}{\boldsymbol{\Phi}}
\newcommand{\bPsi}{\boldsymbol{\Psi}}

\allowdisplaybreaks

\title{Incremental Correction in Dynamic Systems Modelled with Neural Networks for Constraint Satisfaction}

\author{Namhoon Cho\footnote{Research Fellow, Centre for Autonomous and Cyber-Physical Systems, School of Aerospace, Transport and Manufacturing, \texttt{n.cho@cranfield.ac.uk}, AIAA Member}, 
Hyo-Sang Shin\footnote{Professor of Guidance, Navigation and Control, Centre for Autonomous and Cyber-Physical Systems, School of Aerospace, Transport and Manufacturing, \texttt{h.shin@cranfield.ac.uk}}, and 
Antonios Tsourdos\footnote{Professor, Centre for Autonomous and Cyber-Physical Systems, School of Aerospace, Transport and Manufacturing, \texttt{a.tsourdos@cranfield.ac.uk}}}
\affil{Cranfield University, Cranfield, MK43 0AL, Bedfordshire, United Kingdom}
\author{Davide Amato\footnote{Lecturer in Spacecraft Engineering, Department of Aeronautics, Faculty of Engineering, \texttt{d.amato@imperial.ac.uk}}}
\affil{Imperial College London, South Kensington Campus, SW7 2AZ, London, United Kingdom}

\begin{document}
	\maketitle
	
	\begin{abstract}
		This study presents incremental correction methods for refining neural network parameters or control functions entering into a continuous-time dynamic system to achieve improved solution accuracy in satisfying the interim point constraints placed on the performance output variables. The proposed approach is to linearise the dynamics around the baseline values of its arguments, and then to solve for the corrective input required to transfer the perturbed trajectory to precisely known or desired values at specific time points, i.e., the interim points. Depending on the type of decision variables to adjust, parameter correction and control function correction methods are developed. These incremental correction methods can be utilised as a means to compensate for the prediction errors of pre-trained neural networks in real-time applications where high accuracy of the prediction of dynamical systems at prescribed time points is imperative. In this regard, the online update approach can be useful for enhancing overall targeting accuracy of finite-horizon control subject to point constraints using a neural policy. Numerical example demonstrates the effectiveness of the proposed approach in an application to a powered descent problem at Mars.
	\end{abstract}
	
	\section{Introduction}
Recent advances in machine learning theory and methodology have led to many successful applications of deep neural networks (NN) in analysis or control of dynamic systems \cite{Rackauckas_2021}. As a parametric function approximator, NN can be incorporated as an element in the dynamic system represented with a set of Ordinary Differential Equations (ODEs) or even the entire ODE function \cite{Kim_2021}. In model learning problems, the NNs are fitted to the measurement data to establish an ODE function that yields more accurate solution when propagated forward in time. In control problems, the control law given by NN is trained to provide optimal control performance \cite{Cho_2022}.

Typical NN training process may fail to satisfy state constraints within the predicted ODE solution. Regardless of what the NN component represents in the ODE model considered in each application, the usual practice is to train the parameters by leveraging unconstrained optimisation of a cost functional which is usually integrated over time or summed over multiple instances \cite{Rackauckas_2019,Cho_2022}. Unconstrained optimisers including the first-order gradient-descent-based algorithms such as ADAM \cite{Kingma_2015} and NADAM \cite{Dozat_2016} are widely used. In this setup, state constraints are incorporated in a soft form as in the penalty method. As a consequence, even the final optimised parameters may produce ODE solutions that exhibit only a limited accuracy in state constraint satisfaction. However, in problems where the state trajectory is known or desired to attain particular values at certain time points, it is desirable to determine the NN parameters or other adjustable elements in the dynamic system to satisfy the state equality constraints accurately with the solution predicted from the ODE model embedding NN. Incorporating available information given in the form of equality constraints at certain time points can help improving prediction accuracy, as in constrained Kalman filters \cite{Simon_2002,Simon_2010}.

On the other hand, the NN parameters optimised offline using a nominal dynamic model (and assumed uncertainty models) cannot always provide very accurate online prediction due to modelling uncertainties. In finite-horizon control problems, discrepancy between the actual environment and the dynamic model assumed for policy synthesis results in degraded performance of an offline-learned NN feedback policy in achieving the desired final state. Hence, some form of online update is necessary to improve the final state targeting accuracy, however, the high computation bottleneck prohibits a large-sized NN from being trained online in real-time.

To overcome these difficulties, this study presents incremental correction methods to modify the dynamic systems modelled with ODEs embedding pre-determined NNs to improve the accuracy of state constraint satisfaction at given interim points. In principle, constrained optimisation can be performed at the stage of NN weight training. However, it will add complexity to the training process in finding feasible solutions since the number of decision variables (NN parameters) is quite large, and the decision variables affect the state variables only indirectly; in addition, this might result in a limited exploration of the NN parameter space. Instead, this study leaves the NN training pipeline unaltered and takes inspirations from the philosophy of neighbouring optimal control and gradient methods for trajectory optimisation \cite{Bryson_1962,Bryson_1963,Denham_1964}. The proposed approach consists in local linearisation of the system dynamics around the baseline input values followed by the design of an augmentation for the input to enforce state constraints. The correction algorithms are developed in both the NN parameter space and the control function space. In particular, the online-generated correction of parameter or control function on top of the offline-trained NN baseline policy will realise a hybrid offline-online paradigm taking the benefits of both worlds: i) the strength of policy search in learning sophisticated policies by leveraging simulation and data, and ii) the robustness of online correction in providing real-time optimal feedback performance. In summary, the proposed methods can be used as a post-processing step to achieve the state constraints without the necessity to change the pre-existing NN training procedure.

The rest of the paper is organised as follows. Section \ref{Sec:Overview} provides a brief overview of the proposed approach consisting of baseline policy training followed by incremental correction. Section \ref{Sec:IncrCorr} presents the derivation of the incremental correction methods; i) parameter correction method in Sec. \ref{SubSec:ParamCorr}, and ii) control function correction method in Sec. \ref{SubSec:CtrlFuncCorr}. Section \ref{Sec:Applications} demonstrates the efficacy of the incremental correction approach through illustrative examples. Concluding remarks are summarised in Sec. \ref{Sec:Concls}.

\section{Overview} \label{Sec:Overview}
This study considers the problem of designing a NN-based controller for finite-horizon control of a continuous-time nonlinear dynamic system subject to the interim point constraints imposed on the performance output which is linear in the state. Applying constrained parameter optimisation methods for finding the NN policy parameters cannot straightforwardly solve the problem, since the constraints are imposed on the solution of the dynamic system rather than on the NN parameters that are embedded inside the differential equations. Also, random initialisation of the large number of NN parameters is likely to pose difficulties in convergence to a feasible solution with constrained trajectory optimisation methods which require initial guess for the control input. As a workaround, a widely-used methodology is to penalise constraints violation in the cost definition and exploit the NN training tools suitable for unconstrained problems. However, it was observed in numerical experiments that the achieved targeting accuracy tends to leave room for improvement even when very large penalty weights are used.

\begin{figure}[ht!]
	\begin{center}
		\includegraphics[width = 0.25\textwidth]{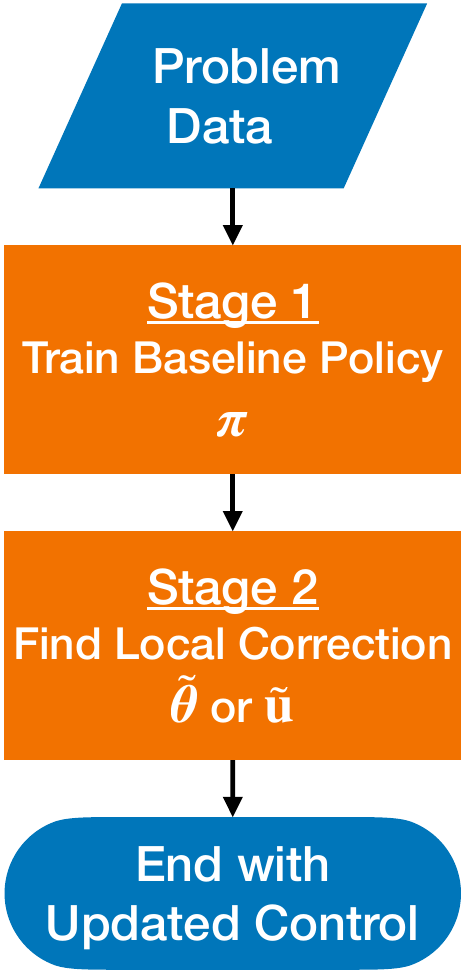}
		\caption{Overview of Proposed Approach} \label{Fig_Approach_Overview}
	\end{center}
\end{figure}

This study presents a two-stage approach to improve the accuracy of satisfying the constraints while circumventing the difficulties in considering hard constraints in NN policy training. Figure \ref{Fig_Approach_Overview} shows the schematic diagram of the proposed two-stage approach. Stage 1 is to train a baseline NN policy with the penalty-based formulation in the same way as it has been done usually. The continuous-time policy gradient method based on adjoint sensitivity previously proposed by some of the authors in Ref. \cite{Cho_2022} can serve this purpose. Stage 2 finds the incremental correction of either the NN parameters or the control input needed to satisfy the interim point constraints with the closed-loop dynamics linearised around the baseline trajectory. The following section presents two different methods for Stage 2.

\section{Incremental Correction Methods} \label{Sec:IncrCorr}
\subsection{Parameter Correction} \label{SubSec:ParamCorr}
Consider the continuous-time system dynamics given by
\begin{equation} \label{Eq:SysDyn_ParamCorr}
	\begin{aligned}
	\dot{\mathbf{x}}\left(t\right) &= \mathbf{f}_{\theta}\left(t, \mathbf{x}\left(t\right), \btheta \right), \quad \mathbf{x}\left(t_{0}\right)=\mathbf{x}_{0}\\
	\mathbf{z}\left(t\right) &= \mathbf{H}\mathbf{x}\left(t\right)
	\end{aligned}
\end{equation}
where $t$, $\mathbf{x}\in\mathbb{R}^{n\times 1}$, $\mathbf{z}\in\mathbb{R}^{p\times 1}$, and $\btheta \in \mathbb{R}^{l \times 1}$ denote the time, the state, the performance output, and the parameter vector, respectively, and $\mathbf{f}_{\theta}$ represents the ODE function with $\btheta$ as its decision variable. Usually, not every state variable needs to be constrained. The performance output is defined as a linear combination of state variables that is subjected to given constraints, and $\mathbf{H}$ is the constant matrix that maps the state to the performance output. The overdot notation stands for the time-derivative of a quantity. Note that a NN embedded in the ODE function is implicit in Eq. \eqref{Eq:SysDyn_ParamCorr}, and its weights are embedded in the parameter vector $\btheta$.

Suppose that the NN is already trained by minimising a loss function and the optimised parameter vector $\btheta^{*}$ is given as the result of Stage 1 training. Let $\mathbf{x}^{*}\left(t\right)$ denote the solution predicted by integrating the ODE model with $\btheta^{*}$ and a given initial condition $\mathbf{x}_{0}^{*}$. That is, 
\begin{equation} \label{Eq:SysDyn_ParamCorr_Nominal}
	\begin{aligned}
	\dot{\mathbf{x}}^{*}\left(t\right) &= \mathbf{f}_{\theta}\left(t, \mathbf{x}^{*}\left(t\right), \btheta^{*} \right), \quad \mathbf{x}^{*}\left(t_{0}\right)=\mathbf{x}_{0}^{*}\\
	\mathbf{z}^{*}\left(t\right) &= \mathbf{H}\mathbf{x}^{*}\left(t\right)
	\end{aligned}
\end{equation}
for $\forall t \in \left[t_{0},t_{f}\right]$. Consider a small perturbation in both the state and the parameter that can be expressed as
\begin{equation} \label{Eq:Perturb_ParamCorr}
	\begin{aligned}
		\mathbf{x}\left(t\right) &= \mathbf{x}^{*}\left(t\right) + \tilde{\mathbf{x}}\left(t\right)\\ 
		\btheta &= \btheta^{*} + \tilde{\btheta}
	\end{aligned}
\end{equation}
where the tilde notation refers to the perturbed terms. Linearising Eq. \eqref{Eq:SysDyn_ParamCorr} around the baseline prediction $\mathbf{x}^{*}\left(t\right)$ and the baseline parameter vector $\btheta^{*}$ yields
\begin{equation} \label{Eq:LinSysDyn_ParamCorr}
	\begin{aligned}
	\dot{\tilde{\mathbf{x}}} \left(t\right) &\approxeq \left. \frac{\partial \mathbf{f}_{\theta}}{\partial \mathbf{x}} \right|_{\mathbf{x}^{*}\left(t\right),\btheta^{*}} \tilde{\mathbf{x}}\left(t\right) + \left. \frac{\partial \mathbf{f}_{\theta}}{\partial \btheta} \right|_{\mathbf{x}^{*}\left(t\right),\btheta^{*}} \tilde{\btheta}\\
	&:=\mathbf{A}_{\theta}\left(t\right)\tilde{\mathbf{x}}\left(t\right) + \mathbf{B}_{\theta}\left(t\right)\tilde{\btheta}
	\end{aligned}
\end{equation}
where the Jacobian matrices $\mathbf{A}_{\theta}\left(t\right)$ and $\mathbf{B}_{\theta}\left(t\right)$ are defined accordingly. 
The perturbed dynamic system in Eq. \eqref{Eq:LinSysDyn_ParamCorr} is generally a linear time-varying (LTV) system whose solution can be written as
\begin{equation} \label{Eq:sol_LinSysDyn_ParamCorr}
	\tilde{\mathbf{x}} \left(t\right) = \bPhi\left(t,t_{0}\right)\tilde{\mathbf{x}} \left(t_{0}\right) + \int_{t_{0}}^{t} \bPhi\left(t,\tau\right)\mathbf{B}_{\theta}\left(\tau\right)\tilde{\btheta}d\tau
\end{equation}
where $\bPhi\left(t_{2},t_{1}\right)$ is the state transition matrix for transfer from $t_{1}$ to $t_{2}$ defined by
\begin{equation} \label{Eq:Phi_transit}
	\dot{\bPhi}\left(t,t_{0}\right) = \mathbf{A}_{\theta}\left(t\right)\bPhi\left(t,t_{0}\right), \quad \bPhi\left(t_{0},t_{0}\right) = \mathbf{I}
\end{equation}
Since $\btheta$ is a constant perturbation of parameters, Eq. \eqref{Eq:sol_LinSysDyn_ParamCorr} can be rewritten as
\begin{equation} \label{Eq:sol_LinSysDyn_ParamCorr_LinTrans}
	\tilde{\mathbf{x}} \left(t\right) = \bPhi\left(t,t_{0}\right)\tilde{\mathbf{x}} \left(t_{0}\right) + \mathbf{M}\left(t\right) \tilde{\btheta}
\end{equation}
with $\mathbf{M}\left(t\right) := \int_{t_{0}}^{t}\bPhi\left(t,\tau\right)\mathbf{B}_{\theta}\left(\tau\right)d\tau$. With the knowledge of $\mathbf{A}_{\theta}\left(t\right)$ and $\mathbf{B}_{\theta}\left(t\right)$, the matrix $\mathbf{M}\left(t\right)$ can be evaluated without explicitly computing the state transition function but instead by integrating the following differential equation once from $t_{0}$ to $t_{f}$:
\begin{equation} \label{Eq:M_Dyn}
	\dot{\mathbf{M}}\left(t\right) = \mathbf{A}_{\theta}\left(t\right)\mathbf{M}\left(t\right) + \mathbf{B}_{\theta}\left(t\right), \quad \mathbf{M}\left(t_{0}\right) = \mathbf{0}
\end{equation}
It is obvious from Eq. \eqref{Eq:sol_LinSysDyn_ParamCorr_LinTrans} that the state perturbation at a time point $t$ is related linearly with the parameter perturbation $\tilde{\btheta}$. Also note from Eq. \eqref{Eq:sol_LinSysDyn_ParamCorr_LinTrans} that the matrix $\mathbf{M}\left(t\right)$ is the first-order sensitivity of the state solution with respect to the perturbation in parameters $\frac{\partial \mathbf{x}\left(t\right)}{\partial \boldsymbol{\theta}}$ at time $t$.

Let the interim point constraints imposed on the performance output variables at some time points $t_{i}$ be given by
\begin{equation} \label{Eq:z_checkpts}
	\mathbf{z}\left(t_{i}\right) = \mathbf{z}_{i} \qquad \left(i = 1, \cdots, N\right).
\end{equation}
The baseline trajectory does not generally satisfy the interim point constraints. The main concept of the proposed approach is to find a parameter correction $\tilde{\btheta}$ such that the updated parameter $\btheta$ will bring the updated prediction $\mathbf{z}\left(t\right)$ on $\mathbf{z}_{i}$ at each $t_{i}$. In other words, $\tilde{\mathbf{z}}\left(t_{i}\right) = \mathbf{H}\tilde{\mathbf{x}}\left(t_{i}\right) = \mathbf{z}_{i} - \mathbf{z}^{*}\left(t_{i}\right)$ should be satisfied by the design of $\tilde{\btheta}$. In view of Eq. \eqref{Eq:sol_LinSysDyn_ParamCorr_LinTrans}, the required correction can simply be obtained by concatenation of constraint relations followed by the Moore-Penrose generalised inverse, i.e.,
\begin{equation} \label{Eq:tilde_btheta_opt}
	\tilde{\btheta} = \mathbf{L}^{\dagger} \left(\mathbf{Z}_{h} - \mathbf{F}_{h}\right)
\end{equation}
where 
\begin{equation} \label{Eq:tilde_btheta_opt_aux}
	\begin{aligned}
		\mathbf{L} 	& := \mathbf{H}\begin{bmatrix} \mathbf{M}\left(t_{1}\right) & \cdots & \mathbf{M}\left(t_{N}\right)\end{bmatrix}\\
		\mathbf{Z}_{h} 		& := \begin{bmatrix} \mathbf{z}_{1} - \mathbf{z}^{*}\left(t_{1}\right) & \cdots & \mathbf{z}_{N} - \mathbf{z}^{*}\left(t_{N}\right) \end{bmatrix}\\
		\mathbf{F}_{h} 		& := \mathbf{H}\begin{bmatrix} \bPhi\left(t_{1},t_{0}\right)\tilde{\mathbf{x}}\left(t_{0}\right) & \cdots & \bPhi\left(t_{N},t_{0}\right)\tilde{\mathbf{x}}\left(t_{0}\right) \end{bmatrix}
	\end{aligned}
\end{equation}
and $\dagger$ refers to the pseudoinverse. 

\begin{rem} \label{Rem:MinNorm_Update}
	In many cases of employing a deep NN, usually $l$ far exceeds $n$. Thus, an underdetermined system of linear equations needs to be solved for the parameter correction. In this regard, it is sensible to utilise the Moore-Penrose generalised inverse as it produces the minimum-Frobenius-norm solution.
\end{rem}

\begin{rem} \label{Rem:Jacobian_calculation}
	In modern scientific computing environments, the Jacobians of the ODE function involving a NN can be computed by using symbolic or automatic differentiation tools. In this study, a source-to-source backward mode automatic differentiation framework called \texttt{Zygote.jl} is utilised in a \texttt{Julia}-based implementation \cite{Innes_2018}.
\end{rem}

\begin{rem} \label{Rem:AD_solver}
	Parameter correction can also be performed by exploiting local forward sensitivity analysis, which directly gives the Jacobian of state solution with respect to parameters along time \cite{Ma_2018}. In modern differentiable programming environment, automatic differentiation through an ODE solver can simply be employed for the purpose of obtaining the sensitivity matrix $\left.\frac{\partial \mathbf{x}^{*}\left(t_{i}\right)}{\partial \boldsymbol{\theta}}\right|_{\boldsymbol{\theta}^{*}}$ instead of directly constructing the dynamics Jacobian matrices $\mathbf{A}_{\theta}\left(t\right)$ and $\mathbf{B}_{\theta}\left(t\right)$ then solving Eq. \eqref{Eq:M_Dyn} for $\mathbf{M}\left(t_{i}\right)$.
\end{rem}

\subsection{Control Function Correction} \label{SubSec:CtrlFuncCorr}
Consider the continuous-time system dynamics given by
\begin{equation} \label{Eq:SysDyn_CtrlFuncCorr}
	\begin{aligned}
		\dot{\mathbf{x}}\left(t\right) &= \mathbf{f}_{u}\left(t, \mathbf{x}\left(t\right), \mathbf{u}\left(t\right) \right), \quad \mathbf{x}\left(t_{0}\right)=\mathbf{x}_{0}\\
		\mathbf{u}\left(t\right) &= \bpi\left(t, \mathbf{x}\left(t\right), \btheta\right)\\
		\mathbf{z}\left(t\right) &= \mathbf{H}\mathbf{x}\left(t\right)
	\end{aligned}
\end{equation}
where $t$, $\mathbf{x}\in\mathbb{R}^{n\times 1}$, $\mathbf{u}\in\mathbb{R}^{m\times 1}$, $\mathbf{z}\in\mathbb{R}^{p \times 1}$, and $\btheta \in \mathbb{R}^{l \times 1}$ denote the time, the state, the control input, the performance output, and the parameter, respectively. In Eq. \eqref{Eq:SysDyn_CtrlFuncCorr}, $\mathbf{f}_{u}$ represents the ODE function with $\mathbf{u}$ as its decision variable, $\boldsymbol{\pi}$ represents the NN policy, and $\mathbf{H}$ is the constant performance output matrix. Note that the system description in Eq. \eqref{Eq:SysDyn_CtrlFuncCorr} can be reduced to the form in Eq. \eqref{Eq:SysDyn_ParamCorr} by defining $\mathbf{f}_{\theta} \left(t,\mathbf{x},\btheta\right):=\mathbf{f}_{u}\left(t,\mathbf{x},\bpi\left(t,\mathbf{x},\btheta\right)\right)$. 

One may consider perturbing $\mathbf{f}_{u}$ to find a control function correction instead of performing parameter correction through the linearisation of the parameter-embedded form ODE $\mathbf{f}_{\theta}$. Suppose that the optimised parameter vector $\btheta^{*}$ is given a priori as the result of NN training in Stage 1 based on unconstrained optimisation. Let $\mathbf{x}^{*}\left(t\right)$ and $\mathbf{u}^{*}\left(t\right)$ denote the baseline state and control input, respectively, that are obtained by forward propagation of Eq. \eqref{Eq:SysDyn_CtrlFuncCorr} with given $\btheta^{*}$ and $\mathbf{x}_{0}^{*}$. Then, the following relation holds
\begin{equation} \label{Eq:SysDyn_CtrlFuncCorr_Nominal}
	\begin{aligned}
		\dot{\mathbf{x}}^{*}\left(t\right) &= \mathbf{f}_{u}\left(t, \mathbf{x}^{*}\left(t\right), \mathbf{u}^{*}\left(t\right) \right), \quad \mathbf{x}^{*}\left(t_{0}\right)=\mathbf{x}_{0}^{*}\\
		\mathbf{u}^{*}\left(t\right) &= \bpi\left(t, \mathbf{x}^{*}\left(t\right), \btheta^{*}\right)\\
		\mathbf{z}^{*}\left(t\right) &= \mathbf{H}\mathbf{x}^{*}\left(t\right)
	\end{aligned}
\end{equation}
for $\forall t \in \left[t_{0},t_{f}\right]$. Now, a small perturbation in both the state and the control input defined by
\begin{equation} \label{Eq:Perturb_CtrlFuncCorr}
	\begin{aligned}
		\mathbf{x}\left(t\right) &= \mathbf{x}^{*}\left(t\right) + \tilde{\mathbf{x}}\left(t\right)\\
		\mathbf{u}\left(t\right) &= \mathbf{u}^{*}\left(t\right) + \tilde{\mathbf{u}}\left(t\right)
	\end{aligned}
\end{equation}
leads to the linearisation of Eq. \eqref{Eq:SysDyn_CtrlFuncCorr} around the baseline state $\mathbf{x}^{*}\left(t\right)$ and the baseline control input $\mathbf{u}^{*}\left(t\right)$ as
\begin{equation} \label{Eq:LinSysDyn_CtrlFuncCorr}
	\begin{aligned}
	\dot{\tilde{\mathbf{x}}} \left(t\right) &\approxeq \left. \frac{\partial \mathbf{f}_{u}}{\partial \mathbf{x}} \right|_{\mathbf{x}^{*}\left(t\right),\mathbf{u}^{*}\left(t\right)} \tilde{\mathbf{x}}\left(t\right) + \left. \frac{\partial \mathbf{f}_{u}}{\partial \mathbf{u}} \right|_{\mathbf{x}^{*}\left(t\right),\mathbf{u}^{*}\left(t\right)} \tilde{\mathbf{u}}\left(t\right)\\
	&:=\mathbf{A}_{u}\left(t\right)\tilde{\mathbf{x}}\left(t\right) + \mathbf{B}_{u}\left(t\right)\tilde{\mathbf{u}}\left(t\right)
	\end{aligned}
\end{equation}
The Jacobian matrices $\mathbf{A}_{u}\left(t\right)$ and $\mathbf{B}_{u}\left(t\right)$ are defined appropriately according to Eq. \eqref{Eq:LinSysDyn_CtrlFuncCorr}. The solution of the perturbed dynamic system in Eq. \eqref{Eq:LinSysDyn_CtrlFuncCorr} is related to the control function correction as
\begin{equation} \label{Eq:sol_LinSysDyn_CtrlFuncCorr}
	\tilde{\mathbf{x}} \left(t\right) = \bPhi\left(t,t_{0}\right)\tilde{\mathbf{x}} \left(t_{0}\right) + \int_{t_{0}}^{t} \bPhi\left(t,\tau\right)\mathbf{B}_{u}\left(\tau\right)\tilde{\mathbf{u}}\left(\tau\right)d\tau
\end{equation}
where $\bPhi\left(t_{2},t_{1}\right)$ represents the state transition matrix associated with $\mathbf{A}_{u}\left(t\right)$.

The proposed approach is to find a minimal amount of control function correction to satisfy the interim point constraints by solving the following function space optimisation problem:
\begin{equation} \label{Eq:Prob_CtrlFuncCorr}
	\begin{aligned}
		&\text{minimise}& J &= \frac{1}{2}\int_{t_{0}}^{t_{f}}\tilde{\mathbf{u}}^{T}\left(\tau\right)\mathbf{R}\left(\tau\right)\tilde{\mathbf{u}}\left(\tau\right) d\tau\\
		&\text{subject to}& \tilde{\mathbf{z}} \left(t_{i}\right) &=  \mathbf{H}\bPhi\left(t_{i},t_{0}\right)\tilde{\mathbf{x}} \left(t_{0}\right)  + \mathbf{H}\int_{t_{0}}^{t_{i}} \bPhi\left(t_{i},\tau\right) \mathbf{B}_{u}\left(\tau\right)\tilde{\mathbf{u}}\left(\tau\right)d\tau\\
		&&  &= \mathbf{z}_{i} - \mathbf{z}^{*}\left(t_{i}\right) \qquad \left(i=1,\cdots,N\right)
	\end{aligned}
\end{equation}
where $\mathbf{R}\left(t\right)=\mathbf{R}^{T}\left(t\right) \succ 0$ is a weighting function. The interim point indices are sorted in the order of increasing time, i.e., $t_{0} \leq t_{1} < \cdots < t_{N} \leq t_{f}$. The problem can be solved by using various linear control approaches. Applying the method of constraint-coupling Lagrange multipliers, one can define the augmented cost function as
\begin{equation} \label{Eq:J_aug}
	\begin{aligned}
	J_{a} &:= \frac{1}{2}\int_{t_{0}}^{t_{f}}\tilde{\mathbf{u}}^{T}\left(\tau\right)\mathbf{R}\left(\tau\right)\tilde{\mathbf{u}}\left(\tau\right) d\tau\\  
	&\quad + \sum_{i=1}^{N}{\bmu_{i}}^{T}\left(\mathbf{z}_{i} - \mathbf{z}^{*}\left(t_{i}\right) - \mathbf{H}\bPhi\left(t_{i},t_{0}\right)\tilde{\mathbf{x}} \left(t_{0}\right) \vphantom{\int} 
	- \mathbf{H}\int_{t_{0}}^{t_{i}} \bPhi\left(t_{i},\tau\right)\mathbf{B}_{u}\left(\tau\right)\tilde{\mathbf{u}}\left(\tau\right)d\tau\right)
	\end{aligned} 
\end{equation}
where $\bmu_{i} \in \mathbb{R}^{p \times 1}$ represents the Lagrange multiplier vector associated with the $i$-th performance output constraint. The upper limit of the integral for each interim point constraint can be lifted up to $t_{f}$ by rewriting Eq. \eqref{Eq:J_aug} as
\begin{equation} \label{Eq:J_aug_activ}
	\begin{aligned}
	J_{a}&:=\sum_{i=1}^{N}{\bmu_{i}}^{T}\left(\mathbf{z}_{i} - \mathbf{z}^{*}\left(t_{i}\right) - \mathbf{H}\bPhi\left(t_{i},t_{0}\right)\tilde{\mathbf{x}} \left(t_{0}\right) \right)\\
	&\quad + \int_{t_{0}}^{t_{f}} \left[\frac{1}{2}\tilde{\mathbf{u}}^{T}\left(\tau\right)\mathbf{R}\left(\tau\right)\tilde{\mathbf{u}}\left(\tau\right) 
	-  \sum_{i=1}^{N}{\bmu_{i}}^{T} 1\left(\tau \leq t_{i}\right)\mathbf{H}\bPhi\left(t_{i},\tau\right)\mathbf{B}_{u}\left(\tau\right)\tilde{\mathbf{u}}\left(\tau\right)\right]d\tau
	\end{aligned}
\end{equation}
with the activator function defined by
\begin{equation} \label{Eq:activ}
	1\left(t\leq t_{i}\right) := 
	\begin{cases}
		1 & \text{if~} t\leq t_{i}\\
		0 & \text{if~} t > t_{i}
	\end{cases}
\end{equation}
The first variation of the augmented cost function $\delta J_{a}$ vanishes at the optimal solution. This leads to the set of equations describing the first-order necessary condition for optimality. More specifically, the first necessary condition requires $\frac{\partial \left(\text{integrand of }J_{a}\right)}{\partial \left(\tilde{\mathbf{u}}\left(\tau\right)\right)} = \mathbf{0}$ to hold, which results in
\begin{equation} \label{Eq:FONC_u}
	\begin{aligned}
		\tilde{\mathbf{u}}\left(\tau\right) &= \mathbf{R}^{-1}\left(\tau\right){\mathbf{B}_{u}}^{T}\left(\tau\right)\sum_{i=1}^{N}\bPhi^{T}\left(t_{i},\tau\right)\mathbf{H}^{T}1\left(\tau\leq t_{i}\right)\bmu_{i}\\
		&= \mathbf{R}^{-1}\left(\tau\right){\mathbf{B}_{u}}^{T}\left(\tau\right)\mathbf{E}\left(\tau\right)\bar{\bmu}
	\end{aligned}
\end{equation}
where 
\begin{equation} \label{Eq:E}
	\mathbf{E}\left(\tau\right) := \begin{bmatrix}\bPhi^{T}\left(t_{1},\tau\right)\mathbf{H}^{T}1\left(\tau\leq t_{1}\right) & \cdots & \bPhi^{T}\left(t_{N},\tau\right)\mathbf{H}^{T}1\left(\tau\leq t_{N}\right)\end{bmatrix}
\end{equation}
\begin{equation} \label{Eq:bar_bmu}
	\bar{\bmu} := \begin{bmatrix} {\bmu_{1}}^{T} & \cdots & {\bmu_{N}}^{T} \end{bmatrix}^{T}
\end{equation}
As the second necessary condition, $\frac{\partial J_{a}}{\partial \bmu_{i}} = \mathbf{0}$ should be satisfied $\forall\, i = 1, \cdots N$, which can be rewritten as
\begin{equation} \label{Eq:FONC_mu}
	\mathbf{H}\int_{t_{0}}^{t_{f}} 1\left(\tau \leq t_{i}\right)\bPhi\left(t_{i},\tau\right)\mathbf{B}_{u}\left(\tau\right)\tilde{\mathbf{u}}\left(\tau\right)d\tau = \mathbf{z}_{i} - \mathbf{z}^{*}\left(t_{i}\right) - \mathbf{H}\bPhi\left(t_{i},t_{0}\right)\tilde{\mathbf{x}} \left(t_{0}\right)
\end{equation}
for $i = 1, \cdots, N$. Substituting Eq. \eqref{Eq:FONC_u} into Eq. \eqref{Eq:FONC_mu} gives
\begin{equation} \label{Eq:constraints_u_opt_subs}
	\begin{aligned}
		&\mathbf{z}_{i} - \mathbf{z}^{*}\left(t_{i}\right)- \mathbf{H}\bPhi\left(t_{i},t_{0}\right)\tilde{\mathbf{x}} \left(t_{0}\right)\\
		&=\int_{t_{0}}^{t_{f}} 1\left(\tau \leq t_{i}\right)\mathbf{H}\bPhi\left(t_{i},\tau\right)\mathbf{B}_{u}\left(\tau\right)\mathbf{R}^{-1}\left(\tau\right) {\mathbf{B}_{u}}^{T}\left(\tau\right)\sum_{j=1}^{N}\bPhi^{T}\left(t_{j},\tau\right)\mathbf{H}^{T}1\left(\tau\leq t_{j}\right)\bmu_{j}d\tau&&\\
		&=\sum_{j=1}^{N}\int_{t_{0}}^{t_{f}} 1\left(\tau \leq t_{i}\right)1\left(\tau\leq t_{j}\right)\mathbf{H}\bPhi\left(t_{i},\tau\right)\mathbf{B}_{u}\left(\tau\right)\mathbf{R}^{-1}\left(\tau\right){\mathbf{B}_{u}}^{T}\left(\tau\right)\bPhi^{T}\left(t_{j},\tau\right)\mathbf{H}^{T}d\tau\bmu_{j}\\
		&:=\sum_{j=1}^{N} \bPsi_{ij}\bmu_{j}
	\end{aligned}
\end{equation}
where
\begin{equation} \label{Eq:Psi_ij}
	\bPsi_{ij}:=\mathbf{H}\int_{t_{0}}^{\min\left(t_{i},t_{j}\right)}\bPhi\left(t_{i},\tau\right)\mathbf{B}_{u}\left(\tau\right)\mathbf{R}^{-1}\left(\tau\right){\mathbf{B}_{u}}^{T}\left(\tau\right)\bPhi^{T}\left(t_{j},\tau\right)d\tau\mathbf{H}^{T}
\end{equation}
Note that the following property of the product of activator functions is used in the derivation of Eq. \eqref{Eq:Psi_ij}.
\begin{equation} \label{Eq:activ_prod}
	\int_{t_{0}}^{t_{f}}1\left(t\leq t_{i}\right)1\left(t\leq t_{j}\right) f\left(t\right) dt = \int_{t_{0}}^{\min\left(t_{i},t_{j}\right)} f\left(t\right)dt
\end{equation}
Although the coefficient matrix $\bPsi_{ij}$ is expressed as a definite integral in Eq. \eqref{Eq:Psi_ij}, it is more convenient to evaluate $\bPsi_{ij}$ through propagation of associated differential equation instead of directly computing the integral using quadrature methods. The integrand contains $\bPhi\left(t_{i}, \tau\right)$ which is already a result of propagation, and therefore, Eq. \eqref{Eq:Psi_ij} can be rewritten as
\begin{equation} \label{Eq:Psi_ij_eval}
	\bPsi_{ij} = \mathbf{H}\mathbf{U}\left(t_{i}\right) \mathbf{N}\left(\min\left(t_{i},t_{j}\right)\right)\mathbf{U}^{T}\left(t_{j}\right)\mathbf{H}^{T}
\end{equation}
with the fundamental solution matrix $\mathbf{U}\left(t\right)$ and the auxiliary matrix $\mathbf{N}\left(t\right)$ defined according to
\begin{align}
	\dot{\mathbf{U}}\left(t\right) &= \mathbf{A}_{u}\left(t\right)\mathbf{U}\left(t\right), 
	& \mathbf{U}\left(t_{0}\right) &= \mathbf{I} \label{Eq:U_dot}\\
	\dot{\mathbf{N}}\left(t\right) &= \mathbf{U}^{-1}\left(t\right)\mathbf{B}_{u}\left(t\right)\mathbf{R}^{-1}\left(t\right){\mathbf{B}_{u}}^{T}\left(t\right)\mathbf{U}^{-T}\left(t\right),
	& \mathbf{N}\left(t_{0}\right) &= \mathbf{0} \label{Eq:N_dot}
\end{align}
Indeed, Eq. \eqref{Eq:Psi_ij_eval} can be evaluated for each combination of $i$ and $j$ through a single forward sweep of the differential equations in Eqs. \eqref{Eq:U_dot} and \eqref{Eq:N_dot}. Therefore, the vertical concatenation of the relation in Eq. \eqref{Eq:constraints_u_opt_subs} followed by the matrix inverse yields the constraint-coupling multipliers as
\begin{equation} \label{Eq:bmu_opt}
	\bar{\bmu} = \bar{\bPsi}^{-1}\left(\mathbf{Z}_{v} - \mathbf{F}_{v}\right)
\end{equation}
with 
\begin{equation} \label{Eq:bmu_opt_aux}
		\bar{\bPsi} := \begin{bmatrix} \bPsi_{11} & \bPsi_{12} & \cdots & \bPsi_{1N}\\ 
					\bPsi_{21}	&	\bPsi_{22}	&	\cdots 	&	\bPsi_{2N}\\ 
					\vdots	&	\vdots	&	\ddots	&	\vdots\\
					\bPsi_{N1}	&	\bPsi_{N2}	&	\cdots 	&	\bPsi_{NN}
				\end{bmatrix},\quad
		\mathbf{Z}_{v} := \begin{bmatrix} \mathbf{z}_{1} - \mathbf{z}^{*}\left(t_{1}\right) \\ \cdots \\ \mathbf{z}_{N} - \mathbf{z}^{*}\left(t_{N}\right) \end{bmatrix},\quad
		\mathbf{F}_{v} := \begin{bmatrix} \mathbf{H}\bPhi\left(t_{1},t_{0}\right)\tilde{\mathbf{x}} \left(t_{0}\right) \\ \cdots \\ \mathbf{H}\bPhi\left(t_{N},t_{0}\right)\tilde{\mathbf{x}} \left(t_{0}\right)\end{bmatrix}
\end{equation}
Finally, substituting Eq. \eqref{Eq:bmu_opt} back into Eq. \eqref{Eq:FONC_u} yields the optimal control function correction as follows:
\begin{equation} \label{Eq:u_opt}
	\tilde{\mathbf{u}}\left(t\right) = \mathbf{R}^{-1}\left(t\right){\mathbf{B}_{u}}^{T}\left(t\right)\mathbf{E}\left(t\right)\bar{\bPsi}^{-1}\left(\mathbf{Z}_{v} - \mathbf{F}_{v}\right)
\end{equation}

\begin{rem} \label{Rem:FinalStateControl}
	If only a single interim point constraint is given for the final time $t_{f}$, the control function correction reduces to
	\begin{equation} \label{Eq:u_FinalStateControl}
		\tilde{\mathbf{u}}\left(t\right) = \mathbf{R}^{-1}\left(t\right){\mathbf{B}_{u}}^{T}\left(t\right)\bPhi^{T}\left(t_{f},t\right) \mathbf{H}^{T}\bPsi^{-1}\left(\mathbf{z}_{f} - \mathbf{z}^{*}\left(t_{f}\right) - \mathbf{H}\bPhi\left(t_{f},t_{0}\right)\tilde{\mathbf{x}} \left(t_{0}\right)\right)
	\end{equation}
	where
	\begin{equation}
		\bPsi = \mathbf{H}\int_{t_{0}}^{t_{f}}\bPhi\left(t_{f},\tau\right)\mathbf{B}_{u}\left(\tau\right)\mathbf{R}^{-1}\left(\tau\right){\mathbf{B}_{u}}^{T}\left(\tau\right)\bPhi^{T}\left(t_{f},\tau\right)d\tau\mathbf{H}^{T}
	\end{equation}
	which is the usual linear quadratic regulator solution for the terminal control problem.
\end{rem}

\begin{rem} \label{Rem:OutputControl}
	The proposed incremental correction method can be slightly modified to solve problems with the interim point constraints imposed on the output $\mathbf{z}\left(t\right)=\mathbf{h}\left(\mathbf{x}\left(t\right)\right)$ for a nonlinear function $\mathbf{h}\left(\cdot\right)$. This can be done by computing the influence functions, i.e., sensitivity functions, associated with the output variables. In the special case where only a single constraint is imposed on the output at the final time, the procedure to determine control function correction becomes identical to the generalised model predictive static programming technique \cite{Maity_2014}.
\end{rem}

\begin{rem} \label{Rem:AlgorithmParams}
	Each correction method has its own computational challenges. As mentioned earlier for the parameter correction method, an automatic differentiation tool can be exploited to compute the sensitivity $\left.\frac{\partial \mathbf{x}\left(t_{i}\right)}{\partial \btheta}\right|_{\btheta^{*}}$ with the capability to differentiate through the ODE solver. One should determine which of the available local sensitivity analysis methods suits the purpose depending on the application. Also, the Jacobian array $\mathbf{L}$ in Eq. \eqref{Eq:tilde_btheta_opt} is a large fat matrix since the number of NN parameters is usually much greater than the number of control inputs in physical systems. Its pseudoinverse computation might show considerable dependence on the relative tolerance setting. Therefore, the tolerance should be selected appropriately to improve constraint satisfaction accuracy reliably through the parameter correction method. 
	
	On the other hand, the control function correction method can suffer from numerical instabilities in propagation of Eqs. \eqref{Eq:U_dot} and \eqref{Eq:N_dot} to obtain matrices $\mathbf{U}\left(t\right)$ and $\mathbf{N}\left(t\right)$, respectively, when the ODE solver employs too tight relative and absolute tolerance setting. Although the choice of weighting function $\mathbf{R}\left(\tau\right)$ itself is not a source of computational challenges, appropriate selection that can bring improvements in the constraint targeting accuracy as expected might be challenging when coupled with the aforementioned numerical instabilities.
\end{rem}

\begin{rem} \label{Rem:x_0}
	Sections \ref{SubSec:ParamCorr} and \ref{SubSec:CtrlFuncCorr} do not necessarily assume $\mathbf{x}_{0} = \mathbf{x}_{0}^{*}$, thus the term $\tilde{\mathbf{x}} \left(t_{0}\right)$ appears in the command equations. The initial state perturbation term realises tracking error feedback when $t_{0}$ and $\mathbf{x}_{0}$ in the expressions are viewed as the current time and state, respectively, at each instance. Therefore, the control laws can be implemented in a feedback form based on continuous re-initialisation at the expense of increased online computation burden.
\end{rem}

\section{Application} \label{Sec:Applications}
This section presents examples to demonstrate the effectiveness of the proposed correction methods in a practical application. A powered descent problem on a planetary surface \cite{Amato_2020,Amato_2021,McMahon_2022} is considered for illustration as the landing accuracy is of interest in this problem.

\subsection{Problem Description and System Model}
This application addresses the problem of powered descent guidance for Mars landing. The objective of this finite-horizon control problem is to find a neural network policy that achieves the desired final position and velocity at a fixed final time while minimising the fuel consumption. 

Consider the vehicle motion with respect to a Mars-centred, Mars-fixed coordinate system. The spacecraft dynamics can be expressed as
\begin{equation} \label{Eq:sys_dyn}
	\begin{aligned}
		\dot{\mathbf{r}} &= \mathbf{v}\\
		\dot{\mathbf{v}} &= -\frac{\mu}{r^{3}}\mathbf{r} + \frac{\mathbf{F}_{L} + \mathbf{F}_{D} + \mathbf{F}_{T}}{m} - 2\mathbf{\Omega} \times \mathbf{v} - \mathbf{\Omega} \times \left(\mathbf{\Omega} \times \mathbf{r}\right)\\
		\dot{m} &= -\frac{T}{I_{sp}g}
	\end{aligned}
\end{equation}
where $\mathbf{r}$ and $\mathbf{v}$ denote the position and velocity vectors with respect to the centre of Mars, respectively, $r \triangleq \left\| \mathbf{r} \right\|$ denotes the radial distance, $\mathbf{\Omega} \triangleq \begin{bmatrix}	0	&	0	&	\Omega \end{bmatrix}^{T}$ denotes the angular velocity for Mars rotation, $m$ denotes the vehicle mass. The lift, drag, and thrust forces are represented by $\mathbf{F}_{L}$, $\mathbf{F}_{D}$, and $\mathbf{F}_{T}$ in Eq. \eqref{Eq:sys_dyn}, respectively, with $T \triangleq \left\| \mathbf{F}_{T} \right\|$, $I_{sp}$, and $g$ representing the thrust magnitude, specific impulse, and gravitational acceleration at the Earth surface, respectively. 

\begin{figure}[!ht]
	\begin{center}
		\includegraphics[width=0.6\textwidth]{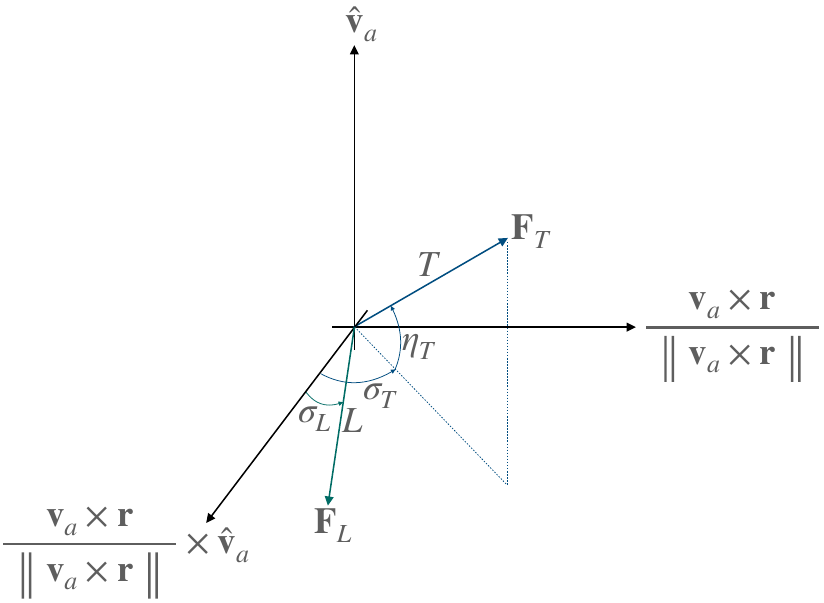}
		\caption{Lift and Thrust Force Resolved in Wind Axes} \label{Fig_F}
	\end{center}
\end{figure}

Consider the wind axes and the force direction angles that are defined according to Fig. \ref{Fig_F}. 

\noindent The three-dimensional force vectors can be described as
\begin{equation} \label{Eq:force}
	\begin{aligned}
		\mathbf{F}_{L} &= L\left(\cos\sigma_{L}  \frac{\mathbf{v}_{a} \times \mathbf{r}}{\left\| \mathbf{v}_{a} \times \mathbf{r}\right\|} \times \hat{\mathbf{v}}_{a} + \sin\sigma_{L} \frac{\mathbf{v}_{a} \times \mathbf{r}}{\left\| \mathbf{v}_{a} \times \mathbf{r}\right\|}\right)\\
		\mathbf{F}_{D} &= -D \hat{\mathbf{v}}_{a}\\
		\mathbf{F}_{T} &= T \left( \cos\eta_{T}\cos\sigma_{T}  \frac{\mathbf{v}_{a} \times \mathbf{r}}{\left\| \mathbf{v}_{a} \times \mathbf{r}\right\|} \times \hat{\mathbf{v}}_{a} + \cos\eta_{T}\sin\sigma_{T} \frac{\mathbf{v}_{a} \times \mathbf{r}}{\left\| \mathbf{v}_{a} \times \mathbf{r}\right\|} + \sin\eta_{T}\hat{\mathbf{v}}_{a} \right)
	\end{aligned}
\end{equation}
where $\mathbf{v}_{a} \triangleq \mathbf{v} - \mathbf{v}_{w}$ denotes the relative velocity with respect to the surrounding air flow of $\mathbf{v}_{w}$, while the hat notation $\hat{\left(\cdot\right)}$ refers to the unit vector in the direction of a quantity. Also, $\sigma_{T}$ and $\eta_{T}$ are the azimuth and elevation angles for the thrust, respectively, and $\sigma_{L}$ is the bank angle. The magnitudes of the lift and drag can be expressed as
\begin{equation} \label{Eq:LD}
	\begin{aligned}
		L &= \frac{1}{2}\rho {\left\| \mathbf{v}_{a} \right\|}^{2} C_{L} S = R_{L/D} D\\
		D &= \frac{1}{2}\rho {\left\| \mathbf{v}_{a} \right\|}^{2} C_{D} S = \frac{m}{2\beta}\rho {\left\| \mathbf{v}_{a} \right\|}^{2} \\
		\rho\left(h\right) &= \rho_{0} \exp\left(-\frac{h}{H}\right)
	\end{aligned}
\end{equation}
where $R_{L/D}$ is the lift-to-drag ratio, $\beta\left(t\right) = \frac{m\left(t\right)}{C_{D}S}$ is the ballistic coefficient of the vehicle, and $\rho$ is the atmospheric density that is approximately modelled as a decreasing exponential function of the altitude $h$ measured above the surface.

\begin{figure}[!ht]
	\begin{center}
		\includegraphics[width=0.4\textwidth]{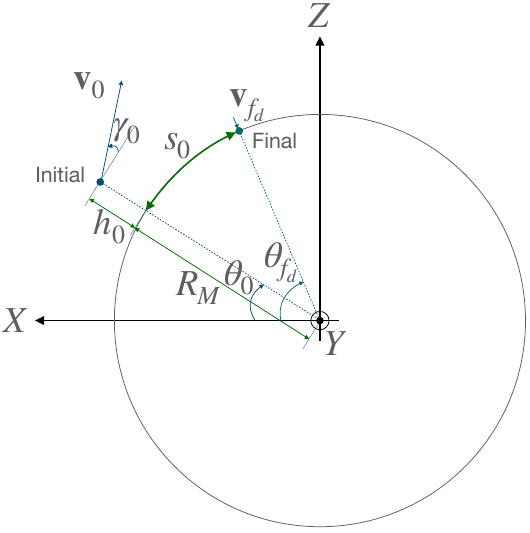}
		\caption{Mars Landing Guidance Problem Geometry (The circle denotes the planetary surface.)} \label{Fig_MCMF}
	\end{center}
\end{figure}

In this example, it is assumed for simplicity that both the vehicle position and velocity initially lie on the $XZ$-plane and the desired final position and velocity are also on the same plane as shown in Fig. \ref{Fig_MCMF}.
In Fig. \ref{Fig_MCMF}, $R_{M}$, $\theta$, $\gamma$, and $s \triangleq R_{M}\cos^{-1}\left(\hat{\mathbf{r}} \cdot \hat{\mathbf{r}}_{f_{d}}\right)$ denote the planet radius, latitude, flight path angle measured in local vertical/horizontal frame, and ground track distance to the target position. The subscripts $0$ and $f_{d}$ refer to the quantities pertaining to the initial and desired final conditions, respectively. The initial and final conditions are given by
\begin{equation} \label{Eq:IC_FC}
	\begin{aligned}
		\mathbf{r}_{0} &= \left(R_{M}+h_{0}\right)\begin{bmatrix}
			\cos\theta_{0}\\
			0\\
			\sin\theta_{0}
		\end{bmatrix}, &
		\mathbf{v}_{0} &= V_{0} \begin{bmatrix}
			-\sin\left(\theta_{0} - \gamma_{0}\right)\\
			0\\
			\cos\left(\theta_{0} - \gamma_{0}\right)
		\end{bmatrix}, &
		\mathbf{r}_{f_{d}} &= R_{M}\begin{bmatrix}
		\cos\theta_{f_{d}}\\
		0\\
		\sin\theta_{f_{d}}
		\end{bmatrix}, &
		\mathbf{v}_{f_{d}} &= -V_{f_{d}} \hat{\mathbf{r}}_{f_{d}}
	\end{aligned}
\end{equation}	
where $V \triangleq \left\|\mathbf{v}\right\|$ denotes the ground speed. Note that one can specify either $\theta_{0}$ or $s_{0}$ to prescribe the initial position since $s_{0}=R_{M}\left(\theta_{f_{d}}-\theta_{0}\right)$. Also note that the performance output constraints at the final time consist of the desired final position and velocity conditions, leaving the final mass unconstrained.

\begin{figure}[!htbp]
	\begin{center}
		\includegraphics[width=0.7\textwidth]{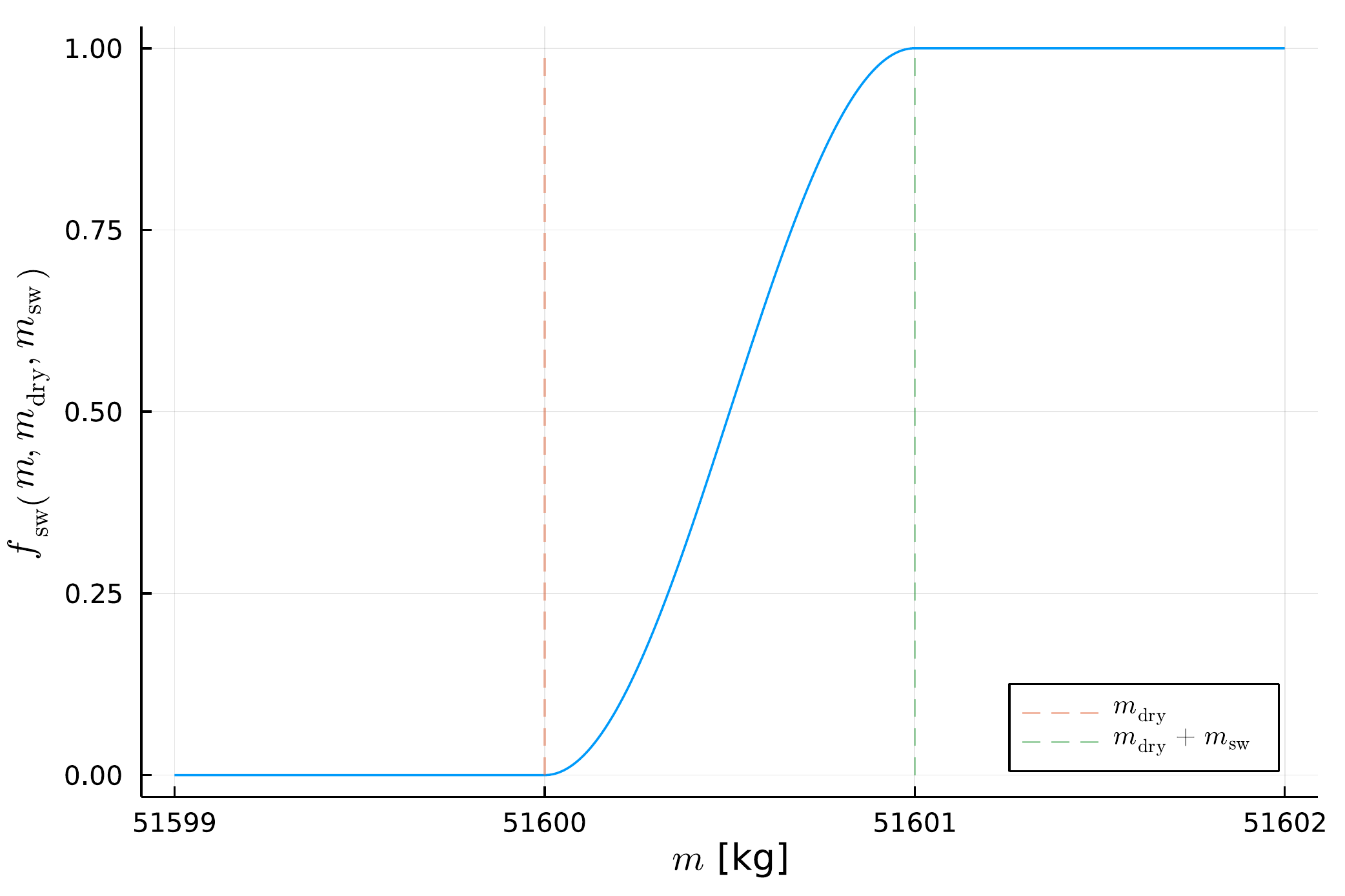}
		\caption{Illustrative Example of Thrust Activation Function}
		\label{Fig_sw}
	\end{center}
\end{figure}

The policy is defined to be a fully-connected feedforward neural network which takes the normalised position and velocity errors as inputs and generates the throttle command, i.e., normalised thrust magnitude, and thrust direction angles as its outputs.
\begin{equation} \label{Eq:pi}
	\begin{bmatrix}
		\delta_{T}\\
		\sigma_{T}\\
		\eta_{T}
	\end{bmatrix}
	= \boldsymbol{\pi}\left(\begin{bmatrix} 
		\displaystyle \frac{\mathbf{r}-\mathbf{r}_{f_{d}}}{s_{0}} & 
		\displaystyle \frac{\mathbf{v}-\mathbf{v}_{f_{d}}}{V_{0}}  
	\end{bmatrix}, \btheta \right)
\end{equation}
The activation function for the output layer is chosen specifically to confine its range to a bounded interval $\left(y_{lb}, y_{ub}\right)$ as
\begin{equation} \label{Eq:f_scale}
	f_{\text{scale}} \left(u\right) = \frac{y_{ub} - y_{lb}}{1+\exp\left(-u\right)} + y_{lb}
\end{equation}
This architectural choice is made to scale the policy network outputs so that the magnitude and the direction angles of the thrust vector naturally satisfy the physical limits given by
\begin{equation} \label{Eq:thrust_lim}
	\begin{aligned}
		\frac{T_{\min}}{T_{\max}} &\leq \delta_{T} \leq 1\\
		\eta_{\min} &\leq  \eta_{T} \leq \eta_{\max}\\
		\sigma_{\min} &\leq  \sigma_{T} \leq \sigma_{\max}
	\end{aligned}
\end{equation}
In Eq. \eqref{Eq:thrust_lim}, $T_{\min}$ and $T_{\max}$ are the minimum and maximum thrust magnitude, respectively. Likewise, $\eta_{\min}$, $\eta_{\max}$, $\sigma_{\min}$, and $\sigma_{\max}$ are limits imposed on the thrust elevation and azimuth angles. The throttle command is related to the thrust magnitude through an activation function which nullifies the thrust in the absence of fuel. A smooth activation function defined by 
\begin{equation} \label{Eq:f_sw}
	f_{\text{sw}}\left(m, m_{\text{dry}}, m_{\text{sw}}\right) = \frac{1-\cos\left( \frac{\min\left(\max\left(m-m_{\text{dry}} ,0\right), m_{\text{sw}}\right)}{m_{\text{sw}}}\pi \right)}{2}
\end{equation}
is introduced instead of the discrete switching function to keep the gradient computation required in backpropagation through the ODE dynamics well-defined. In Eq. \eqref{Eq:f_sw}, $m_{\text{dry}}$ denotes the dry mass, and $m_{\text{sw}} > 0$ denotes the design parameter which determines the steepness of the activation function. Figure \ref{Fig_sw} shows an example plot of the thrust activation function in Eq. \eqref{Eq:f_sw} with $m_{\text{dry}} = 51\,600$kg and $m_{\text{sw}} = 1$kg. The thrust model is thus given by
\begin{equation} \label{Eq:T_deltaT}
	T =  f_{\text{sw}}\left(m, m_{\text{dry}}, m_{\text{sw}}\right) T_{\max} \delta_{T}
\end{equation}

\subsection{Simulation Setup}
For the powered descent landing guidance problem, the purpose of the policy is to achieve the desired position and velocity at the given final time as accurately as possible while minimising the fuel expenditure. This objective is encoded into the cost function for baseline policy optimisation which is defined in the form of a weighted sum as
\begin{equation} 
		J = k_{r_{f}} \frac{\left\|\mathbf{r}\left(t_{f}\right)-\mathbf{r}_{f_{d}}\right\|^{2}}{{s_{0}}^{2}} + k_{v_{f}} \frac{\left\|\mathbf{v}\left(t_{f}\right)-\mathbf{v}_{f_{d}}\right\|^{2}}{{V_{0}}^{2}}  + k_T \int_{t_{0}}^{t_{f}} \delta_{T}\left(\tau\right) d\tau + k_{\theta} \left\|\boldsymbol{\theta}\right\|^{2}
\end{equation}
where $k_{r_{f}}$, $k_{v_{f}}$, $k_{T}$, and $k_{\theta}$ are the positive constant weights. The $L_{2}$-norm of the neural network parameters is included in the cost for regularisation. The cost function is chosen to account for soft constraints on the final state by penalising the errors with very large weights $k_{r_{f}}$ and $k_{v_{f}}$.

The continuous-time policy gradient method developed in \cite{Cho_2022} based on adjoint sensitivity analysis techniques is used for training of the policy. The configuration of the policy neural network and the optimiser setup used for its training are summarised in Table \ref{Table:TrainConfig}. The policy parameter vector $\boldsymbol{\theta}$ is optimised by simulating the controlled trajectory at each training iteration with the stopping criteria for integration of individual trajectory given by $\left(t\geq43\text{s}\right) ~\lor~\left(\left\|\mathbf{r} - \mathbf{r}_{f_{d}}\right\| \leq 100\text{m}~\land~\mathbf{v} \cdot \left( \mathbf{r} - \mathbf{r}_{f_{d}} \right) \geq 0\right)$. Incremental correction for computing either $\tilde{\btheta}$ or $\tilde{\mathbf{u}}\left(t\right)$ is applied only once at the initial time to clearly compare the two different correction methods (see Remark \ref{Rem:x_0} for the note on continuous re-initialisation). A relative tolerance of $0.005$ is used for computing the pseudoinverse in the parameter correction method. The model data for the dynamical system found in \cite{Lu_2019, Amato_2021}, including the Martian environment and the vehicle initial conditions, are considered for the illustrative example. The initial and final conditions defining the mission scenario, the parameters for vehicle dynamics model and policy design, and the parameters for environmental physics model are described in Tables \ref{Table:SimParam_Mission}-\ref{Table:SimParam_Environment}, respectively. 

Two different cases are considered for simulation in order to consider different test purposes. The main difference in the simulation scenarios is in the initial condition used for testing:
\begin{itemize}
	\item Case 1. single identical initial condition for both baseline policy training and testing
	\item Case 2. multiple initial conditions for testing that are perturbed from initial condition for baseline policy training
\end{itemize}

Case 1 in Sec. \ref{SubSubSec:Sim_Case1} aims to test the effectiveness of two incremental correction methods in improving the actual accuracy in satisfying the performance output constraint at the final time. The secondary purpose of Case 1 is to demonstrate the full process of NN policy design using the proposed two-stage approach specifically for the powered descent application. The performance of the NN policy is evaluated with various random seeds used for NN weight initialisation considering fixed initial and final conditions. 

Case 2 in Sec. \ref{SubSubSec:Sim_Case2} aims to show the utility of incremental correction in dealing with the adverse effect of the inital position dispersion in landing accuracy. In practice, imperfect handover from the entry guidance phase and uncertainties such as environmental disturbances cause a dispersion in the initial condition for the powered descent phase. As a simple model for the dispersion, the closed-loop system is propagated for a set of initial positions obtained along the circle of radius $100$m which lies on the plane perpendicular to the initial velocity and is centred at the nominal initial position. The performance of the baseline policy alone, the policies updated with parameter correction and control function correction are compared with each other. The NN training configuration and simulation parameters are identical across both cases, except for the initial position.
	
\begin{table}[h!tb]
	\begin{center}
		\caption{Training Configuration} \label{Table:TrainConfig}
		\begin{tabular}{cc}
			\hline
			Object 				&	Value\\
			\hline
			$\boldsymbol{\pi}$ Input layer			&	10 tanh\\
			$\boldsymbol{\pi}$ Hidden 1st layer		&	10 tanh\\
			$\boldsymbol{\pi}$ Hidden 2nd layer		&	3 linear\\
			$\boldsymbol{\pi}$ Output layer			&	3 scale\\
			1st optimiser				&	ADAM\\
			2nd optimiser				&	BFGS\\
			$\left(k_{r_{f}}, k_{v_{f}}, k_T, k_\theta \right)$	&	$\left(10^{6},10^{5},1,10^{-6}\right)$\\
			\hline
		\end{tabular}
	\end{center}
\end{table}

\begin{table}[h!tb]
	\begin{center}
		\caption{Simulation Parameters: Initial and Final Conditions} \label{Table:SimParam_Mission}
		\begin{tabular}{ccc|ccc}
			\hline
			Parameter 			&	Value					&	Unit & 		Parameter 			&	Value					&	Unit\\
			\hline
			$h_{0}$				&	$2480$					&	m &
			$V_{0}$				&	$505$					&	m/s\\
			$\gamma_{0}$		& 	$0$						&	deg &
			$s_{0}$				&	$11500$					&	m\\
			$m_{0}$				&	$62000$					&	kg &
			$V_{f_{d}}$			&	$2.5$					&	m/s\\
			$\theta_{f_{d}}$	&	$45$					&	deg &
			$t_{f}$				& 	$43$					& s\\
			\hline
		\end{tabular}
	\end{center}
\end{table}

\begin{table}[h!tb]
	\begin{center}
		\caption{Simulation Parameters: Vehicle Dynamics and Policy} \label{Table:SimParam_Vehicle}
		\begin{tabular}{ccc|ccc}
			\hline
			Parameter 			&	Value					&	Unit & 		Parameter 			&	Value					&	Unit\\
			\hline
			$m_{\text{dry}}$	&	$51600$					&	kg &
			$m_{\text{sw}}$		&	$1$						&   kg	\\
			$I_{sp}$			& 	$360$					&	s  &
			$\beta_{0}$			&	$379$					&	kg/m$^2$ \\
			$R_{L/D}$			&	$0.54$					&	-  &
			$\sigma_{L}$		&	$0$						&	deg \\
			$T_{\max}$			&	$8\times10^{5}$			&	N &
			$T_{\min}$			&	$0.2T_{\max}$			&	N\\
			$\eta_{\max}$		&	$90$					&	deg &
			$\eta_{\min}$		&	$-90$					&	deg\\
			$\mathbf{R}\left(t\right)$ & $\mathrm{diag}\left(10,1,1\right)$ 		&   - \\
			\hline
		\end{tabular}
	\end{center}
\end{table}

\begin{table}[h!tb]
	\begin{center}
		\caption{Simulation Parameters: Environmental Physics} \label{Table:SimParam_Environment}
		\begin{tabular}{ccc}
			\hline
			Parameter 			&	Value					&	Unit\\
			\hline
			$\mathbf{v}_{w}$ 	& 	$\mathbf{0}$ 			& m/s \\
			$\mu$				&	$4.282837\times 10^{13}$&	m$^3$/s$^2$\\
			$R_{M}$					&	$3389.5 \times 10^{3}$	&	m \\
			$\Omega$			&	$2\pi / 1.025957$		&	rad/d\\
			$\rho_{0}$			&	$0.0263$				&	kg/m$^3$ \\
			$H$					&	$10153.6$				&	m\\
			$g$					&	$9.805$					&	m/s$^2$\\
			\hline
		\end{tabular}
	\end{center}
\end{table}

\begin{rem} \label{Rem:NumericIssues}
	There are several points to note here in relation to the computational issues observed by trial and error that should be acknowledged for better practical implementation. First, when single-precision floating point data type is used for real numbers to reduce memory burden, the quantities which involve values in the order of planet radius $R_{M}$ in calculation should be avoided from being included in the policy representation. The lack of precision due to low resolution manifests itself in the policy training as high numerical sensitivity in the gradient computation and in the variable-step integration of the ODEs. Second, parametrising the policy with respect to the groundtrack distance $s$ defined by $R_{M}\cos^{-1}\left(\hat{\mathbf{r}}\cdot\hat{\mathbf{r}}_{f_{d}}\right)$ may trigger faulty numerical behaviours during training due to the inaccuracy of arccosine function in some scientific computing systems. Third, the input and output variables of the neural network would better be of a similar order of magnitude to avoid the gradient computed through the neural network being sensitive only to a certain variable. This is the reason for introducing normalisation in the inputs and choosing throttle command $\delta_{T}$ instead of thrust magnitude $T$ as an output.
\end{rem}

\subsection{Simulation Results}
\subsubsection{Case 1. Single Matched Initial Condition} \label{SubSubSec:Sim_Case1}
\begin{figure}[!ht]
	\begin{center}
		\includegraphics[width=0.7\textwidth]{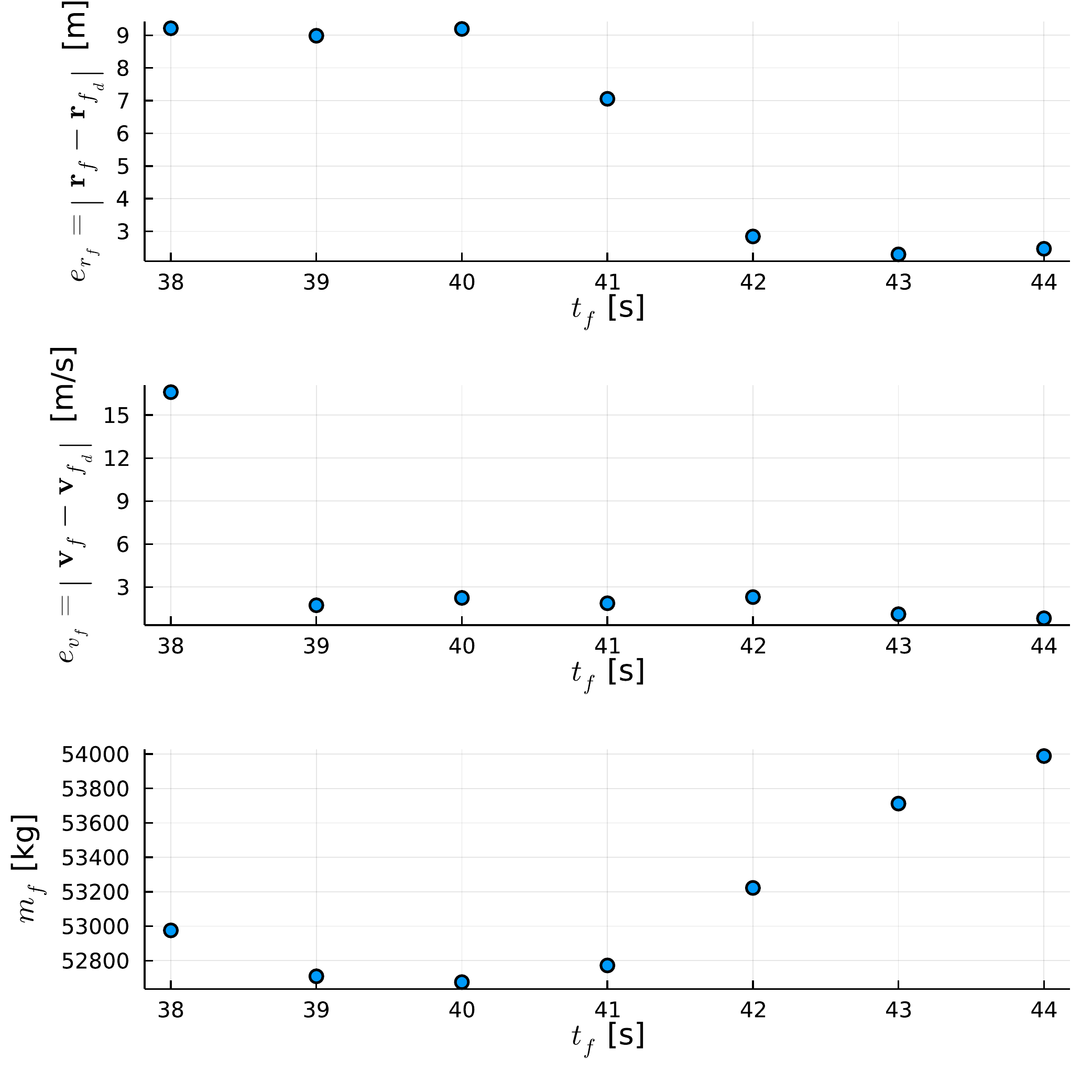}
		\caption{Final Errors and Mass with Baseline Policy} \label{Fig_e_rvm_f_baseline}
	\end{center}
\end{figure}

\begin{figure}[!ht]
	\begin{center}
		\includegraphics[width=0.6\textwidth]{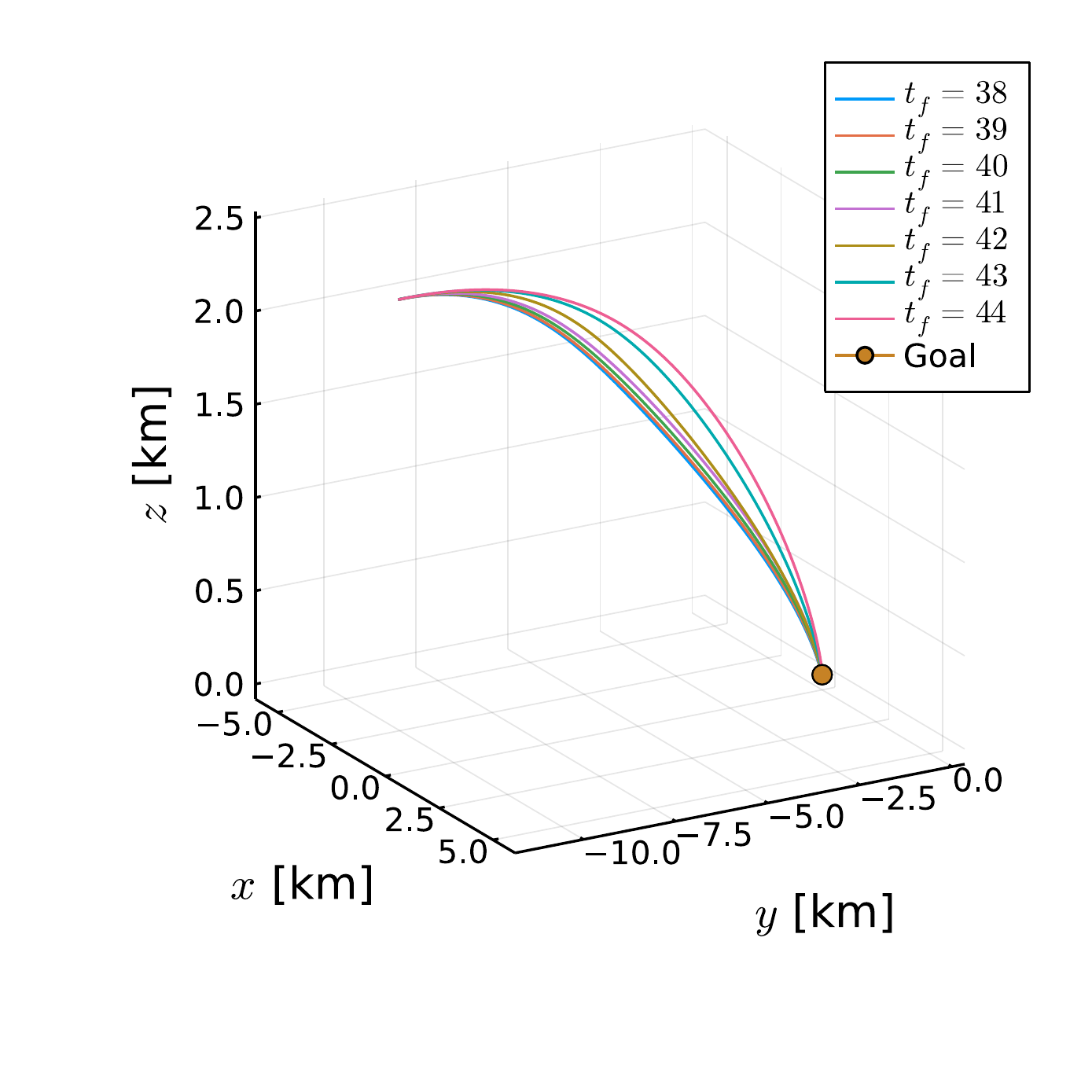}
		\caption{Three-Dimensional Trajectory with Baseline Policy} \label{Fig_3D_baseline}
	\end{center}
\end{figure}

\begin{figure}[!ht]
	\begin{center}
		\includegraphics[width=0.7\textwidth]{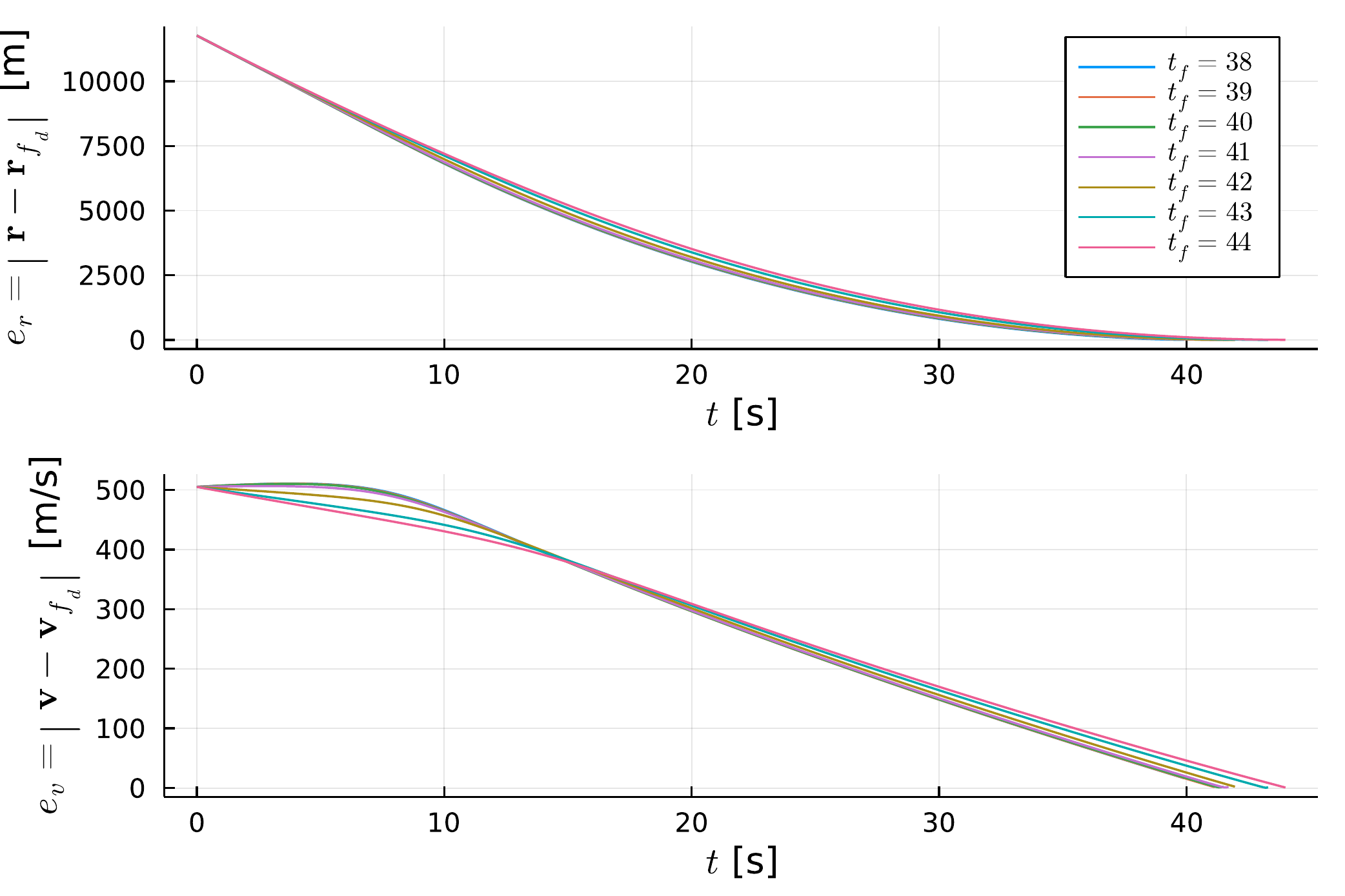}
		\caption{Error History with Baseline Policy} \label{Fig_e_baseline}
	\end{center}
\end{figure}

\begin{figure}[!ht]
	\begin{center}
		\includegraphics[width=0.7\textwidth]{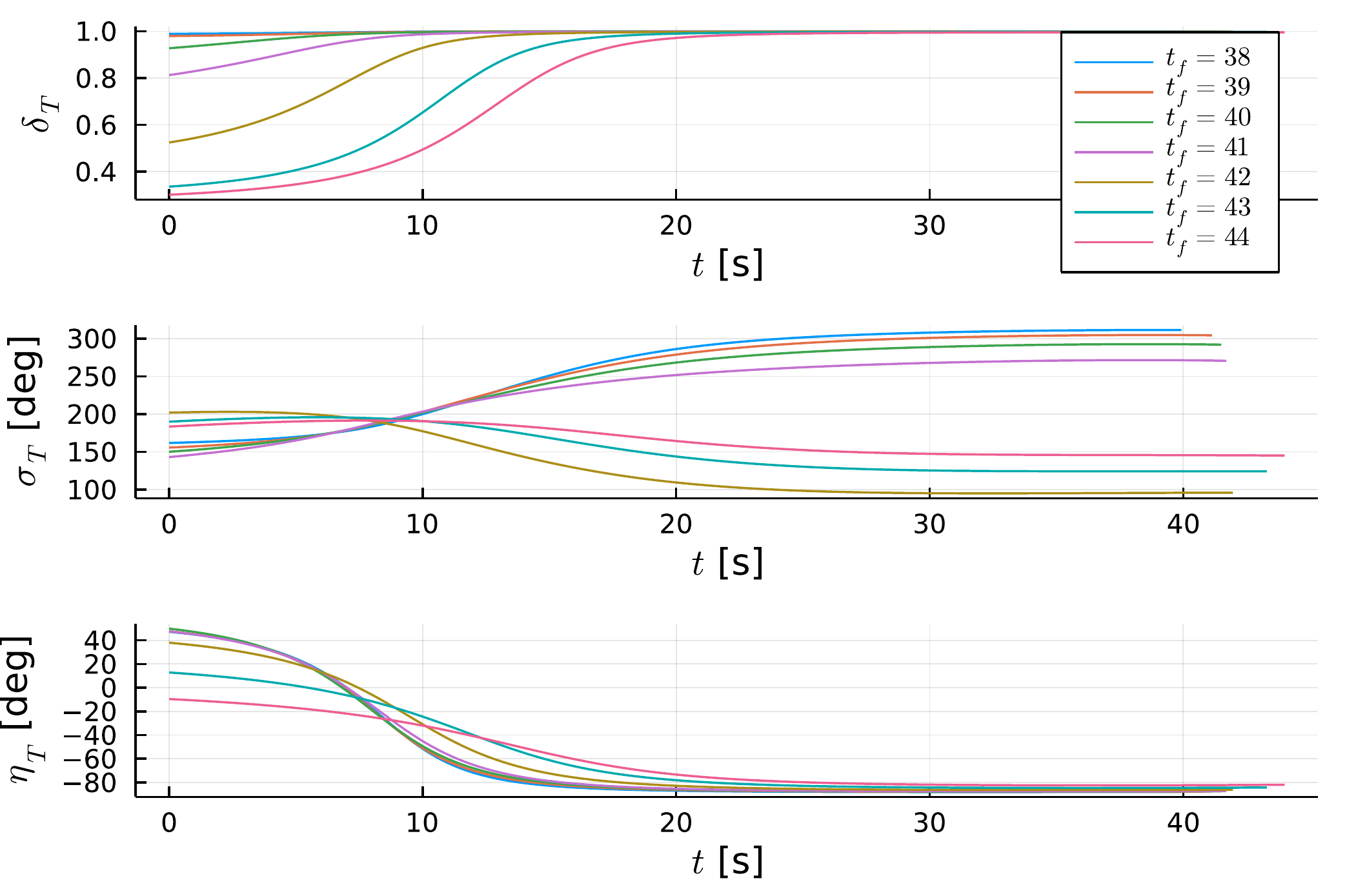}
		\caption{Control History with Baseline Policy} \label{Fig_u_baseline}
	\end{center}
\end{figure}

\paragraph{Stage 1. Baseline Policy Training with Various Values of $t_{f}$}
The available thrust magnitude of the spacecraft is limited, leading to limitations on the region reachable within a fixed flight time. Baseline policy training for Stage 1 is performed with various values of $t_{f}$ ranging from $38$s to $44$s to determine a physically feasible solution. Figures \ref{Fig_e_rvm_f_baseline}-\ref{Fig_u_baseline} show final errors and mass, three-dimensional trajectory, error history, and control history, respectively, obtained with the baseline policies for each $t_{f}$. Figure \ref{Fig_3D_baseline} shows that the vehicle approaches close to target without unnecessarily consuming much control effort in the horizontal direction with the baseline policies for each $t_{f}$. The position and the velocity errors that are defined by the norm of the vector difference from their respective desired final value tend to zero as shown in Fig. \ref{Fig_e_baseline}. The time histories of the thrust magnitude and direction angles shown in Fig. \ref{Fig_u_baseline} describe a continuous trend dependent upon $t_{f}$. In all cases, the optimised policy exhibits maximal throttle command $\delta_{T}$ and minimal thrust elevation angle $\eta_{T}$ as the vehicle approaches the end of flight. The final position error as well as the velocity error are the least at $t_{f} = 43$s as shown in Fig. \ref{Fig_e_rvm_f_baseline}. Therefore, $t_{f} = 43$s is chosen for testing Stage 2 methods in the following.

\paragraph{Stage 2. Incremental Correction with Fixed $t_{f}$ and Various Random Seeds}
The result of baseline policy optimisation in Stage 1 substantially depends on initialisation of the NN parameters, which is done by random sampling. Here, the entire process of the proposed two-stage approach is repeated with different random seeds while fixing $t_{f} = 43$s. Both the parameter and the control function correction methods are tested to demonstrate their characteristics. 

\begin{figure}[!ht]
	\begin{center}
		\includegraphics[width=0.7\textwidth]{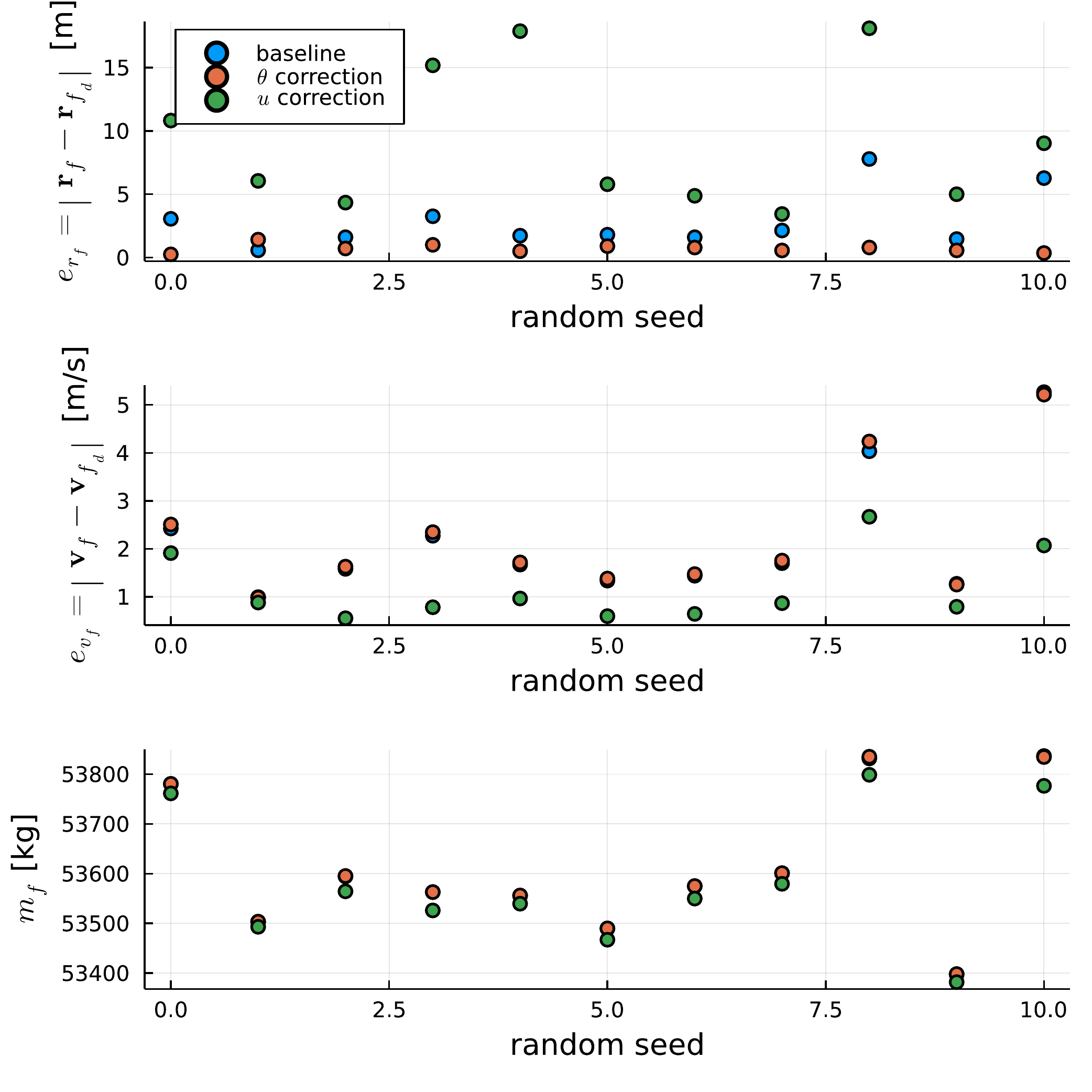}
		\caption{Final Errors and Mass with Incremental Correction (Case 1)} \label{Fig_e_rvm_f_C1}
	\end{center}
\end{figure}

\begin{figure}[!ht]
	\begin{center}
		\includegraphics[width=0.6\textwidth]{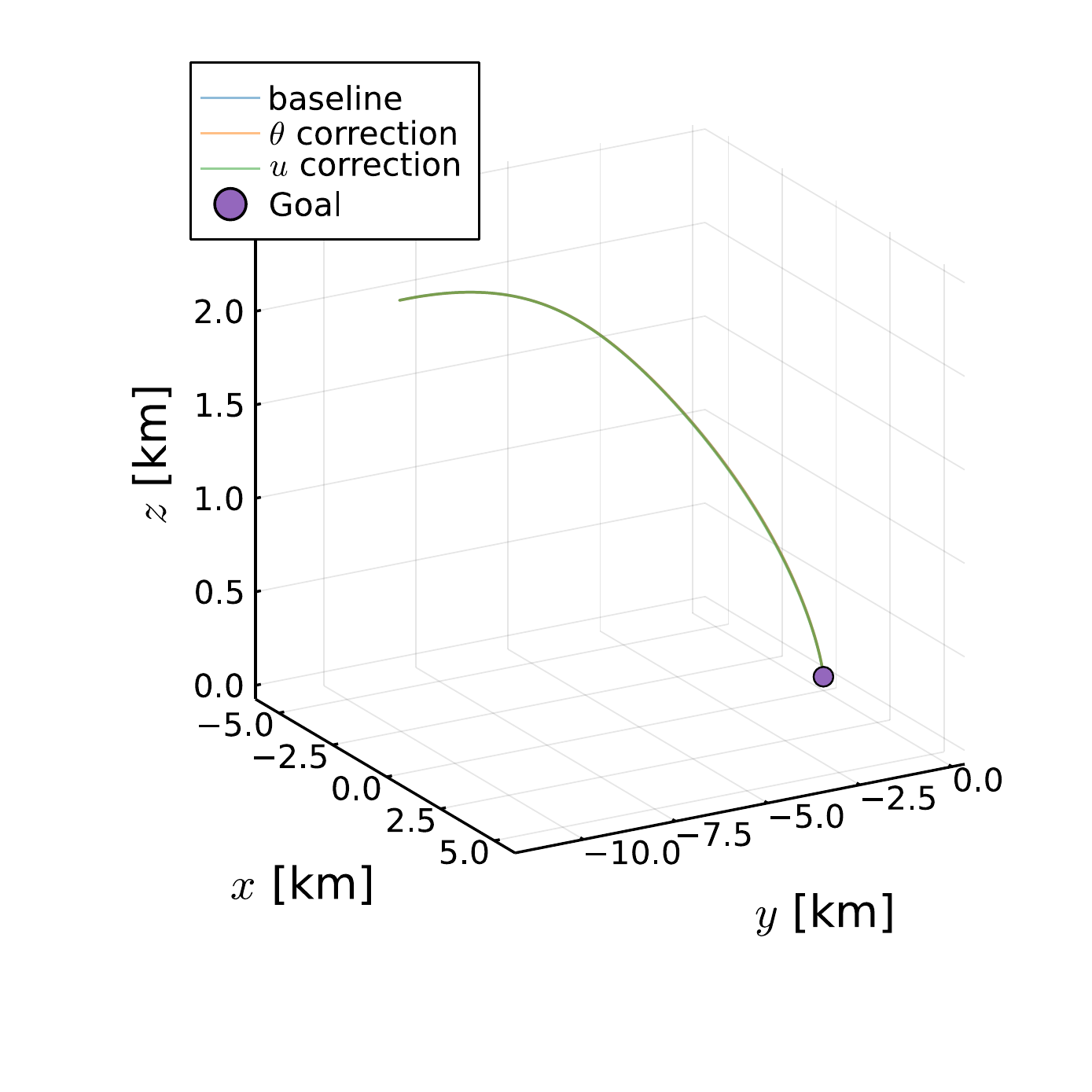}
		\caption{Three-Dimensional Trajectory with Incremental Correction - Single Random Seed (Case 1)} \label{Fig_3D_C1}
	\end{center}
\end{figure}

\begin{figure}[!ht]
	\begin{center}
		\includegraphics[width=0.7\textwidth]{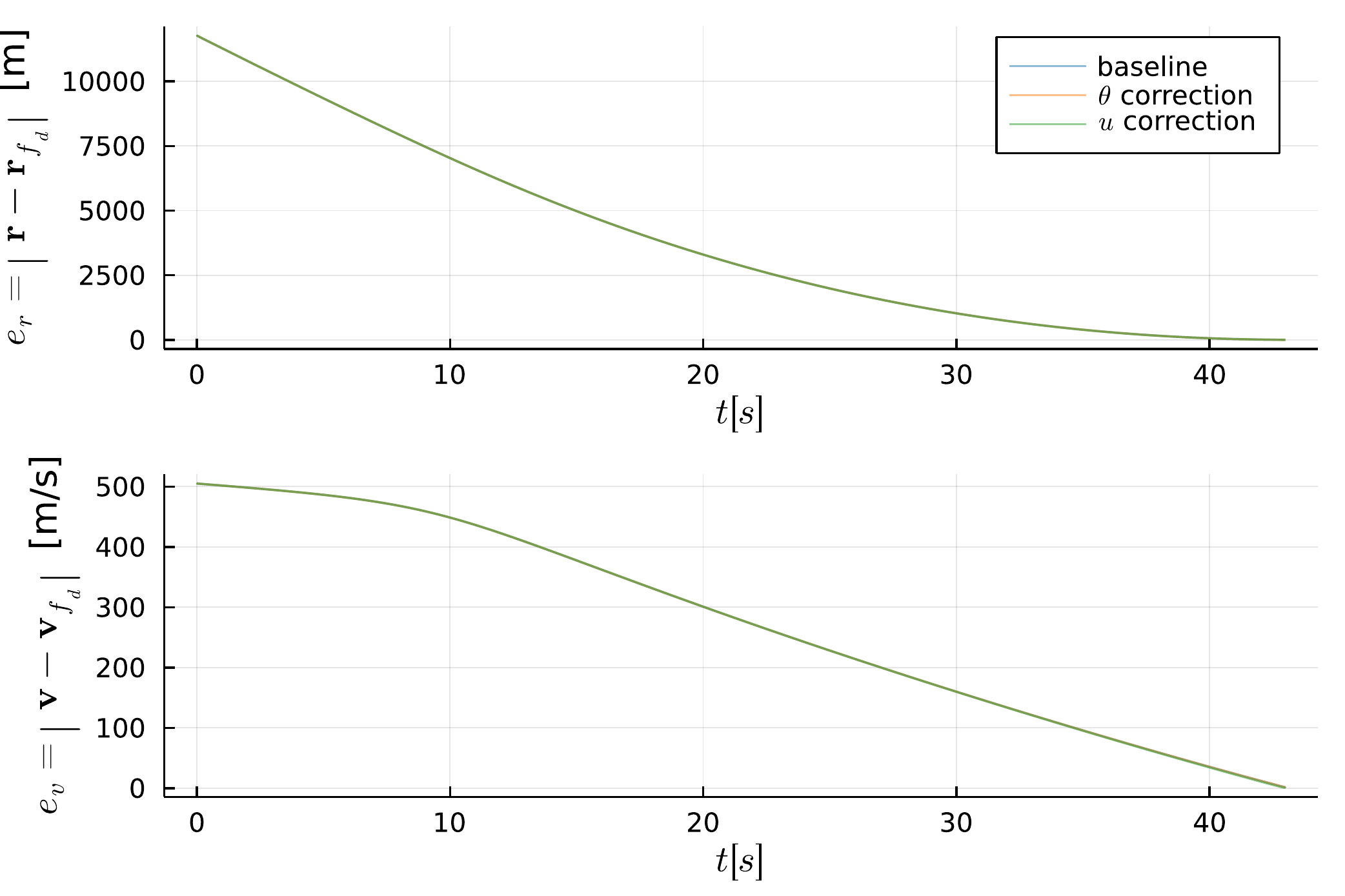}
		\caption{Error History with Incremental Correction - Single Random Seed (Case 1)} \label{Fig_e_C1}
	\end{center}
\end{figure}

\begin{figure}[!ht]
	\begin{center}
		\includegraphics[width=0.7\textwidth]{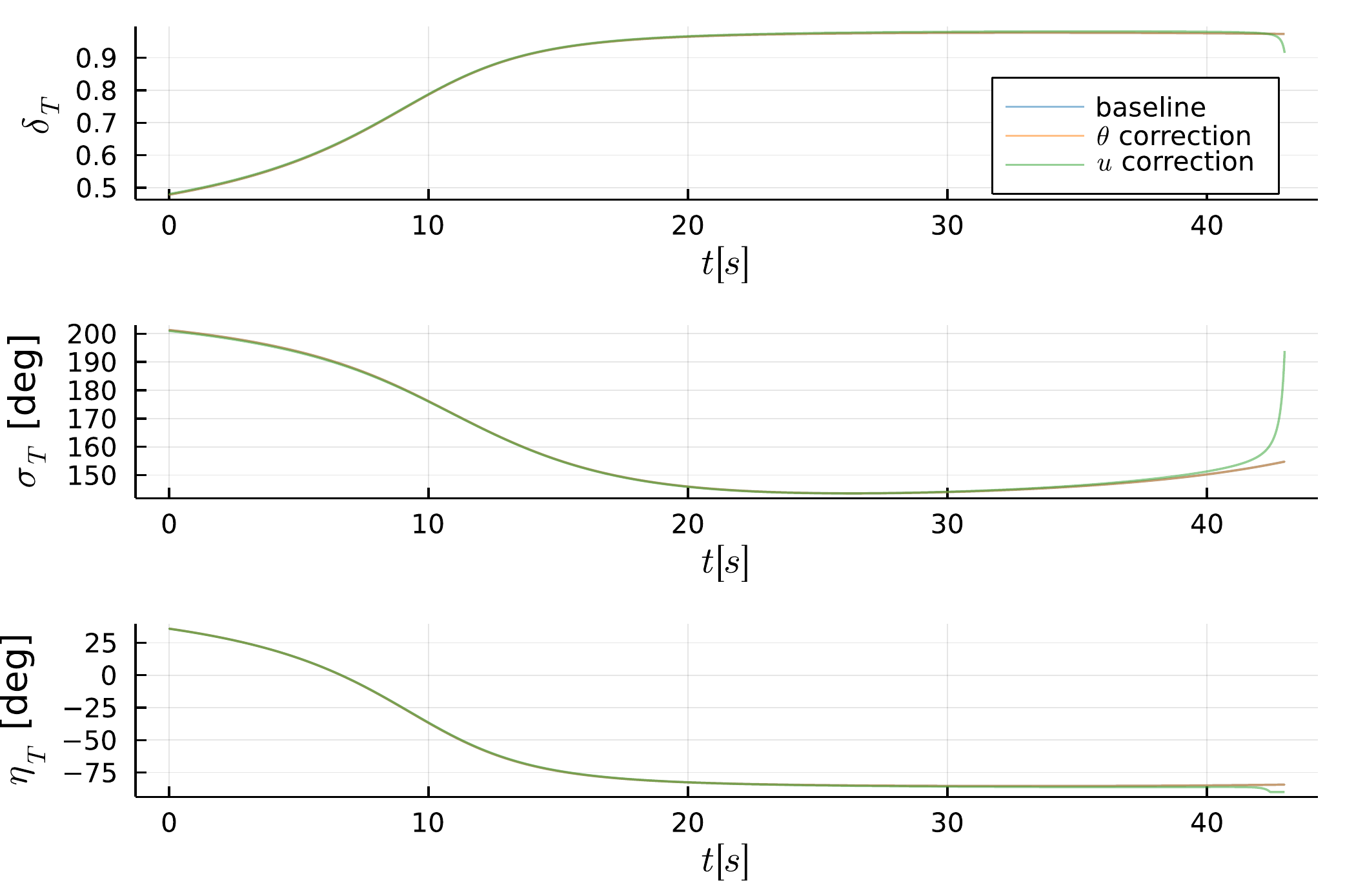}
		\caption{Control History with Incremental Correction - Single Random Seed (Case 1)} \label{Fig_u_C1}
	\end{center}
\end{figure}

Figure \ref{Fig_e_rvm_f_C1} shows final errors and mass for different random seeds. Figures \ref{Fig_3D_C1}-\ref{Fig_u_C1} show three-dimensional trajectory, error history, and control history, respectively, for a single random seed. The results shown are obtained i) with no correction, ii) with parameter correction, and iii) with control function correction. Figure \ref{Fig_e_rvm_f_C1} indicates that the two incremental correction methods exhibit substantially different characteristics in their performance in terms of the capability to reduce final position and/or velocity errors. The final errors are neither completely nullified nor always reduced by performing policy correction once at the initial time. According to Fig. \ref{Fig_e_rvm_f_C1} and additional numerical experiments, the parameter correction method yields significant reduction in the final position error $e_{r_{f}} \triangleq \left\|\mathbf{r}\left(t_{f}\right)-\mathbf{r}_{f_{d}}\right\|$ while showing less effectiveness in reducing the final velocity error $e_{v_{f}}\triangleq \left\|\mathbf{v}\left(t_{f}\right)-\mathbf{v}_{f_{d}}\right\|$. On the other hand, the control function correction method provides improved final velocity targeting accuracy with increased fuel expenditure, however, the final position accuracy is degraded even in comparison to the case with baseline policy alone. Also, a noticeable rapid change in the control input is often observed around the end of flight as shown in the history of $\sigma_{T}$ in Fig. \ref{Fig_u_C1}. 

In summary, different incremental correction algorithms employed in Stage 2 lead to considerable differences in the satisfaction of constraints on the performance output. Although the trend understandably depends on design parameters being used in each method, the numerical experiments suggest that the parameter correction method tends to be more consistent and thus reliable in satisfying the performance output constraints.

\clearpage

\subsubsection{Case 2. Dispersion in Initial Conditions} \label{SubSubSec:Sim_Case2}
\paragraph{Stage 1. Baseline Policy Training}
The initial and final conditions considered in Case 1 are taken as the nominal conditions for Case 2, and the baseline policy trained in Case 1 for $t_{f}=43$s is re-used in Case 2.

\paragraph{Stage 2. Incremental Correction with Fixed $t_{f}$ and Various Random Seeds}
In a similar way as done in Case 1, both correction techniques are applied to update the baseline policies that are trained by using different random seeds for NN parameter initialisation. The initial position perturbation is modelled deterministically as
\begin{equation} \label{Eq:init_pos_dispersion}
	\bar{\mathbf{r}}_{0} = \mathbf{r}_{0} + \bar{r}\left(\cos\alpha  \frac{\mathbf{v}_{0} \times \mathbf{r}_{0}}{\left\| \mathbf{v}_{0} \times \mathbf{r}_{0}\right\|} \times \hat{\mathbf{v}}_{0} + \sin\alpha \frac{\mathbf{v}_{0} \times \mathbf{r}_{0}}{\left\| \mathbf{v}_{0} \times \mathbf{r}_{0}\right\|}\right)
\end{equation}
where $\alpha$ is the angle introduced to parametrise the circle of radius $\bar{r} = 100$m centred at $\mathbf{r}_{0}$ on the plane perpendicular to $\hat{\mathbf{v}}_{0}$. 16 equally spaced points are obtained by considering a range of $\alpha$ values in $\left[0, 2\pi\right]$ with the interval of $\pi/8\,$rad.

\begin{figure}[!htbp]
	\begin{center}
		\begin{subfigure}[ht!]{0.32\textwidth}
			\includegraphics[width=\hsize]{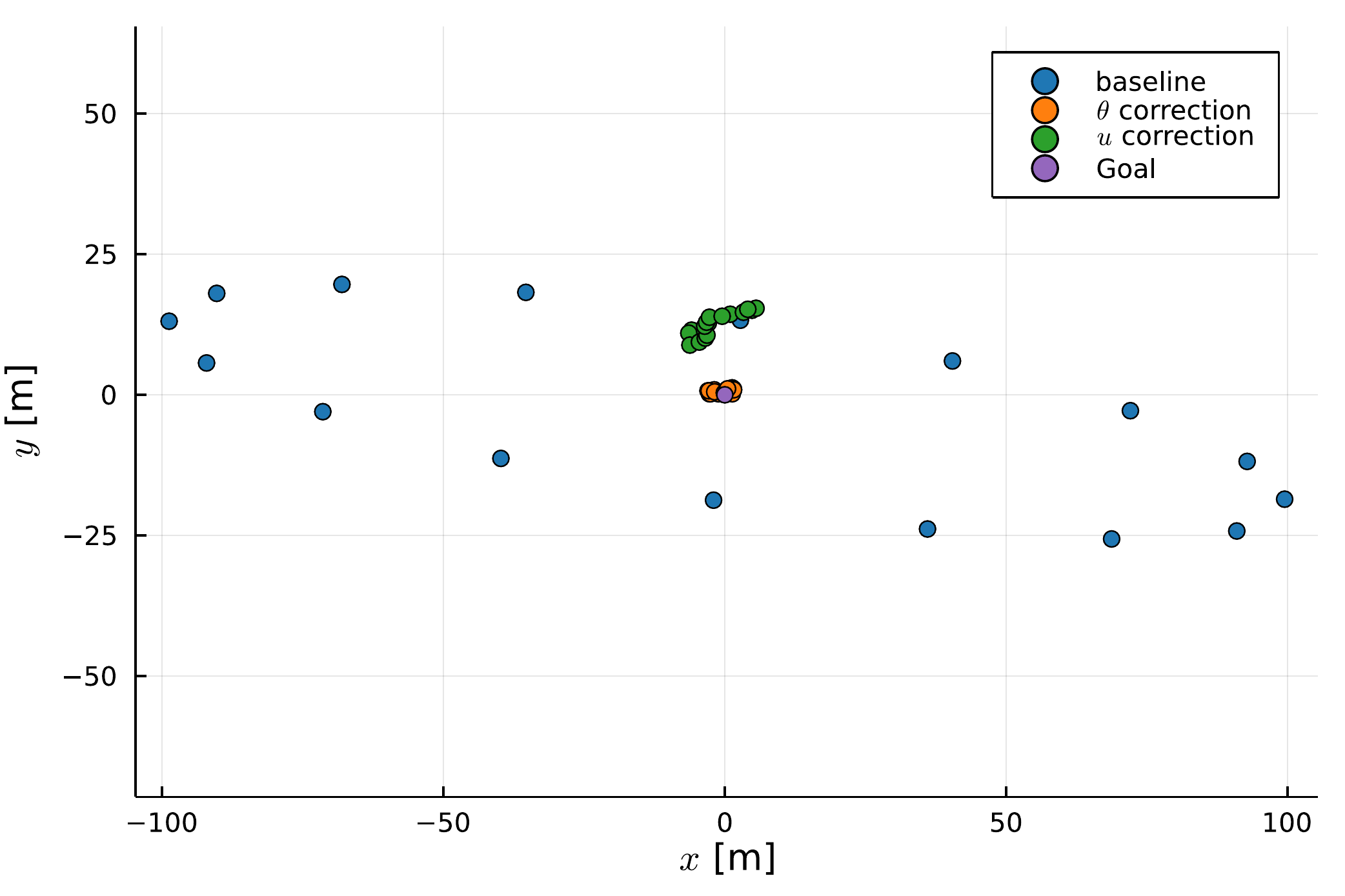}
			\caption{Random Seed $0$}
			\label{Fig_2D_f_C2_1}
		\end{subfigure}
		~
		\begin{subfigure}[ht!]{0.32\textwidth}
			\includegraphics[width=\hsize]{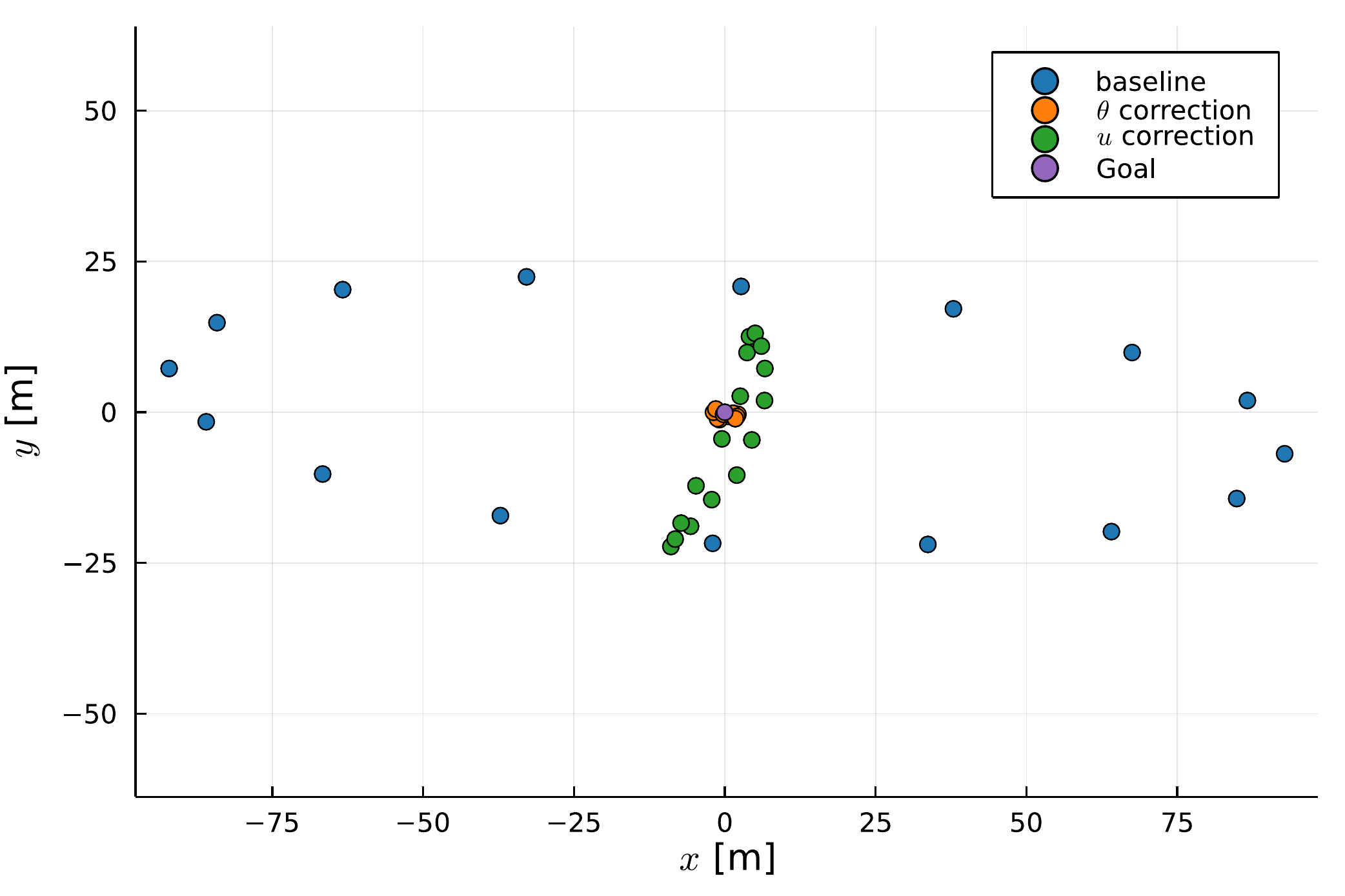}
			\caption{Random Seed $1$}
			\label{Fig_2D_f_C2_2}
		\end{subfigure}
		~
		\begin{subfigure}[ht!]{0.32\textwidth}
			\includegraphics[width=\hsize]{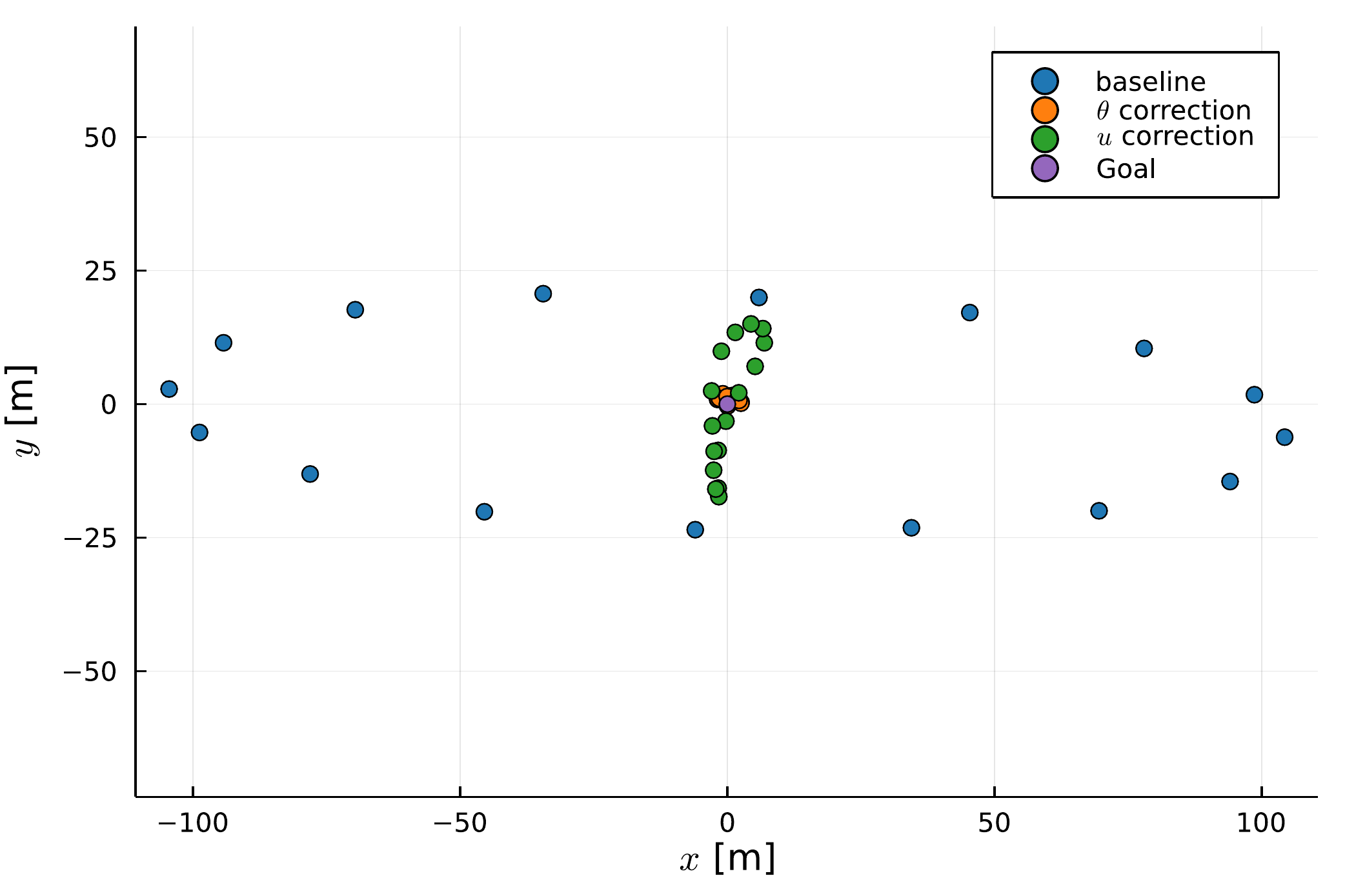}
			\caption{Random Seed $2$}
			\label{Fig_2D_f_C2_3}
		\end{subfigure}
		\vspace{1em}
		
		\begin{subfigure}[ht!]{0.32\textwidth}
			\includegraphics[width=\hsize]{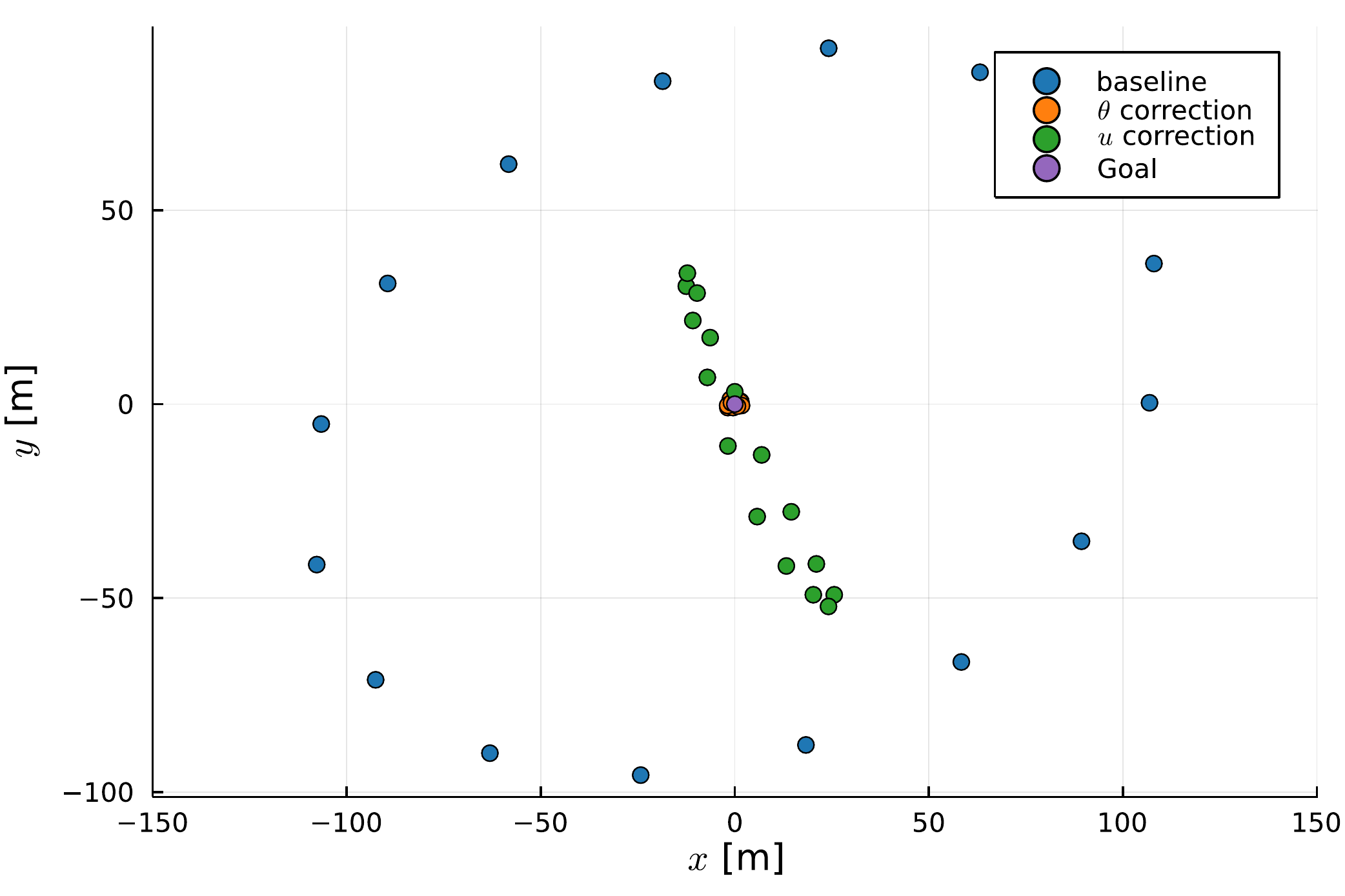}
			\caption{Random Seed $3$}
			\label{Fig_2D_f_C2_4}
		\end{subfigure}
		~
		\begin{subfigure}[ht!]{0.32\textwidth}
			\includegraphics[width=\hsize]{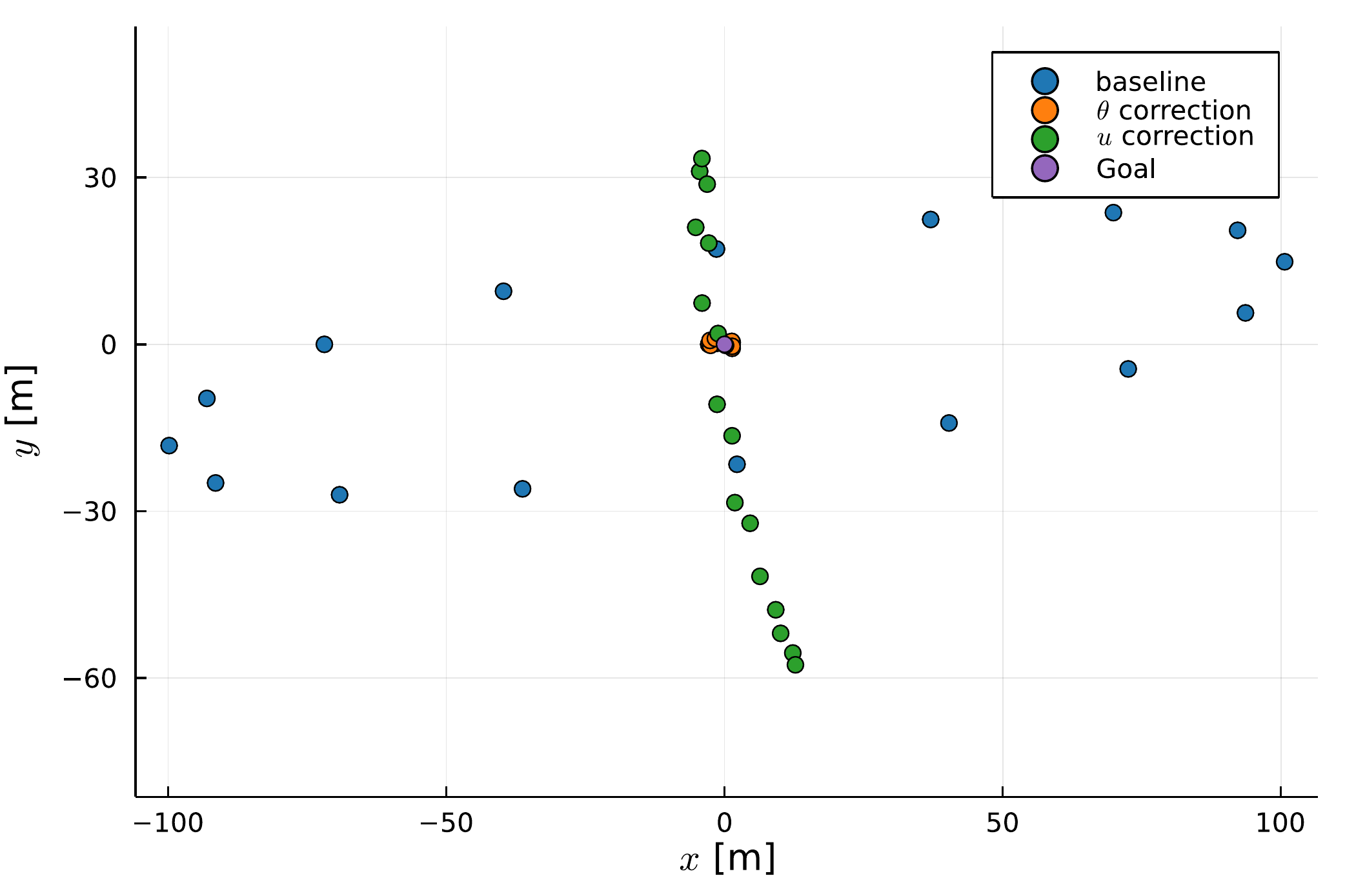}
			\caption{Random Seed $4$}
			\label{Fig_2D_f_C2_5}
		\end{subfigure}
		~
		\begin{subfigure}[ht!]{0.32\textwidth}
			\includegraphics[width=\hsize]{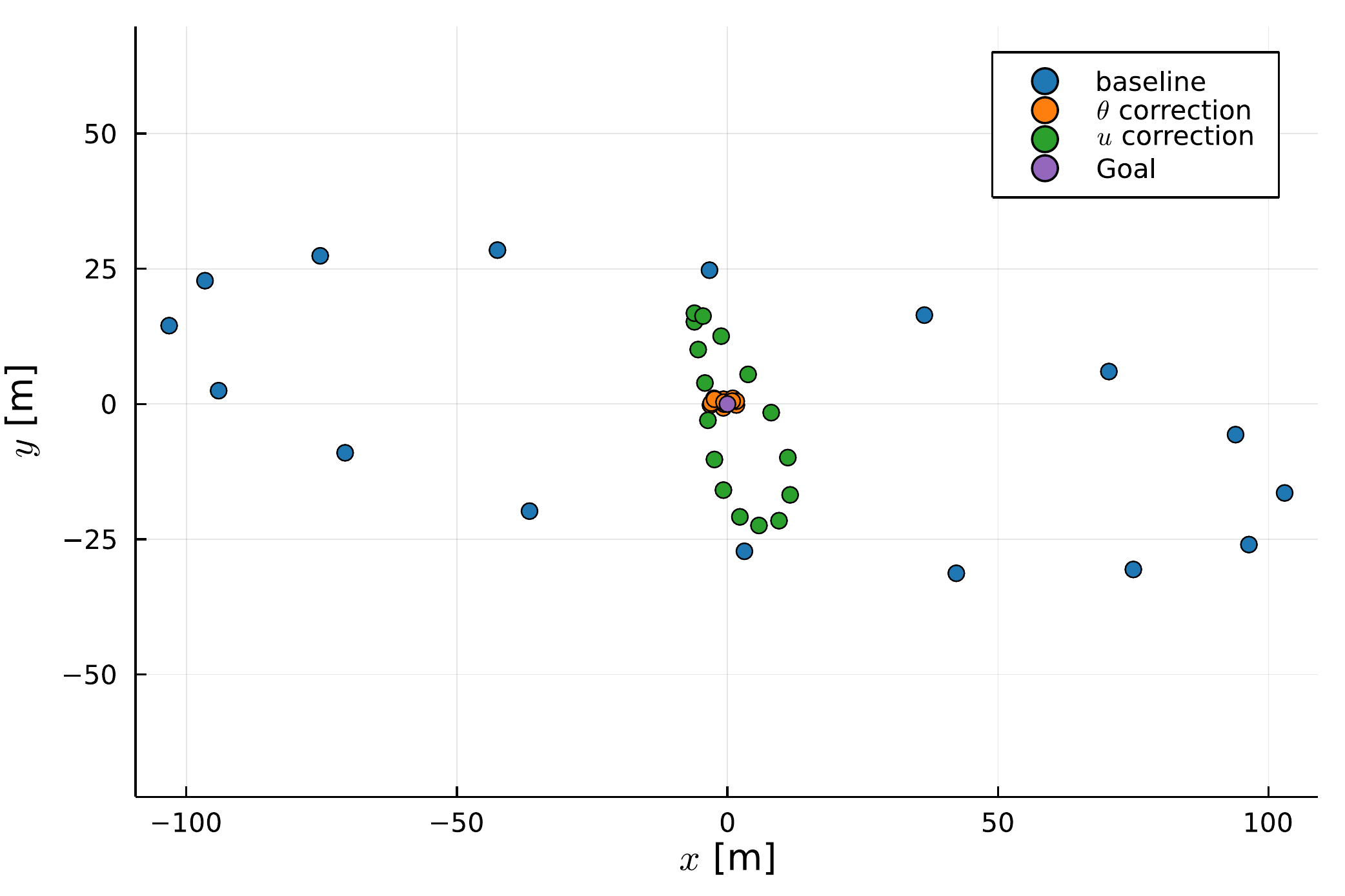}
			\caption{Random Seed $5$}
			\label{Fig_2D_f_C2_6}
		\end{subfigure}
		\vspace{1em}
		
		\begin{subfigure}[ht!]{0.32\textwidth}
			\includegraphics[width=\hsize]{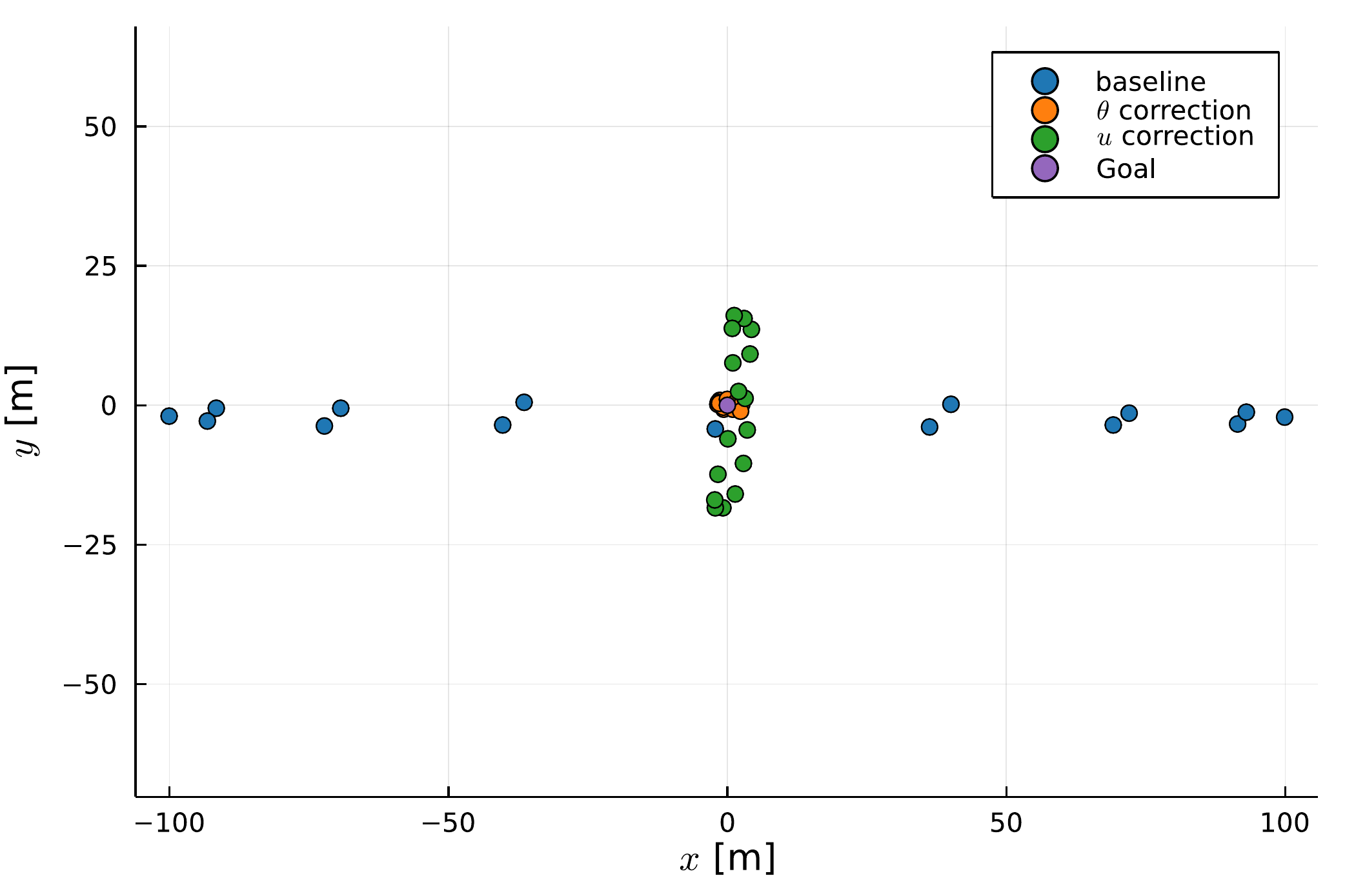}
			\caption{Random Seed $6$}
			\label{Fig_2D_f_C2_7}
		\end{subfigure}
		~
		\begin{subfigure}[ht!]{0.32\textwidth}
			\includegraphics[width=\hsize]{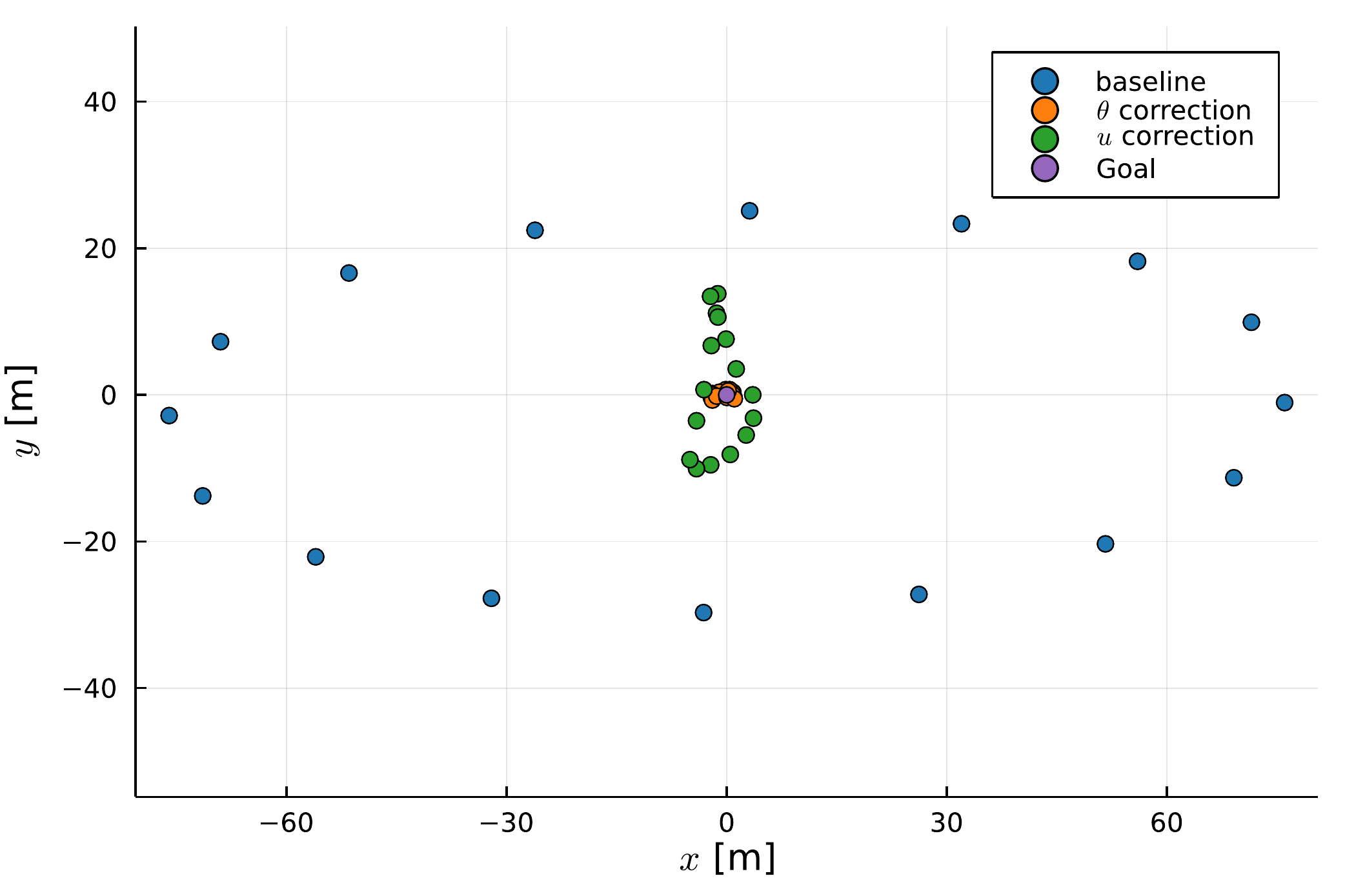}
			\caption{Random Seed $7$}
			\label{Fig_2D_f_C2_8}
		\end{subfigure}
		~
		\begin{subfigure}[ht!]{0.32\textwidth}
			\includegraphics[width=\hsize]{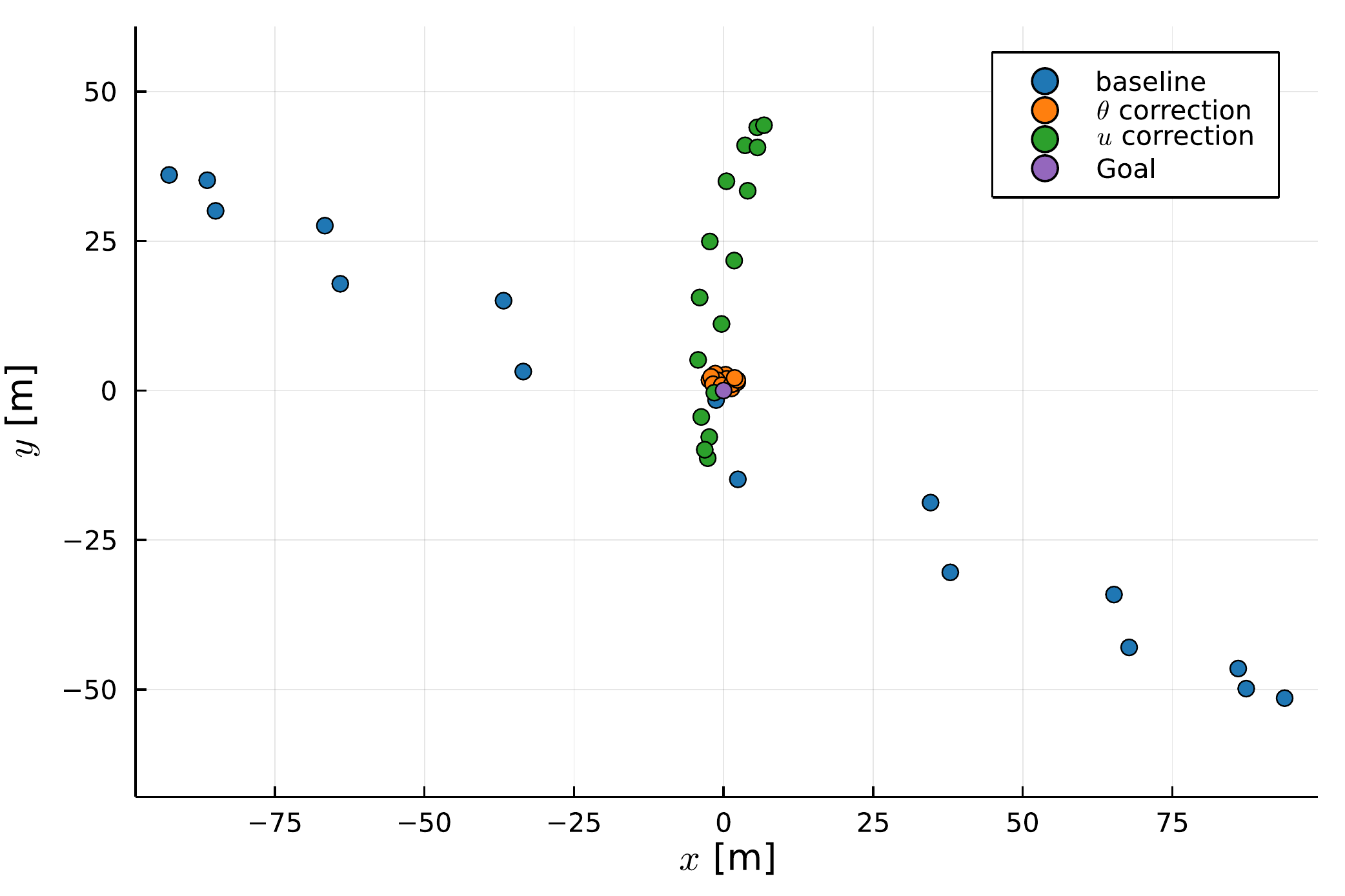}
			\caption{Random Seed $8$}
			\label{Fig_2D_f_C2_9}
		\end{subfigure}
		\vspace{1em}
		
		\begin{subfigure}[ht!]{0.32\textwidth}
			\includegraphics[width=\hsize]{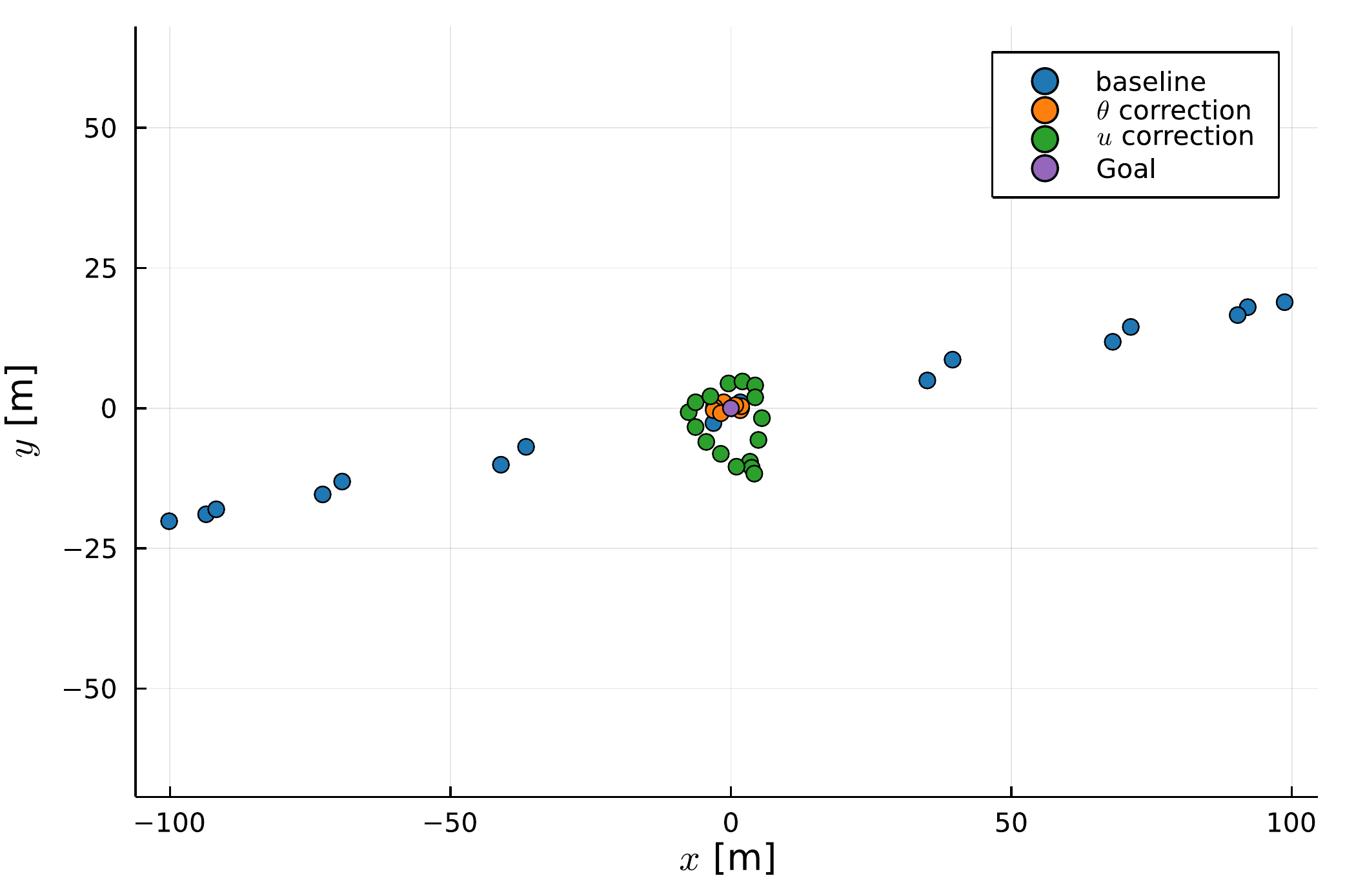}
			\caption{Random Seed $9$}
			\label{Fig_2D_f_C2_10}
		\end{subfigure}
		~
		\begin{subfigure}[ht!]{0.32\textwidth}
			\includegraphics[width=\hsize]{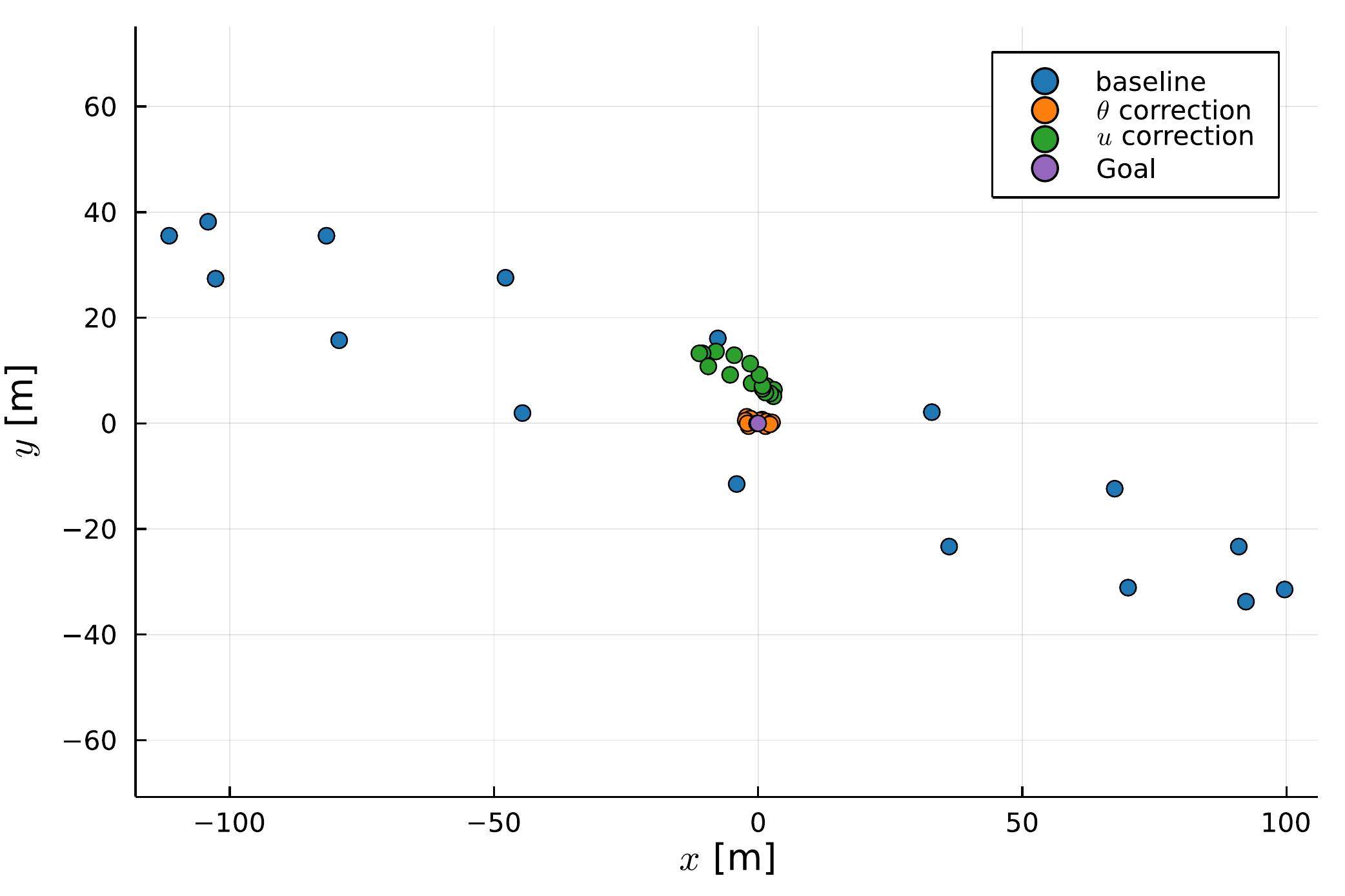}
			\caption{Random Seed $10$}
			\label{Fig_2D_f_C2_11}
		\end{subfigure}
		
		\caption{Final Position on Horizontal Plane (Case 2)}
		\label{Fig_2D_f_C2}
	\end{center}
\end{figure}

Each subfigure in Fig. \ref{Fig_2D_f_C2} shows the final position projected on the horizontal plane attached to the desired landing position for all combinations of initial condition and incremental correction method, with a fixed baseline policy. In Fig. \ref{Fig_2D_f_C2}, the clusters of points are scattered around the desired landing position with different size and shape depending on the correction method. The size of the dispersion in terms of the area contained within the convex hull over the trace of final positions is the least with the parameter correction method regardless of the random seed. The centroid of the dispersion resulting from parameter correction is also well-aligned with the goal position for all random seeds. The control function correction method could also reduce the size of the landing point cluster as compared to the case with no correction, however, the irregularly-shaped trace is not always centred around the desired landing position.

\begin{figure}[!htbp]
	\vspace{-2cm}
	\begin{center}
		\begin{subfigure}[ht!]{0.32\textwidth}
			\includegraphics[width=\hsize]{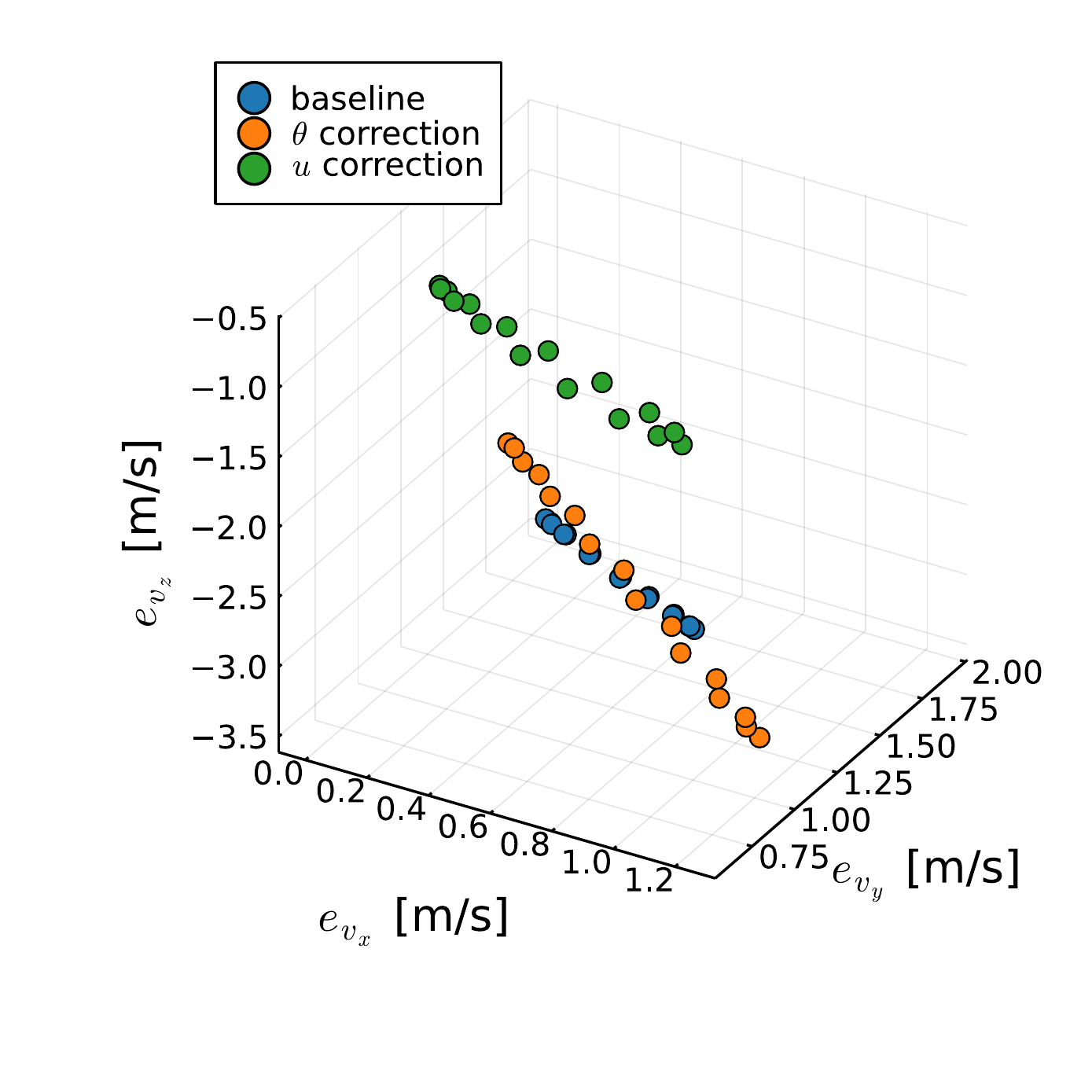}
			\caption{Random Seed $0$}
			\label{Fig_3D_e_v_C2_1}
		\end{subfigure}
		~
		\begin{subfigure}[ht!]{0.32\textwidth}
			\includegraphics[width=\hsize]{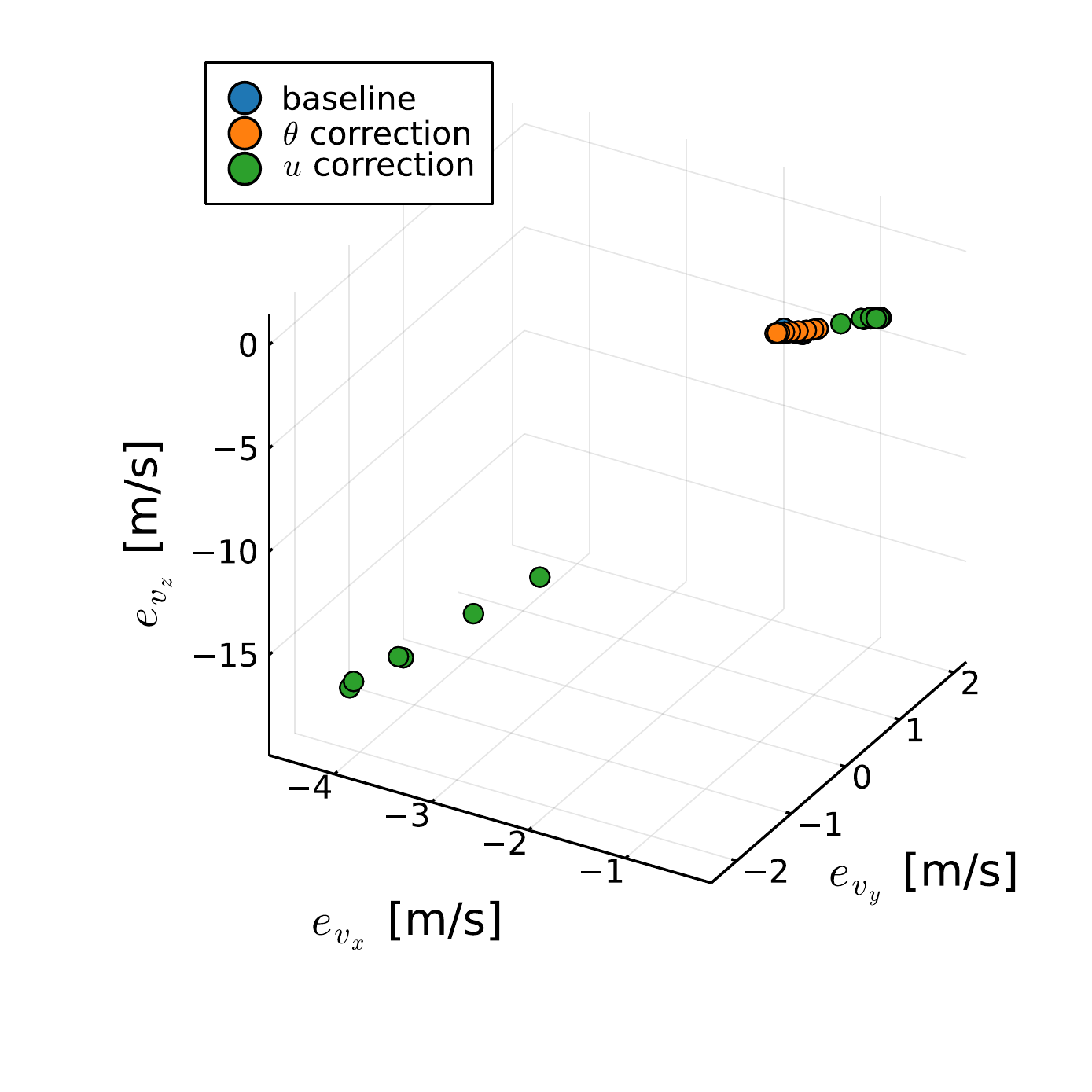}
			\caption{Random Seed $1$}
			\label{Fig_3D_e_v_C2_2}
		\end{subfigure}
		~
		\begin{subfigure}[ht!]{0.32\textwidth}
			\includegraphics[width=\hsize]{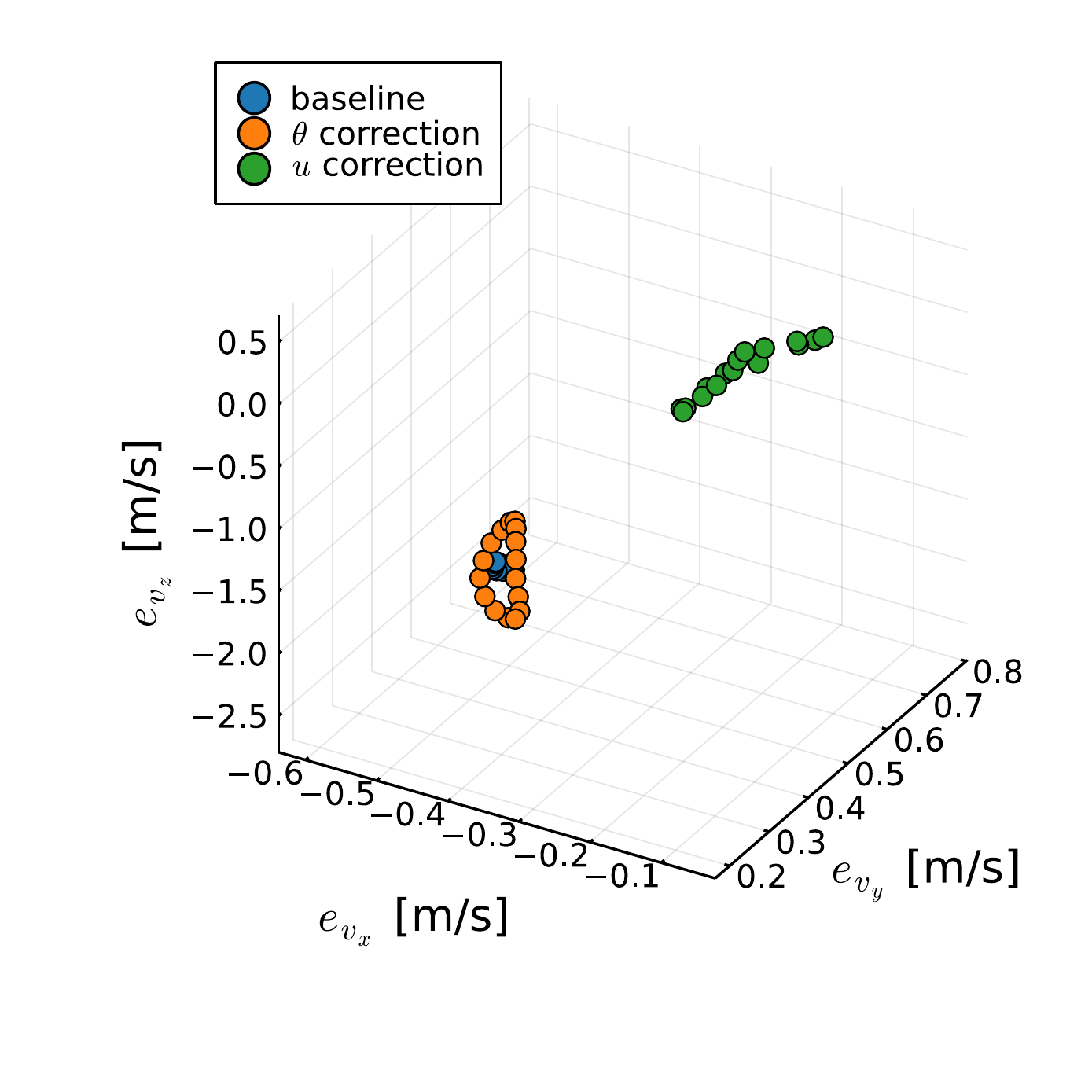}
			\caption{Random Seed $2$}
			\label{Fig_3D_e_v_C2_3}
		\end{subfigure}
		\vspace{1em}
		
		\begin{subfigure}[ht!]{0.32\textwidth}
			\includegraphics[width=\hsize]{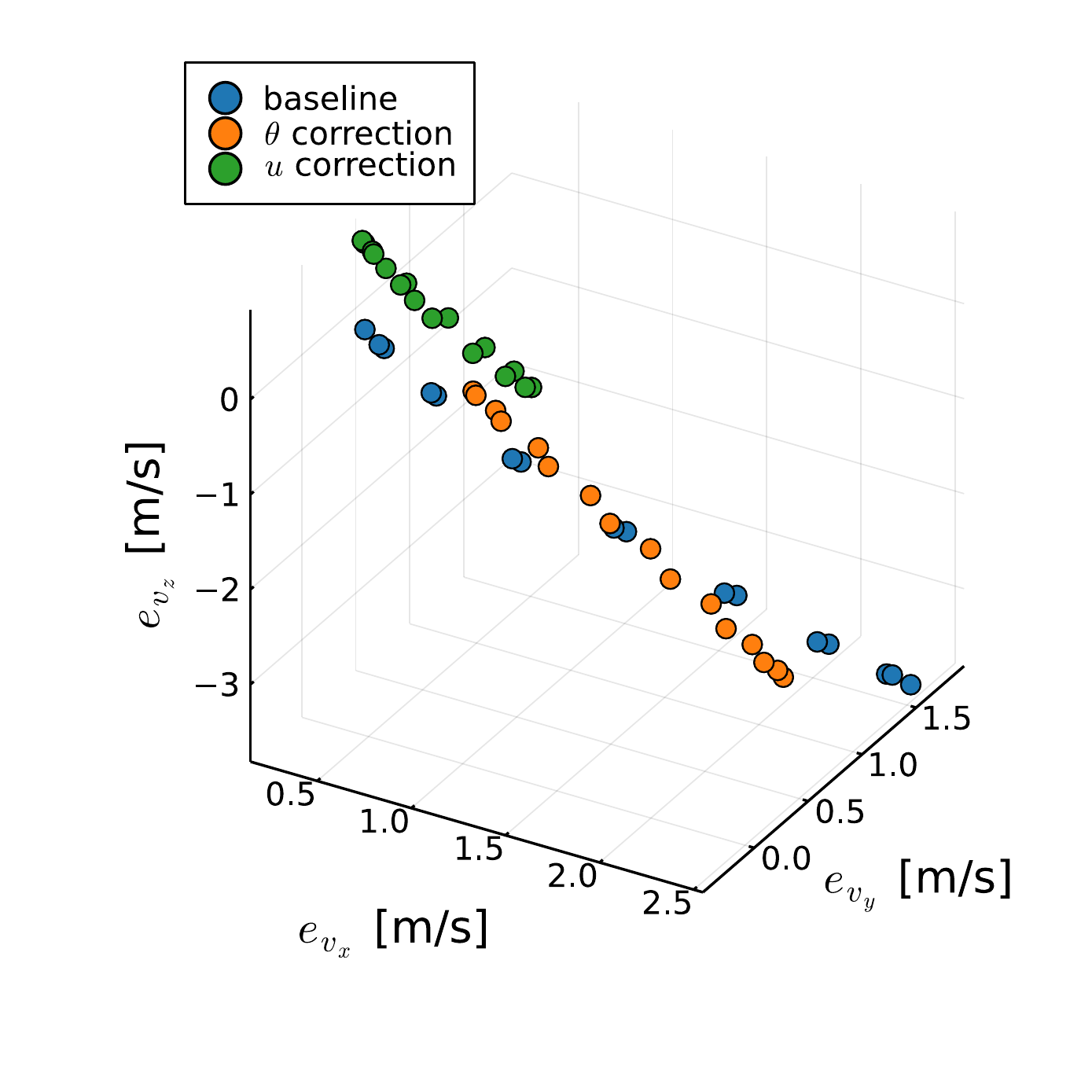}
			\caption{Random Seed $3$}
			\label{Fig_3D_e_v_C2_4}
		\end{subfigure}
		~
		\begin{subfigure}[ht!]{0.32\textwidth}
			\includegraphics[width=\hsize]{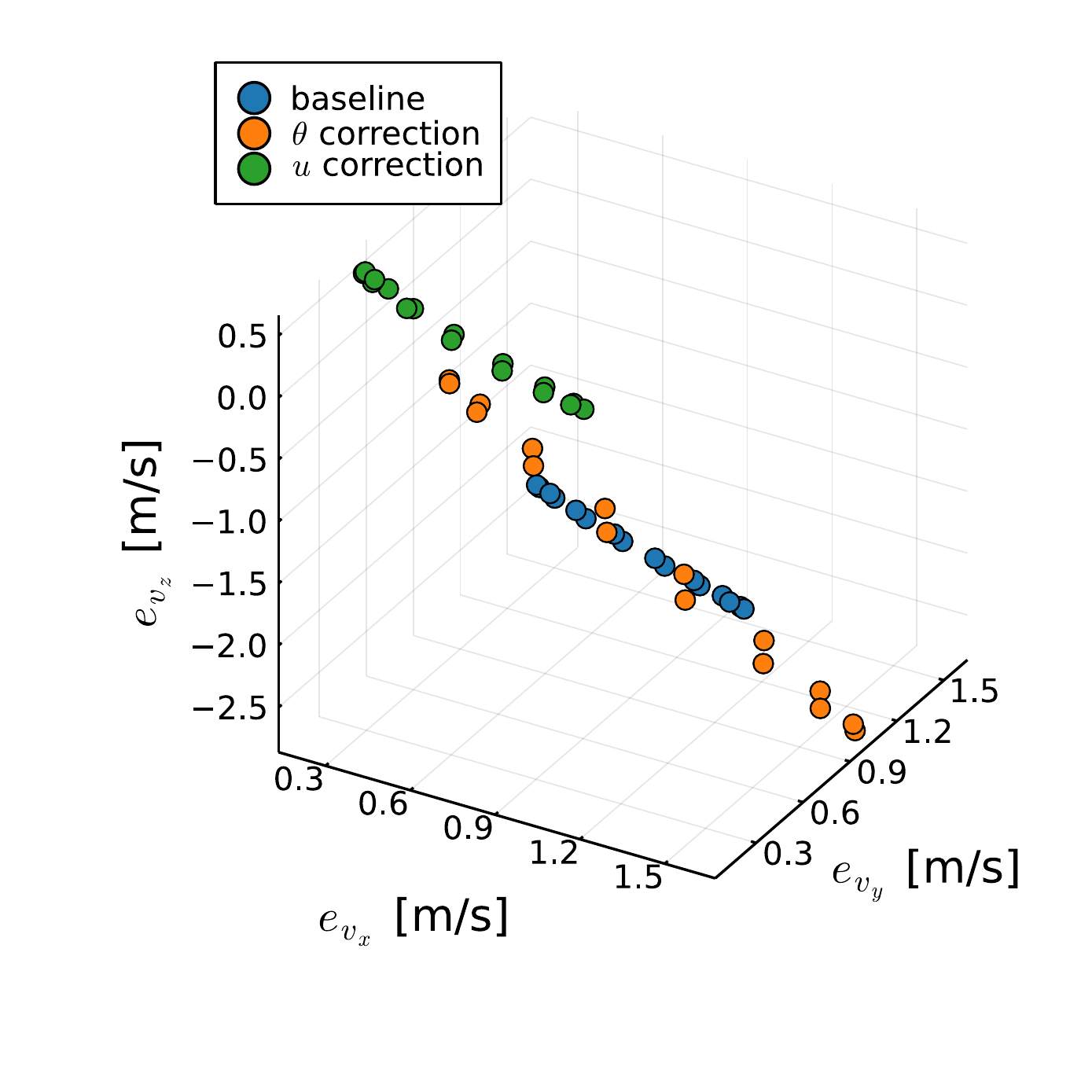}
			\caption{Random Seed $4$}
			\label{Fig_3D_e_v_C2_5}
		\end{subfigure}
		~
		\begin{subfigure}[ht!]{0.32\textwidth}
			\includegraphics[width=\hsize]{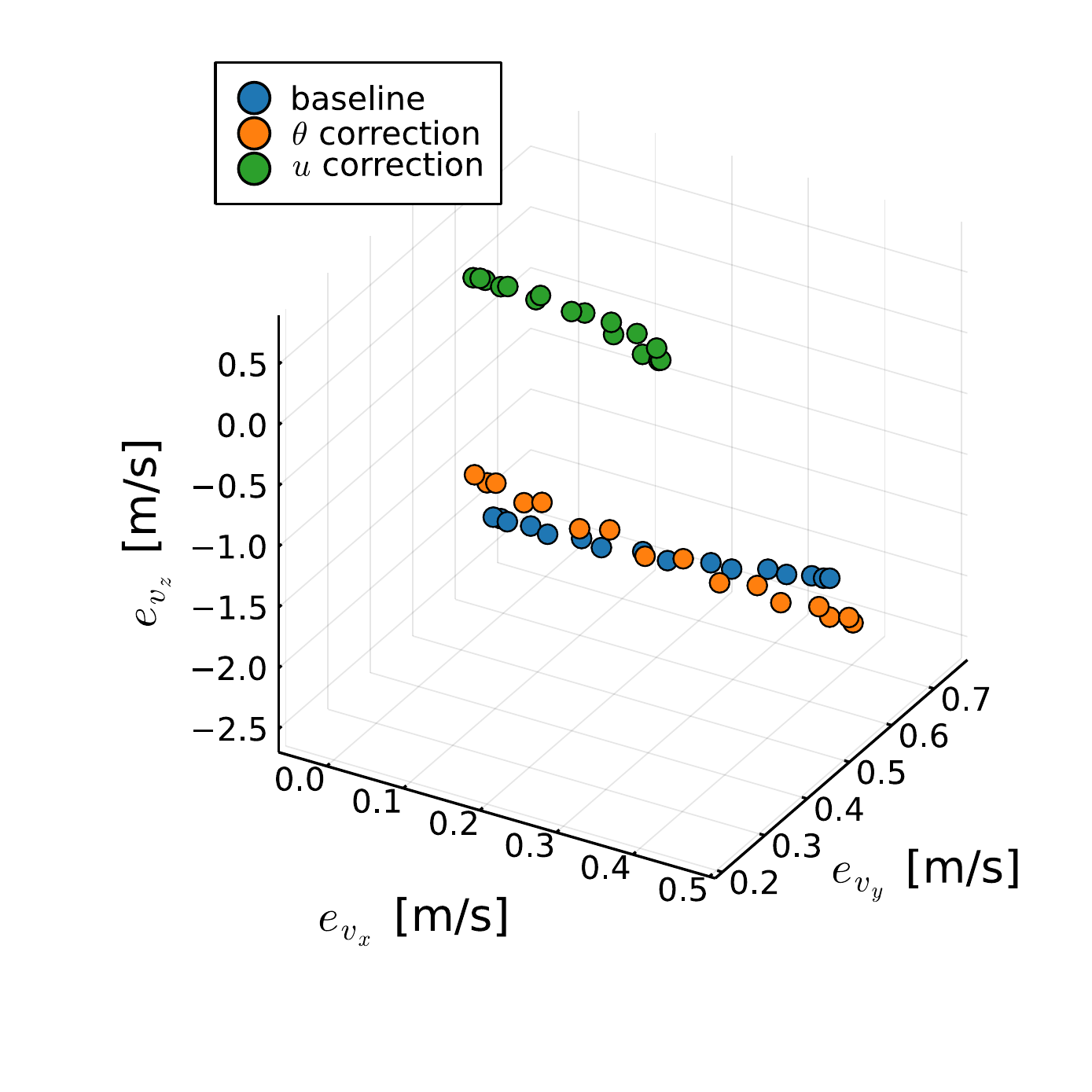}
			\caption{Random Seed $5$}
			\label{Fig_3D_e_v_C2_6}
		\end{subfigure}
		\vspace{1em}
		
		\begin{subfigure}[ht!]{0.32\textwidth}
			\includegraphics[width=\hsize]{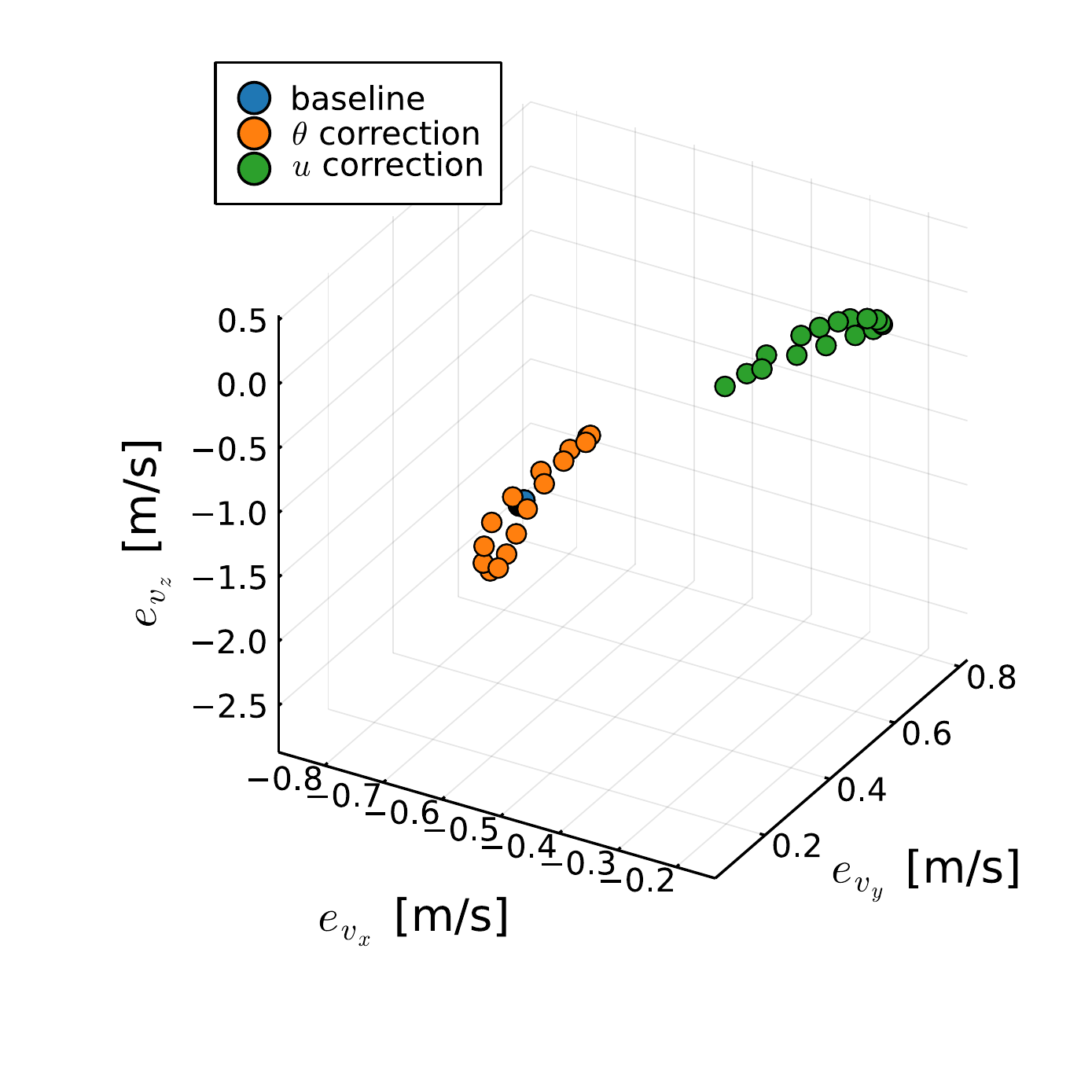}
			\caption{Random Seed $6$}
			\label{Fig_3D_e_v_C2_7}
		\end{subfigure}
		~
		\begin{subfigure}[ht!]{0.32\textwidth}
			\includegraphics[width=\hsize]{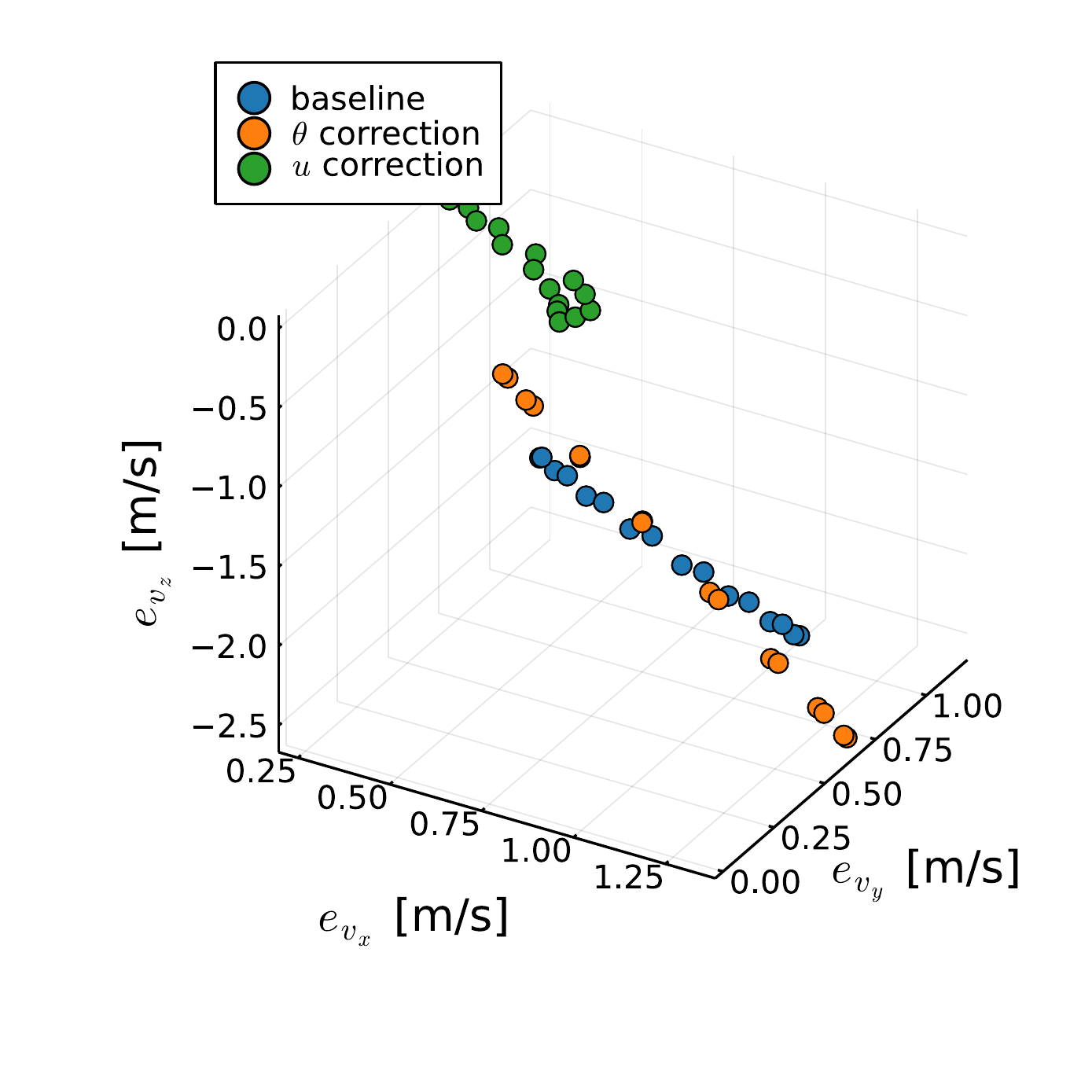}
			\caption{Random Seed $7$}
			\label{Fig_3D_e_v_C2_8}
		\end{subfigure}
		~
		\begin{subfigure}[ht!]{0.32\textwidth}
			\includegraphics[width=\hsize]{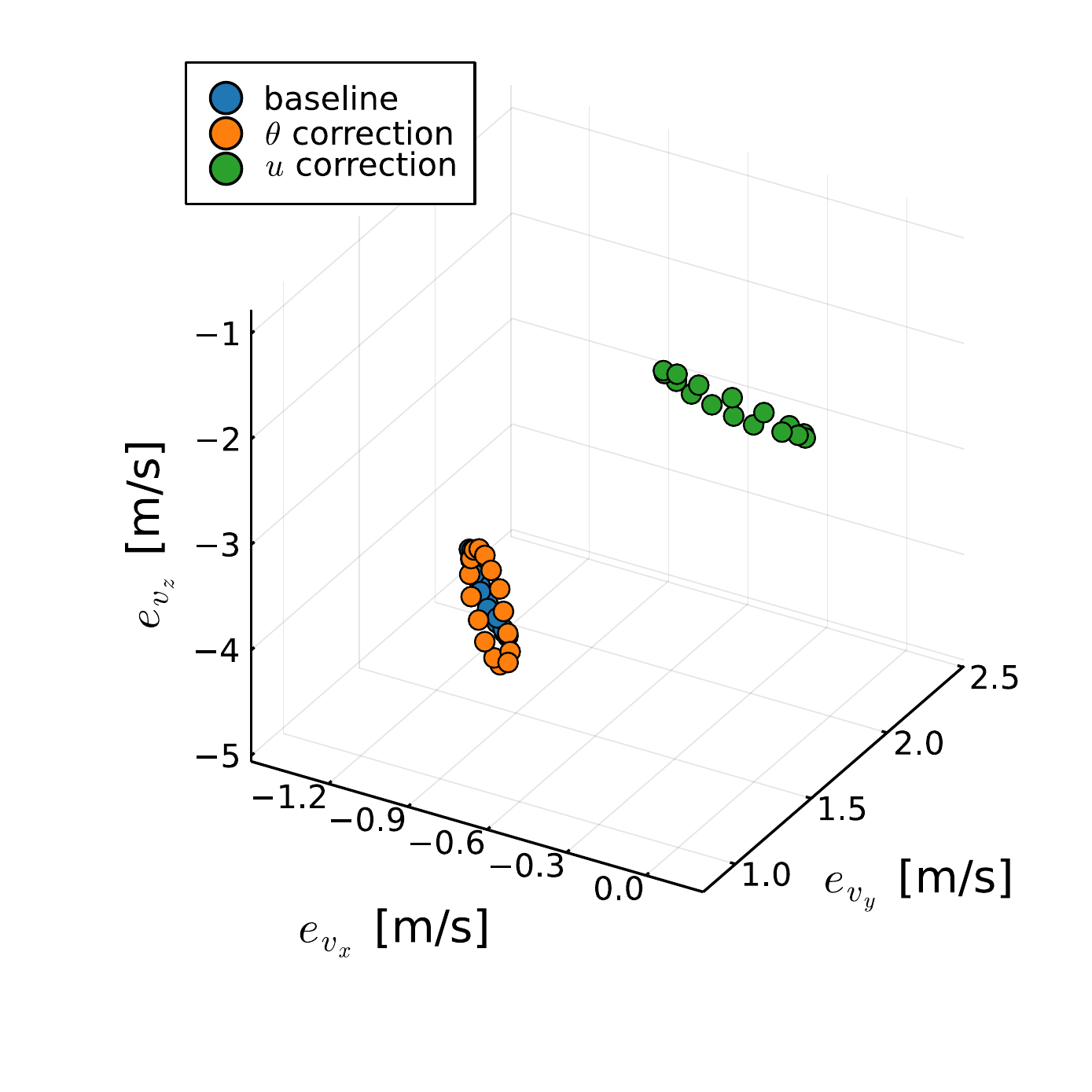}
			\caption{Random Seed $8$}
			\label{Fig_3D_e_v_C2_9}
		\end{subfigure}
		\vspace{1em}
		
		\begin{subfigure}[ht!]{0.32\textwidth}
			\includegraphics[width=\hsize]{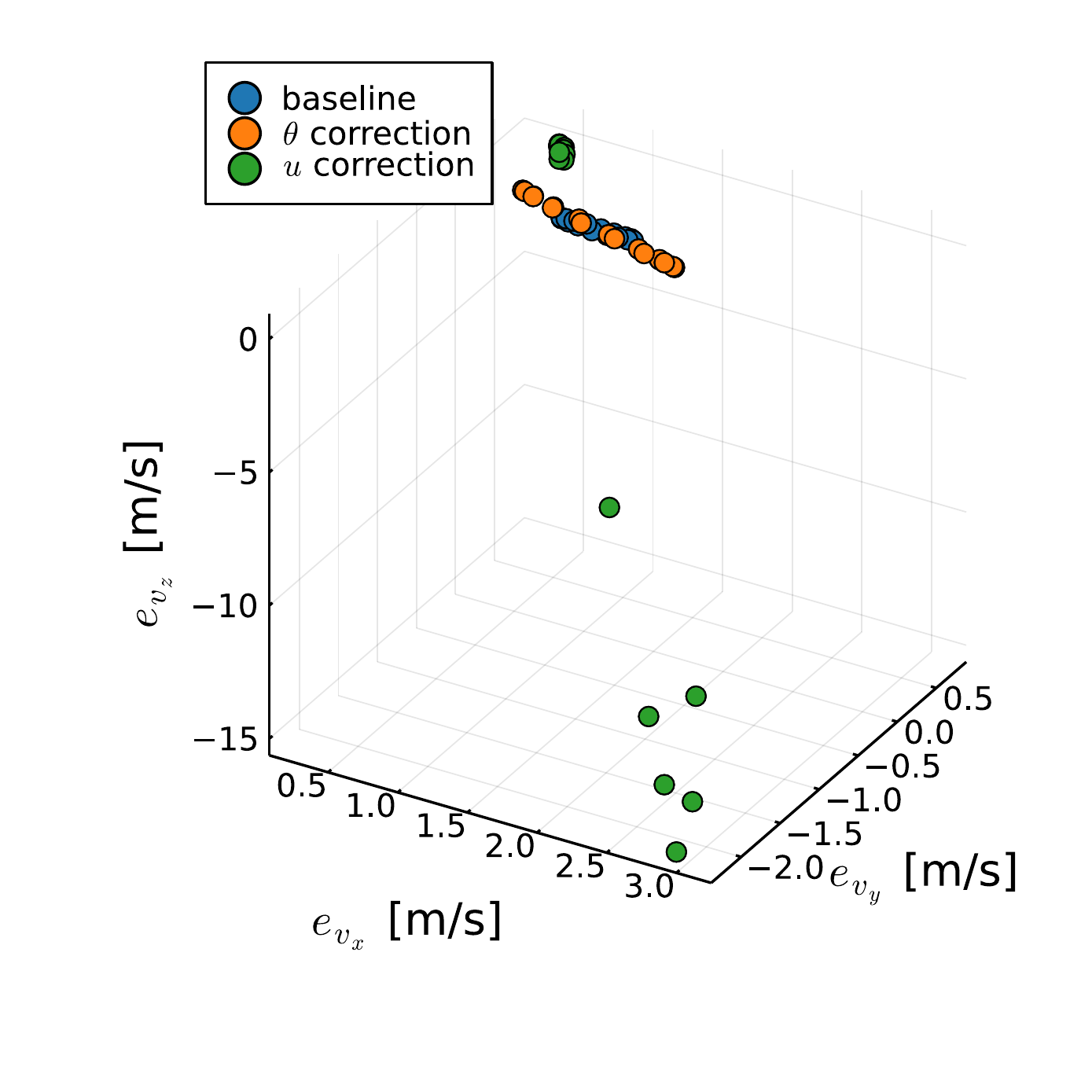}
			\caption{Random Seed $9$}
			\label{Fig_3D_e_v_C2_10}
		\end{subfigure}
		~
		\begin{subfigure}[ht!]{0.32\textwidth}
			\includegraphics[width=\hsize]{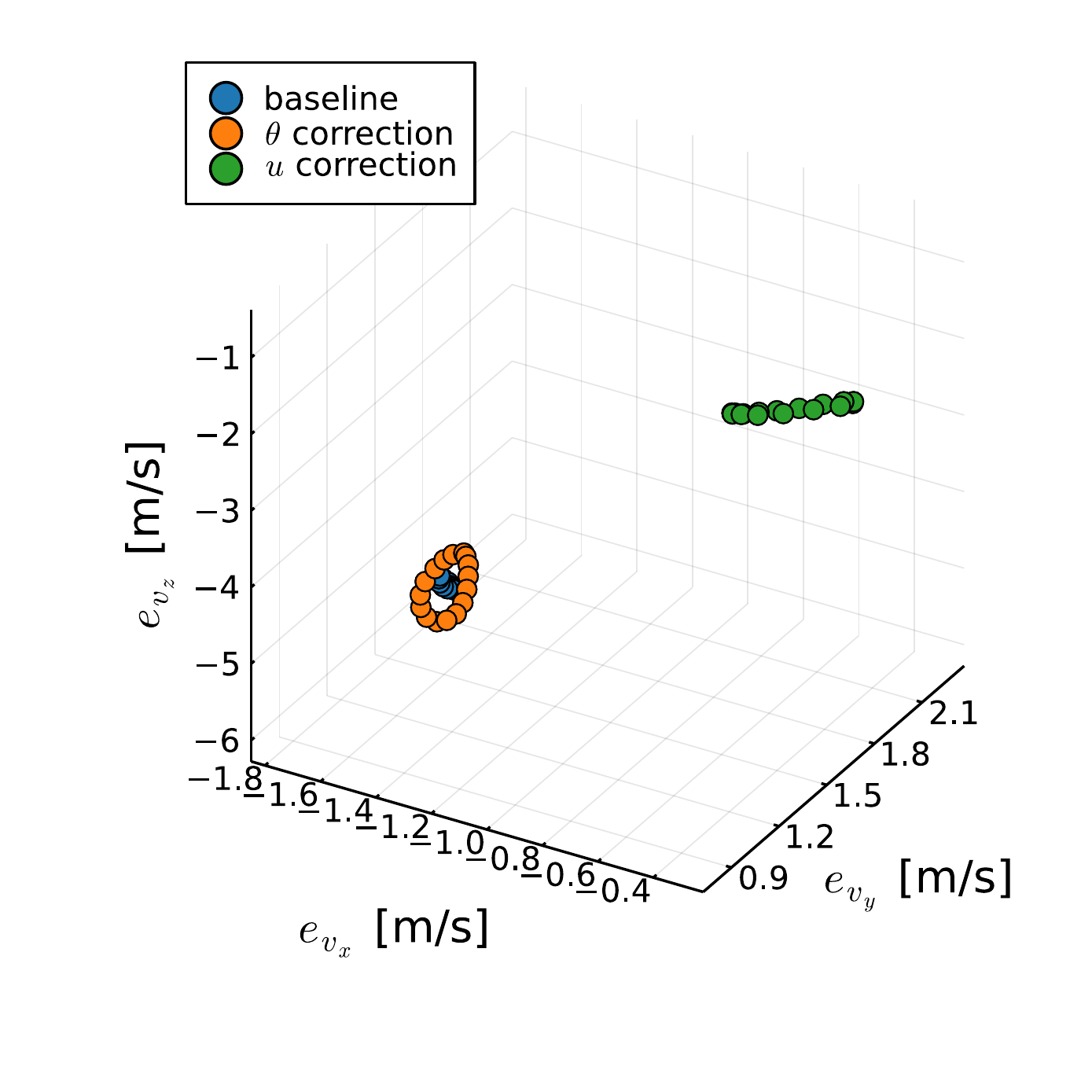}
			\caption{Random Seed $10$}
			\label{Fig_3D_e_v_C2_11}
		\end{subfigure}
		
		\caption{Final Velocity Error (Case 2)}
		\label{Fig_3D_e_v_C2}
	\end{center}
\end{figure}

Likewise, Fig. \ref{Fig_3D_e_v_C2} shows the final velocity error represented in the local horizontal and vertical coordinate system with respect to the desired landing position for all combination of baseline policy, initial condition, and incremental correction method. Here, the final velocity error is defined by the difference between the achieved and the desired final velocity vectors, i.e., $\mathbf{v}\left(t_{f}\right) - \mathbf{v}_{f_{d}}$. The trend observed in the final velocity error with control function correction is inconsistent across different random seeds and tends to show high variance depending on initial position in each random seed. The irregular tendency manifests as large magnitude of final velocity error in some cases. As for the parameter correction method, the improvement in constraint satisfaction accuracy is relatively limited for the final velocity as compared to that for the final position. However, the distribution of final velocity error resulting from parameter correction does not involve serious outliers that may pose high risk at the final time.

\begin{figure}[!ht]
	\begin{center}
		\includegraphics[width=0.7\textwidth]{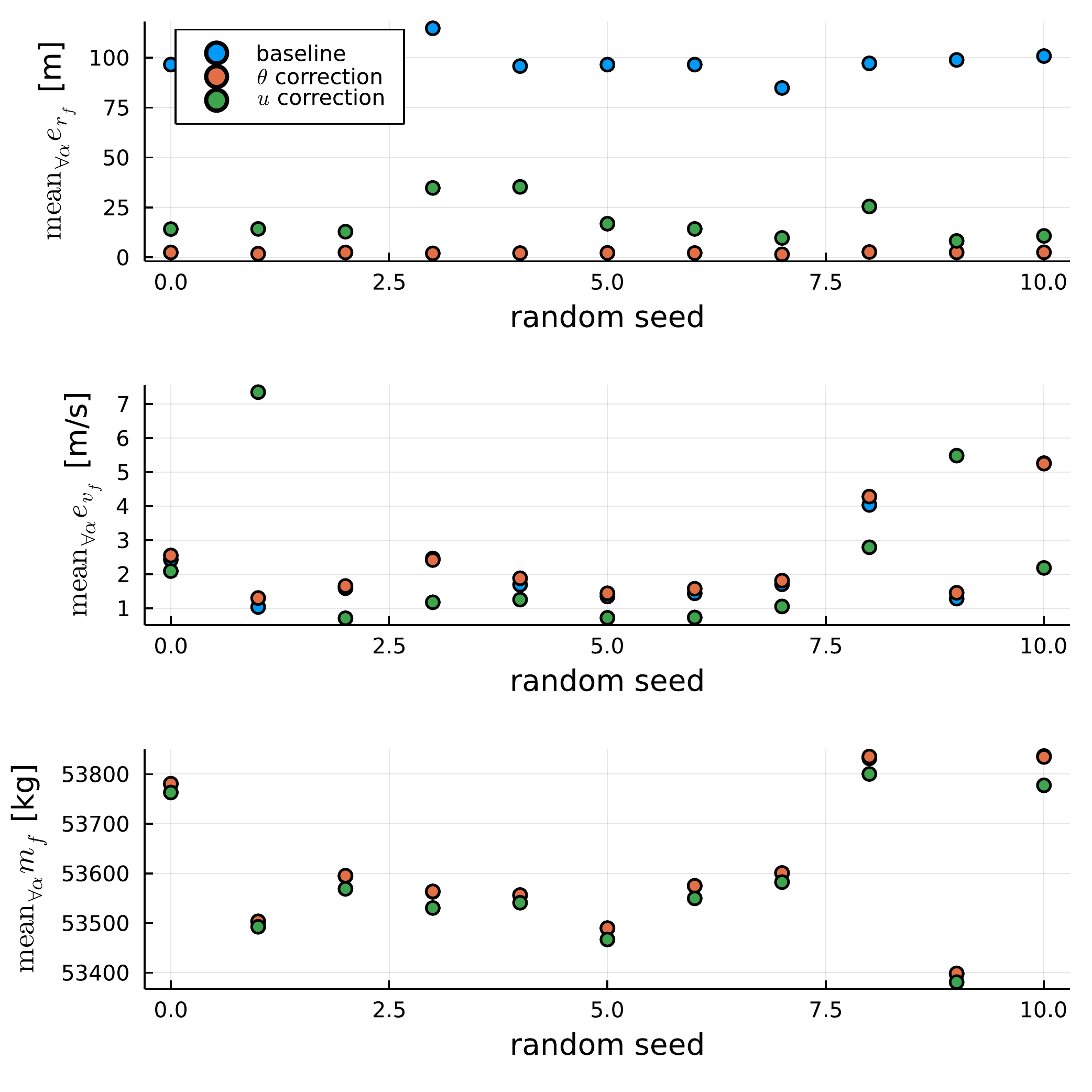}
		\caption{Mean of Final Errors and Mass with Incremental Correction (Case 2)} \label{Fig_mean_e_rvm_f_C2}
	\end{center}
\end{figure}

\begin{figure}[!ht]
	\begin{center}
		\includegraphics[width=0.7\textwidth]{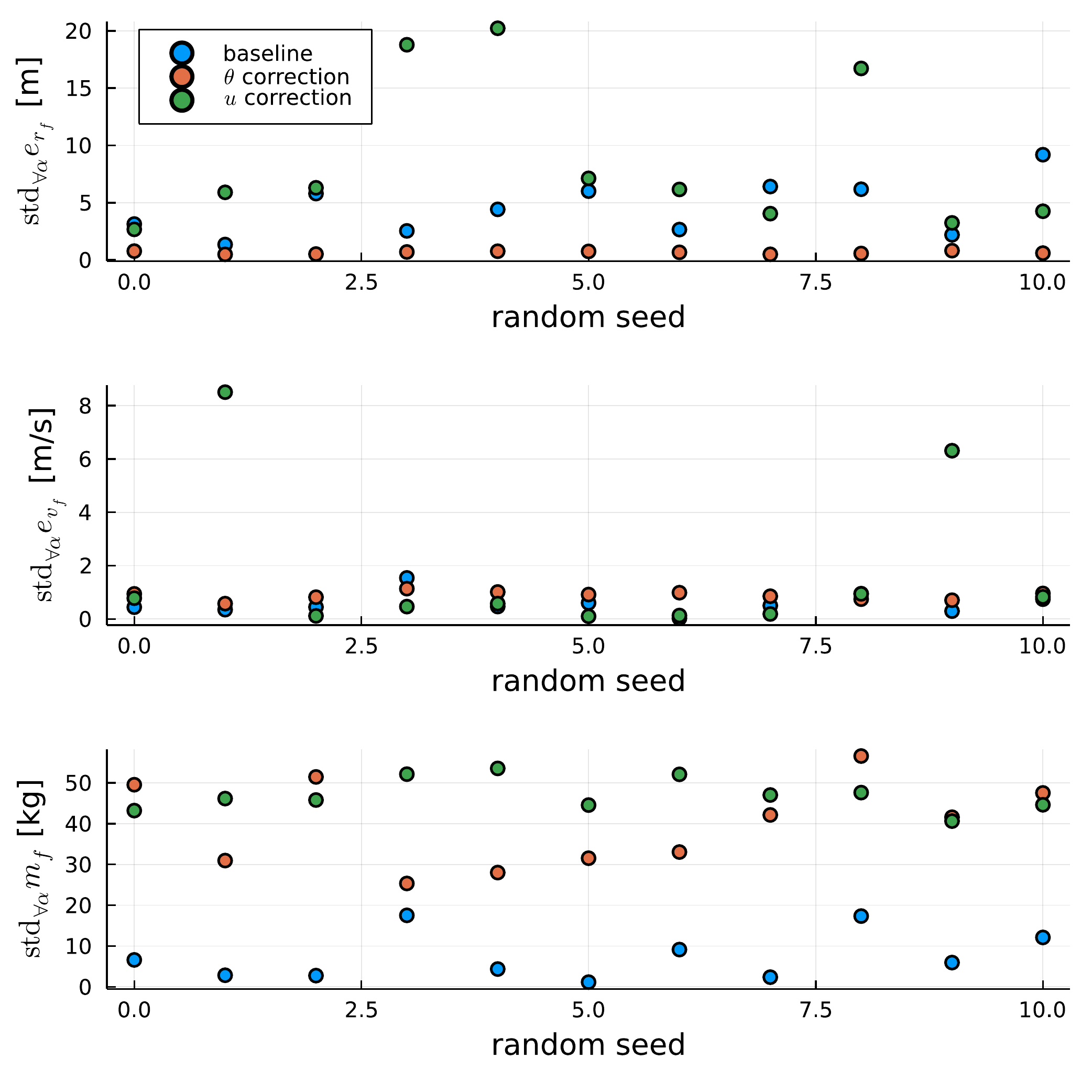}
		\caption{Standard Deviation of Final Errors and Mass with Incremental Correction (Case 2)} \label{Fig_std_e_rvm_f_C2}
	\end{center}
\end{figure}

To facilitate more quantitative analysis of errors, Figs. \ref{Fig_mean_e_rvm_f_C2} and \ref{Fig_std_e_rvm_f_C2} show the statistics of final position and velocity errors, and final mass in terms of the mean and standard deviation over the range of initial conditions for each random seed, respectively. The mean and standard deviation depicted in Figs. \ref{Fig_mean_e_rvm_f_C2} and \ref{Fig_std_e_rvm_f_C2} are consistent with the observation from Fig. \ref{Fig_2D_f_C2}. The accuracy and precision of landing position represented by the mean and the standard deviation of $e_{r_{f}}$, respectively, are improved most significantly by the parameter correction method. The performance benefit of the control function correction method in final position is not entirely clear as it might result in larger standard deviation in $e_{r_{f}}$. The final velocity error reduction is not apparent with the parameter correction as in Case 1.

\begin{figure}[!ht]
	\vspace{-0.5cm}
	\begin{center}
		\includegraphics[width=0.6\textwidth]{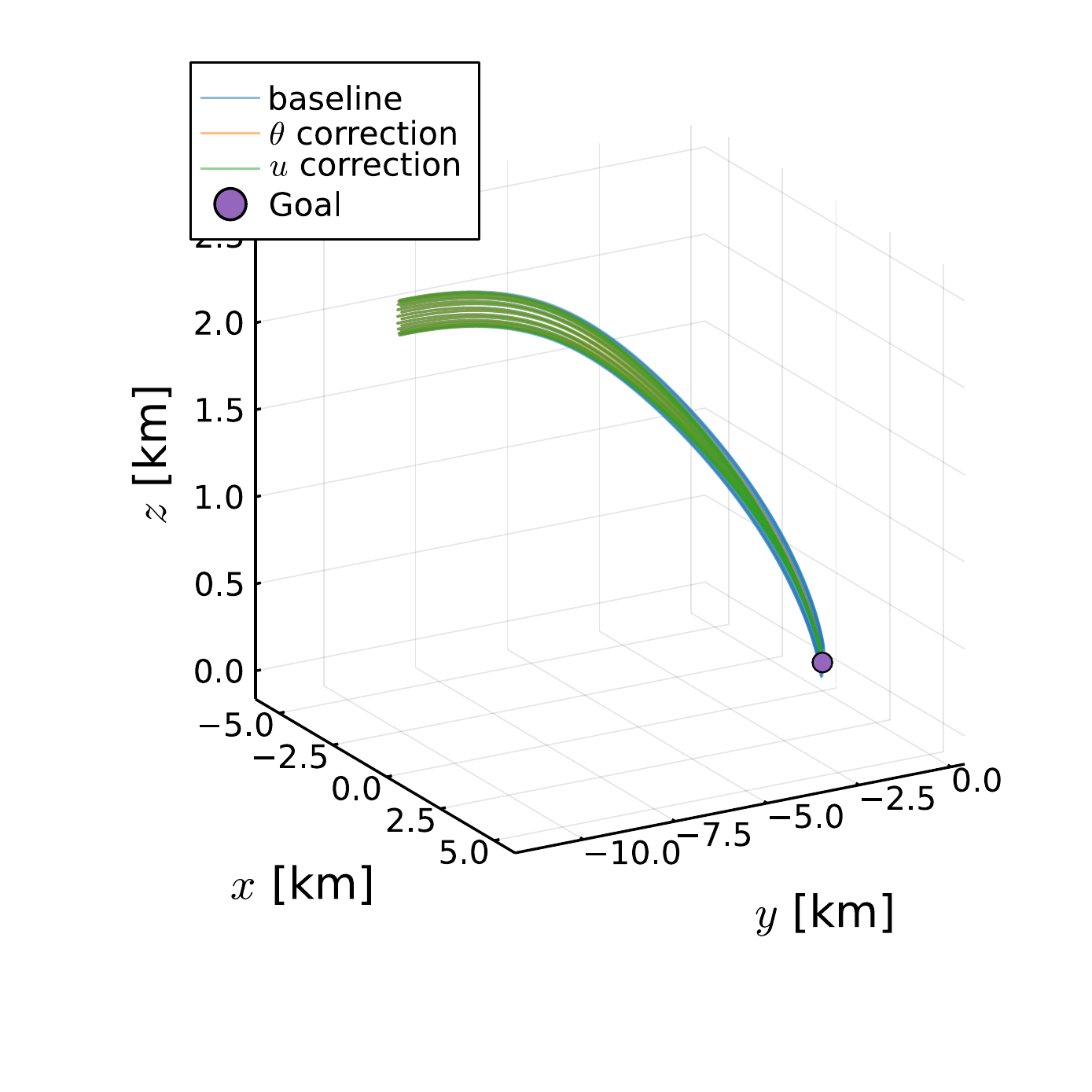}
		\vspace{-1cm}
		\caption{Three-Dimensional Trajectory with Incremental Correction - Single Random Seed (Case 2)} \label{Fig_3D_C2}
	\end{center}
\end{figure}

\begin{figure}[!ht]
	\begin{center}
		\includegraphics[width=0.7\textwidth]{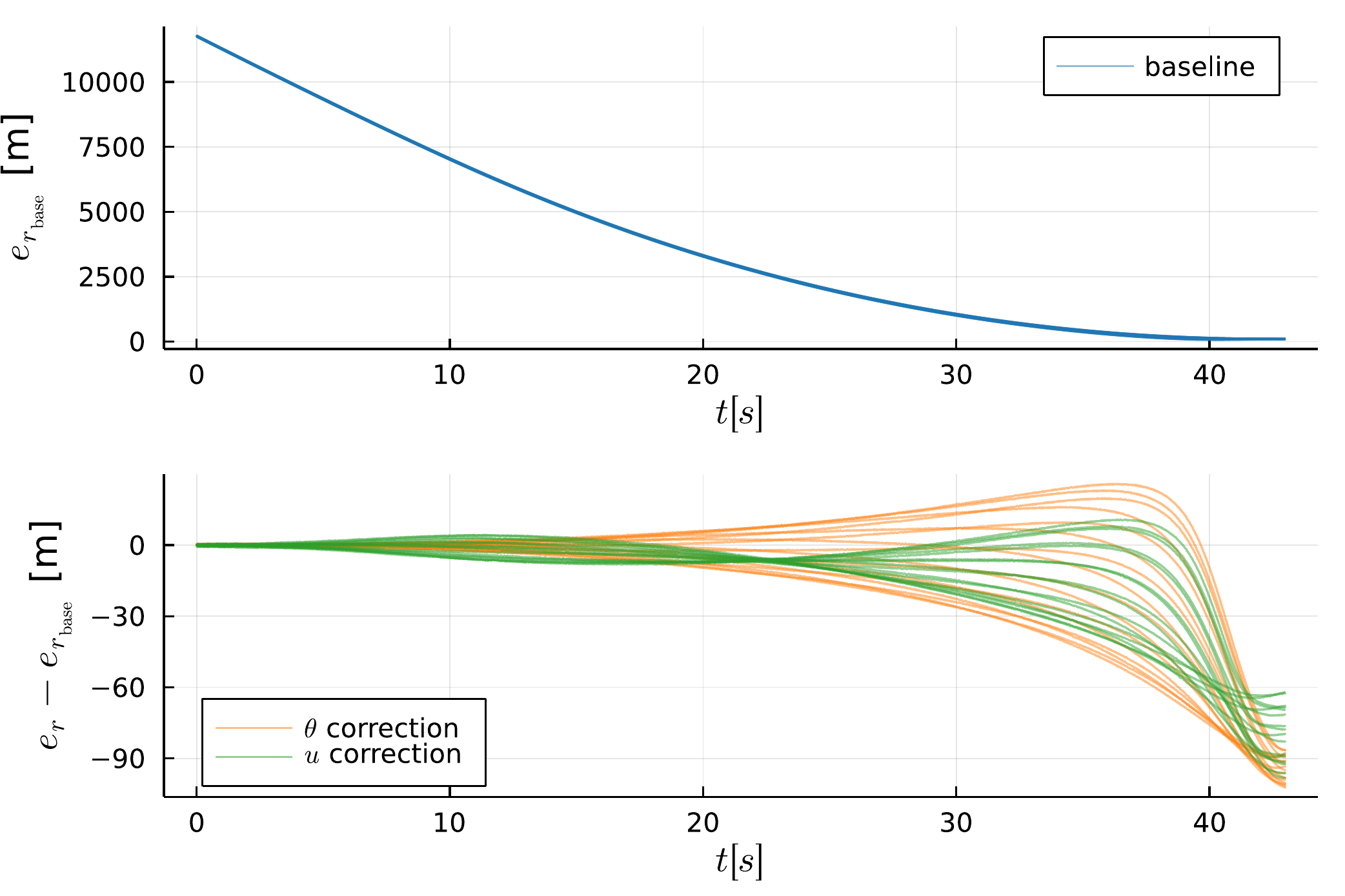}
		\caption{Position Error History with Incremental Correction - Single Random Seed (Case 2)} \label{Fig_e_r_basediff_C2}
	\end{center}
\end{figure}

\begin{figure}[!ht]
	\begin{center}
		\includegraphics[width=0.7\textwidth]{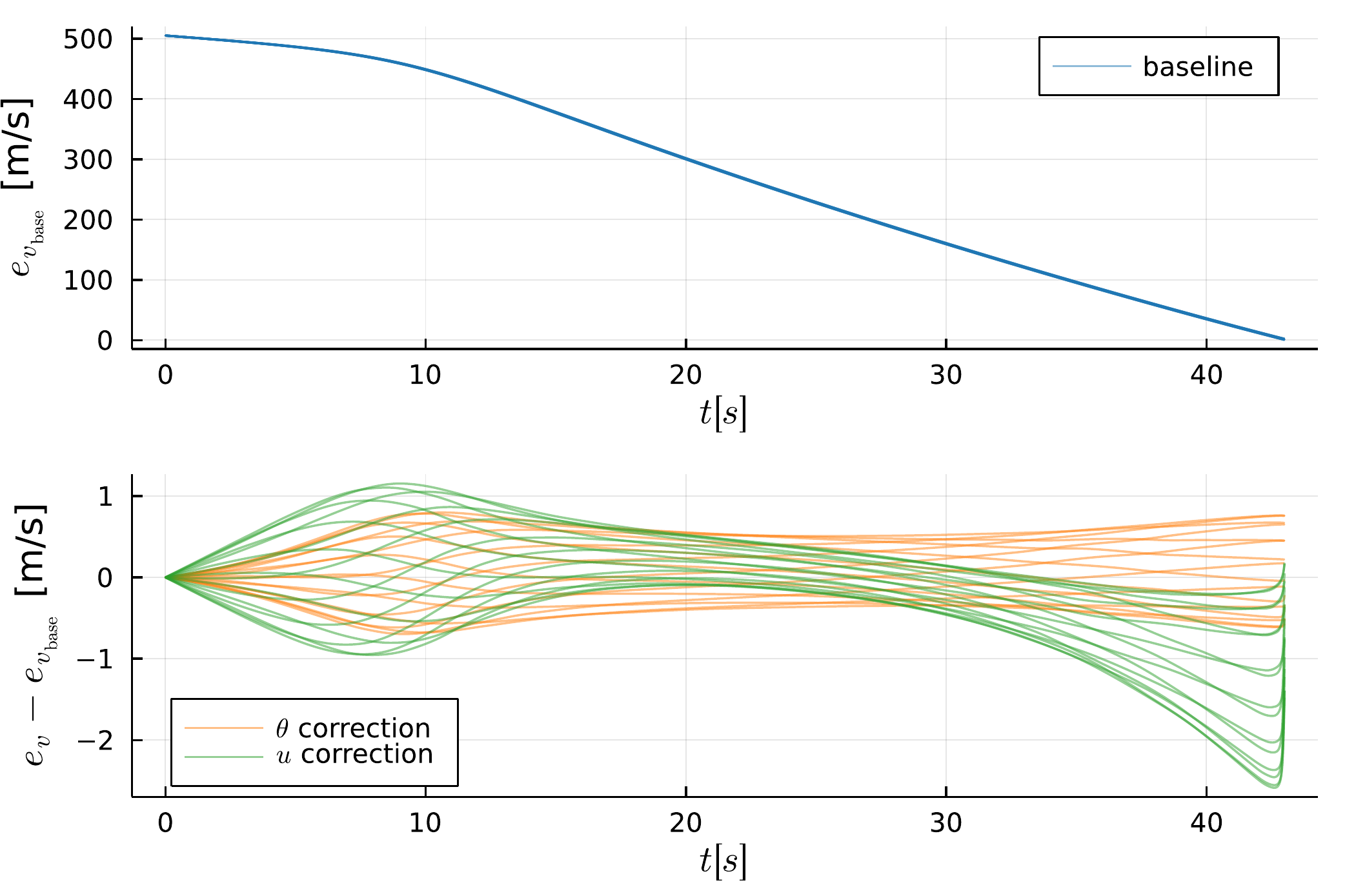}
		\caption{Velocity Error History with Incremental Correction - Single Random Seed (Case 2)} \label{Fig_e_v_basediff_C2}
	\end{center}
\end{figure}

\begin{figure}[!ht]
	\begin{center}
		\includegraphics[width=0.7\textwidth]{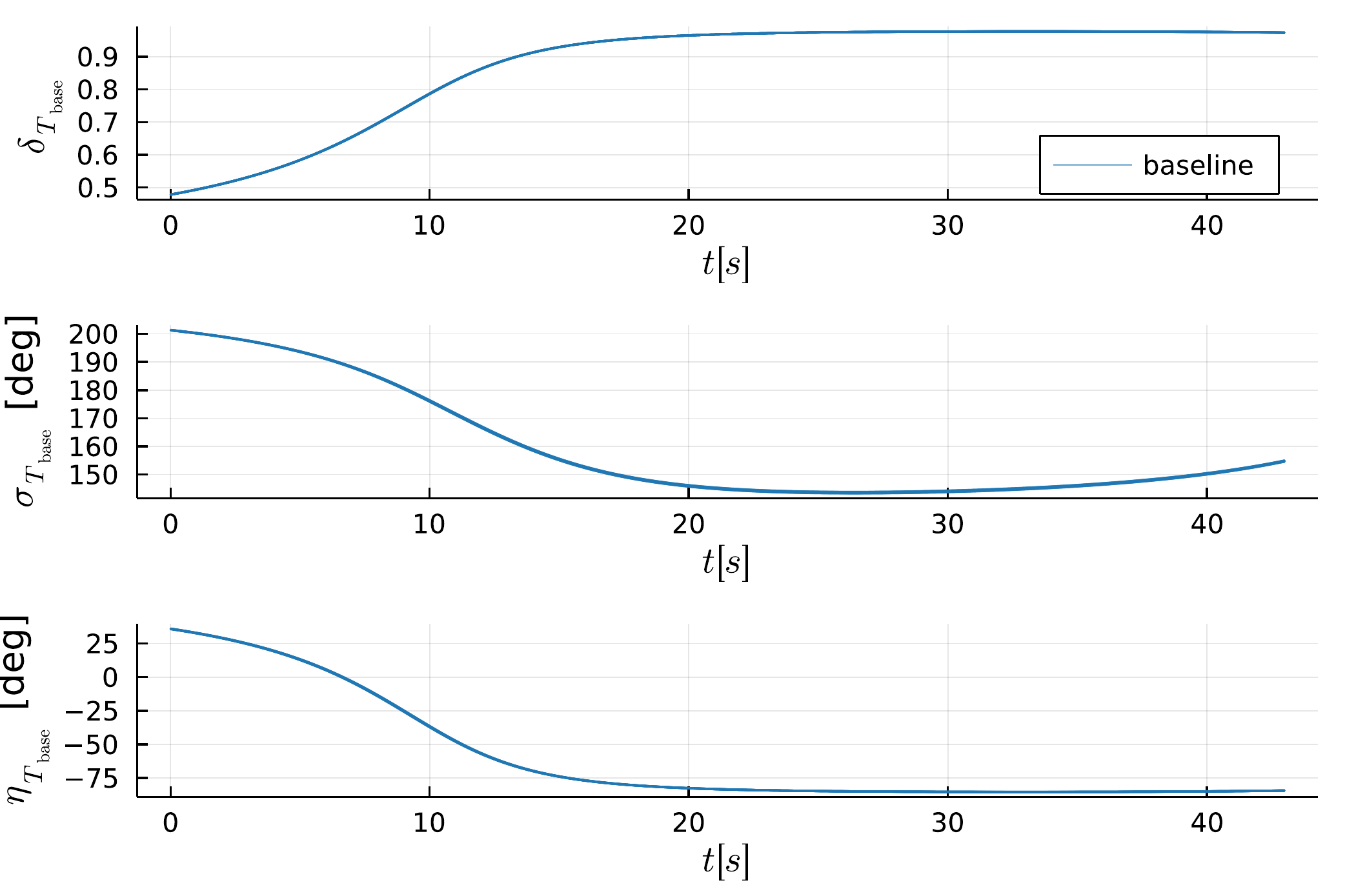}
		\caption{Baseline Control History with Incremental Correction - Single Random Seed (Case 2)} \label{Fig_u_base_C2}
	\end{center}
\end{figure}

\begin{figure}[!ht]
	\begin{center}
		\includegraphics[width=0.7\textwidth]{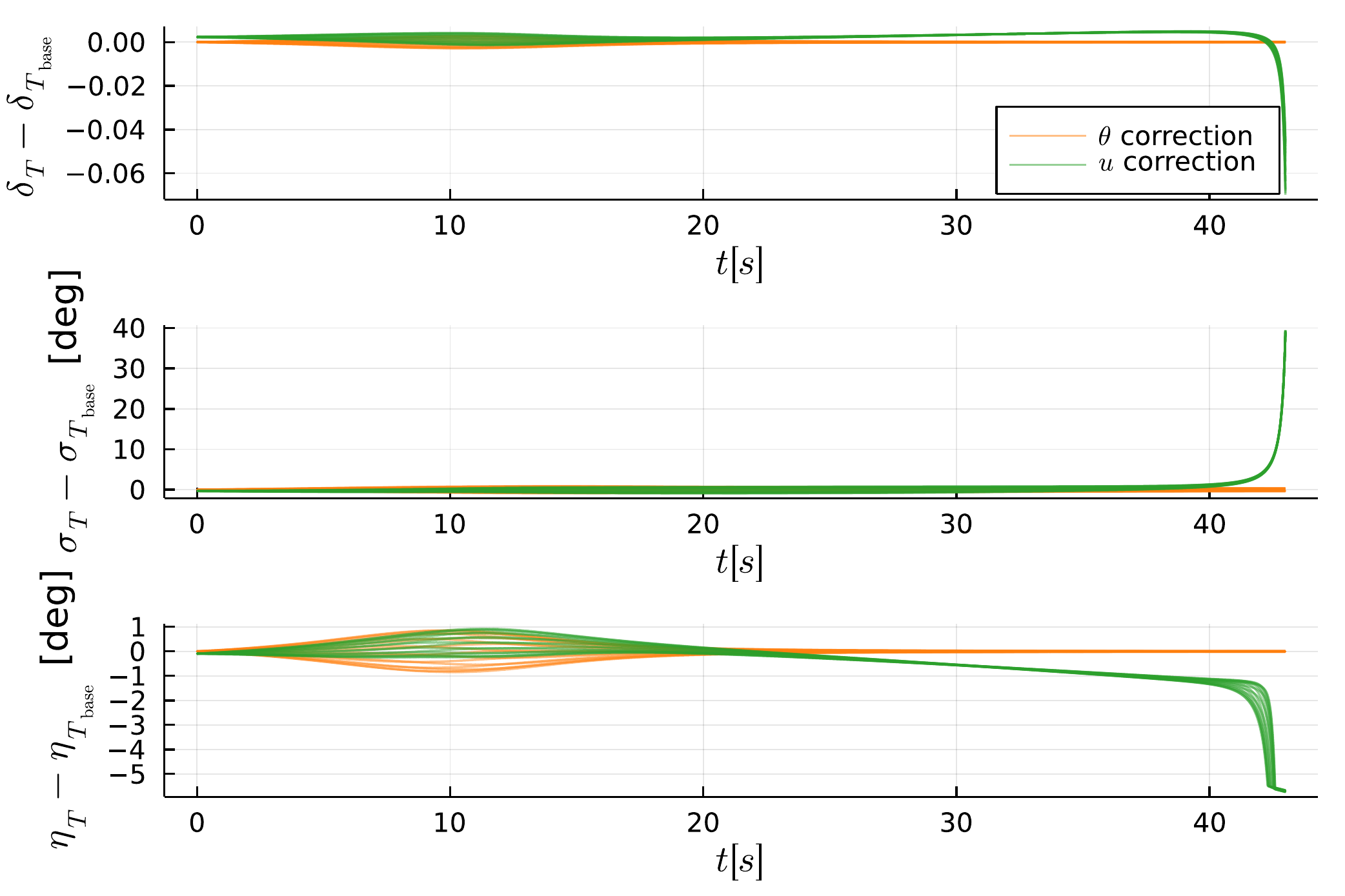}
		\caption{Control Update History with Incremental Correction - Single Random Seed (Case 2)} \label{Fig_u_diff_C2}
	\end{center}
\end{figure}

Figures \ref{Fig_3D_C2}-\ref{Fig_u_diff_C2} show three-dimensional trajectory, position error history $e_{r}\left(t\right) = \left\|\mathbf{r}\left(t\right)-\mathbf{r}_{f_{d}}\right\|$, velocity error history $e_{v}\left(t\right) = \left\|\mathbf{v}\left(t\right)-\mathbf{v}_{f_{d}}\right\|$, baseline control history $\mathbf{u}_{base}\left(t\right)$, and control update history $\tilde{\mathbf{u}}\left(t\right) = \mathbf{u}\left(t\right) - \mathbf{u}_{base}\left(t\right)$, respectively, with different Stage 2 correction strategies for a fixed random seed. For clear comparison of correction methods with respect to the baseline performance, each of the lower panels of Figs. \ref{Fig_e_r_basediff_C2}, \ref{Fig_e_v_basediff_C2}, along with Fig. \ref{Fig_u_diff_C2} show the pointwise difference between each quantity of interest and its baseline counterpart. In the presence of initial position perturbation, the trajectories obtained with each method form a family of solutions that continuously depend on the initial condition, meaning that the actual trajectories depicted in Fig. \ref{Fig_3D_C2} do not involve any disruptive irregularities due to the incremental correction. The convergence of position and velocity errors over time follows the similar general trend at the large scale with or without incremental correction as it can be seen from Figs. \ref{Fig_e_r_basediff_C2} and \ref{Fig_e_v_basediff_C2}. Figure \ref{Fig_u_base_C2} confirms that the baseline policy applied to each perturbed initial condition generated smooth control input histories in a small region around the baseline control history for the nominal initial condition shown in Fig. \ref{Fig_u_C1}. However, with the control function correction method, a rapid change in $\sigma_{T}$ near the final time appears in Fig. \ref{Fig_u_diff_C2} for all initial conditions, with a similar pattern as observed in Fig. \ref{Fig_u_C1} for Case 1. This undesirable behaviour might be attributed to the discrepancy between the actual trajectory and the baseline trajectory that is predicted at the instance of computing the corrective input and taken as the reference for dynamics linearisation. The discrepancy is usually the largest at the final time as it accumulates over time, leading to the amplification of corrective input $\tilde{\mathbf{u}}\left(t\right)$ as $t \rightarrow t_{f}$. Another reason lies in the essential difference between the two correction methods, which differ in the choice of decision variables. The control function correction method not only enforces constraint satisfaction, but also minimises the weighted $\mathcal{L}_{2}$-norm of the corrective input $\tilde{\mathbf{u}}\left(t\right)$. The parameter correction method also solves for the minimum $l_{2}$-norm solution $\tilde{\btheta}$ for the linear system of equations, but it does not necessarily translate into minimal amount of change at the level of control input $\mathbf{u}$ with respect to the signal norm.

Table \ref{Table:SimTime} gives the average simulation time for each combination of random seed and correction method as measured on a Macbook Pro 15-inch 2017 with 2.8 GHz Quad-Core Intel Core i7 CPU and 16GB 2133 MHZ LPDDR3 RAM. An ensemble ODE problem was constructed for each random seed and correction method to simulate the same dynamics for 16 different initial conditions through multithreading. Then, the average simulation time is obtained by repeating the ensemble ODE simulation multiple times and then taking the average of elapsed times to reduce noise. In Table \ref{Table:SimTime}, the simulation time for the baseline case does not include the time consumed in baseline policy training. The computational load in terms of average simulation time is shown to be heavier for the parameter correction than the control function correction, at least in the current implementation. The point that mainly contributes to the increase in the simulation time for both correction methods is thought to be the Jacobian calculation using automatic differentiation (See Remarks \ref{Rem:Jacobian_calculation} and \ref{Rem:AD_solver}). In the parameter correction method, automatic differentiation is called only once when the correction is triggered to compute the sensitivity matrix of a large dimension and also its pseudoinverse. On the other hand, in the control function correction method, automatic differentiation is called at two different points; i) once when the correction is triggered to obtain linearised system matrices $\mathbf{A}_{u}\left(t\right)$ and $\mathbf{B}_{u}\left(t\right)$ for $\forall t \in \left[t_{0}, t_{f}\right]$ around the \emph{stored} baseline state and input that are required to precompute the matrix $\bar{\bPsi}$ in the update equation given by Eq. \eqref{Eq:u_opt}, and ii) at each instance to linearise the system dynamics around the \emph{measured} current state and the \emph{stored} baseline input to obtain $\mathbf{B}_{u}\left(t\right)$ in Eq. \eqref{Eq:u_opt}. Note that the average simulation time shown in Table \ref{Table:SimTime} incorporates the times due to multithreading initialisation and data transfer, memory allocation, correction computation done at the beginning of simulation, simulation execution, and garbage collection, except the compilation time spent at the first run because of the Just-In-Time compilation behaviour of \texttt{Julia}. Therefore, the data provided here should be regarded as an indicator of relative computational burden. More detailed profiling along with code optimisation will be needed to assess the actual performance on a real-time computer.

Overall, the results of Cases 1 and 2 together imply that the proposed two-stage approach with parameter correction can effectively compensate for small perturbation from the condition considered for baseline policy training. Note that the parameter correction method can be seen as the deliberate induction of overfitting of neural network parameters to satisfy the performance output constraints at given time points. Therefore, the improved landing position accuracy/precision with the parameter correction method can be viewed as an advantage of having an overparametrised function as the control law.

\begin{table}[!ht]
	\begin{center}
		\caption{Average Simulation Time in Seconds (Case 2)} \label{Table:SimTime}
		\begin{tabular}{cccc}
			\hline
			Random Seed	&	baseline	&	$\theta$ correction		& $u$ correction\\
			\hline
			0	&	$0.92290187$	&	$22.150549$	&	$16.141281$\\
			1	&	$0.9858486$		&	$21.831982$	&	$16.177732$\\
			2	&	$0.9294312$		&	$21.921776$	&	$16.727678$\\
			3	&	$0.85724115$	&	$21.274307$	&	$16.34164$\\
			4	&	$0.79484373$	&	$21.282663$	&	$16.316748$\\
			5	&	$1.0670316$		&	$21.866957$	&	$16.353233$\\
			6	&	$0.8742441$		&	$20.86112$ 	&	$14.892452$\\
			7	&	$0.8459512$		&	$19.810734$	&	$14.906239$\\
			8	&	$0.83680046$	&	$19.761679$	&	$14.979348$\\
			9	&	$0.8883742$		&	$20.30517$ 	&	$17.351496$\\
			10	&	$0.93385476$	&	$20.065922$	&	$14.895322$\\
			\hline
		\end{tabular}
	\end{center}
\end{table}

\clearpage
\section{Conclusion} \label{Sec:Concls}
This study presents refinement methods for improving the accuracy of satisfying performance output constraints at given time points in continuous-time dynamic systems in which the control policy is provided by a neural network. Provided a baseline neural network policy constructed a priori, the incremental correction can be performed either at the level of neural network parameters or at the level of control input variables, so that the updated controller can enforce point constraints on performance output. The effectiveness of the proposed two-stage approach consisting of baseline policy optimisation followed by incremental correction was illustrated on a Mars landing guidance problem in the powered descent phase. The two types of incremental correction methods exhibited different performance characteristics, demanding a comparative study before deciding which method will be more appropriate for implementation depending on the application. For the Mars landing example addressed in this study, the parameter correction method showed a relative advantage in computational reliability and landing position targeting accuracy. As long as the required sensitivity matrices can be computed accurately and efficiently, the parameter correction method can be useful as a plug-and-play extension for improving the constraint satisfaction accuracy to any baseline neural network policy trained by using off-the-shelf package relying on unconstrained optimisation. 

	

	
	\bibliographystyle{new-aiaa-nhcho}
	\bibliography{NODE_Correction}
\end{document}